\documentclass[10pt]{article}

\usepackage{extarrows}
\usepackage{amsmath,amssymb,amsthm,amscd,mathrsfs}    
\usepackage{latexsym,amsmath,amssymb,graphicx}   
\usepackage[fleqn]{mathtools}
\usepackage[all]{xy} 
\usepackage{multirow}
\usepackage{color}
\usepackage{hyperref}
\usepackage{epsf}
\usepackage{tikz}
\usetikzlibrary{shapes.geometric, arrows}

\usepackage{graphicx}

\makeatletter
\newcommand\opteq[1]{\mathrel{\mathpalette\opt@eq{#1}}}
\newcommand{\opt@eq}[2]{%
  \begingroup
  \sbox\z@{$#1#2$}%
  \sbox\tw@{\resizebox{!}{.5\ht\z@}{$\m@th#1($}}%
  \nonscript\hskip-\wd\tw@
  \mkern1mu
  \raisebox{-.35\ht\z@}[0pt][0pt]{\resizebox{!}{.5\ht\z@}{$\m@th#1($}}%
  \mkern-1mu
  {#2}%
  \mkern-1mu
  \raisebox{-.35\ht\z@}[0pt][0pt]{\resizebox{!}{.5\ht\z@}{$\m@th#1)$}}%
  \mkern1mu
  \nonscript\hskip-\wd\tw@
  \endgroup
}
\makeatother

\xyoption{v2} \xyoption{2cell} \xyoption{dvips}
\CompileMatrices \numberwithin{equation}{section}
\newtheorem{prop}{Proposition}[section]
\newtheorem{theo}[prop]{Theorem}
\newtheorem{lemm}[prop]{Lemma}
\newtheorem{coro}[prop]{Corollary}
\newtheorem{rema}[prop]{Remark}

\newtheorem{defi}[prop]{Definition}
\newtheorem{conj}[prop]{Conjecture}

\numberwithin{equation}{section}
\makeatletter
\newcommand{\subalign}[1]{%
  \vcenter{%
    \Let@ \restore@math@cr \default@tag
    \baselineskip\fontdimen10 \scriptfont\tw@
    \advance\baselineskip\fontdimen12 \scriptfont\tw@
    \lineskip\thr@@\fontdimen8 \scriptfont\thr@@
    \lineskiplimit\lineskip
    \ialign{\hfil$\m@th\scriptstyle##$&$\m@th\scriptstyle{}##$\crcr
      #1\crcr
    }%
  }
}
\makeatother

\newcommand{\be}{\begin{equation}}
\newcommand{\ee}{\end{equation}}
\newcommand{\IP}{\mathbb{P}}

\newcommand\IZ{\mathbb {Z}}

\newcommand{\IC}{\mathbb{C}}

\newcommand{\CU}{{\mathcal U}}
\newcommand{\ba}{\begin{array}}
\newcommand{\ea}{\end{array}}

\newcommand{\CV}{{\mathcal V}}

\newcommand{\CS}{{\mathcal S}}
\newcommand{\CB}{{\mathcal B}}

\newcommand{\CN}{{\mathcal N}}

\newcommand{\CW}{{\mathcal W}}

\newcommand{\bal}{\begin{aligned}}
\newcommand{\eal}{\end{aligned}}

\newcommand{\longto}{\longrightarrow}

\newcommand{\ch}{{\mathrm{ch}}}

\newcommand{\CO}{{\mathcal O}}

\newcommand{\CA}{{\mathcal A}}
\newcommand{\CH}{{\mathcal H}}

\newcommand{\CM}{{\mathcal M}}
\newcommand{\CL}{{\mathcal L}}

\newcommand{\CR}{{\mathcal R}}
\newcommand{\CC}{{\mathcal C}}
\newcommand{\CI}{{\mathcal I}}

\newcommand{\CQ}{{\mathcal Q}}

\newcommand{\calP}{{\mathcal P}}

\newcommand{\IA}{{\mathbb A}}

\CompileMatrices \LaTeXdiagrams \UseAllTwocells
\newdimen\tableauside\tableauside=1.0ex
\newdimen\tableaurule\tableaurule=0.4pt
\newdimen\tableaustep
\def\phantomhrule#1{\hbox{\vbox to0pt{\hrule height\tableaurule width#1\vss}}}
\def\phantomvrule#1{\vbox{\hbox to0pt{\vrule width\tableaurule height#1\hss}}}
\def\sqr{\vbox{%
  \phantomhrule\tableaustep
  \hbox{\phantomvrule\tableaustep\kern\tableaustep\phantomvrule\tableaustep}%
  \hbox{\vbox{\phantomhrule\tableauside}\kern-\tableaurule}}}
\def\squares#1{\hbox{\count0=#1\noindent\loop\sqr
  \advance\count0 by-1 \ifnum\count0>0\repeat}}
\def\tableau#1{\vcenter{\offinterlineskip
  \tableaustep=\tableauside\advance\tableaustep by-\tableaurule
  \kern\normallineskip\hbox
    {\kern\normallineskip\vbox
      {\gettableau#1 0 }%
     \kern\normallineskip\kern\tableaurule}%
  \kern\normallineskip\kern\tableaurule}}
\def\gettableau#1 {\ifnum#1=0\let\next=\null\else
  \squares{#1}\let\next=\gettableau\fi\next}
\allowdisplaybreaks 

\DeclareMathAlphabet{\mathpzc}{OT1}{pzc}{m}{it}
\begin{document}

\title{Hilbert schemes of nonreduced divisors in Calabi-Yau threefolds and $W$-algebras} 
\author{Wu-yen Chuang, Thomas Creutzig, D.-E. Diaconescu, Yan Soibelman}
\date{}  
\maketitle
 \begin{abstract} 
A $W$-algebra action is constructed via Hecke transformations on the equivariant Borel-Moore homology of the Hilbert scheme of points on a 
nonreduced plane in three dimensional affine space. The resulting $W$-module  is then identified to the vacuum module. The construction is based on a generalization of the ADHM construction as well as the $W$-action on the 
equivariant Borel-Moore homology of the moduli space of instantons constructed by Schiffmann and  Vasserot. 
\end{abstract}

\tableofcontents

\section{Introduction}

The AGT conjecture \cite{AGT} leads to conjectural geometric constructions of $W$-algebra  actions on the equivariant Borel-Moore homology of moduli spaces of sheaves on surfaces. Independent proofs of this conjecture for moduli spaces of rank $r\geq 1$ framed sheaves on the complex projective plane have been given by Schiffmann and Vasserot \cite{Cherednik_W} and 
Maulik, Okounkov 
\cite{Quantum_Y}. In both cases one first constructs an action of the affine Yangian of ${\mathfrak gl}_1$ on the equivariant Borel-Moore homology of the corresponding moduli space, while the $W$-algebra action is obtained through a free field realization. 
A generalization to moduli spaces of framed $G$-instantons on $\IC^2$ for more general reductive group $G$ using Donaldson-Uhlenbeck compactifications was proven by Braverman, Finkelberg and Nakajima in \cite{Inst_W_alg}. 

An alternative geometric construction of the $W$-algebra action was 
carried out by Negut 
\cite{AGT_Ext} using the shuffle algebra algebra realization of the affine Yangian. Using the Ext operators of \cite{Extop,Fivedvertex},
this yields a proof of the AGT conjecture for rank two quiver gauge theories on
$\mathbb{C}^2$ as well as a proof \cite{AGT_shuffle} of the 
five-dimensional analog of the AGT conjecture for any quiver gauge theory with
gauge group ${SU}(r)$, $r\geq 1$. 
The latter involves a $q$-deformed $W$-algebra action on the equivariant $K$-theory of moduli spaces. This construction has been further generalized to 
moduli spaces of stable  sheaves on smooth projective surfaces 
in  \cite{AGT_surfaces}. 

Finally, motivated by work of Gaiotto and Rapcak
\cite{Vertex_corner}, Nekrasov \cite{BPS_CFT} and Nekrasov and Prabhakar \cite{Spiked_Inst}, an action of a more general class of vertex algebras on 
the dual of the compactly supported equivariant vanishing cycle cohomology of certain quiver moduli spaces was constructed 
by Rapcak, Soibelman, Yang and Zhao in \cite{Vertex_COHA}. 

A central element in the constructions of \cite{Cherednik_W,Vertex_COHA} 
is the action of certain cohomological Hall algebras on the cohomology of moduli spaces. 
Cohomological Hall algebra ($3d$ COHA for short) 
was introduced by  Kontsevich and Soibelman \cite{COHA_DT} for categories of modules associated to Quillen-smooth algebras with potential e.g. path algebras of quivers with potential. Those are the hearts of $t$-structures of Calabi-Yau categories of dimension three.   An independent construction of a certain cohomological Hall algebra ($2d$ COHA for short) for categories of modules over preprojective algebras, which are hearts of $t$-structures of two dimensional Calabi-Yau categories,  was developed by 
Schiffmann and Vasserot \cite{Cherednik_W,COHA_quivers_I,COHA_quivers_II}, and also Yang and Zhao \cite{COHA_preproj}. A comparison between the two approaches was carried out in \cite{COHA_2d} where it was shown that the $2d$ COHA can be obtained by dimensional reduction of a special case of $3d$ COHA.  Further developments on $3d$ COHAs and their representations include \cite{Critical_COHA, COHA_enveloping,Remarks_COHA} while 
geometric constructions of $2d$ COHAs for various categories of coherent sheaves on surfaces have been developed in 
\cite{COHA_Higgs, COHA_triples,COHA_surface, Categ_COHA}. 
It should be also noted that the proof of the AGT conjecture given in \cite{Cherednik_W} uses the two dimensional variant, while the generalization proven in \cite{Vertex_COHA} employs three-dimensional cohomological Hall algebras. 

From this perspective, the geometric framework of the present paper consists of torsion sheaves on a 
Calabi-Yau threefold with set theoretic support on a 
given divisor. Such a generalization of the AGT framework was first proposed 
in \cite{Vertex_COHA}. 

In more detail, let $D_r\subset \IA^3:={\rm Spec}\, \IC[x_1, x_2, x_3]$ be the divisor $x_3^r=0$. 
Let also ${Hilb}(r,n)$ be the Hilbert scheme of zero dimensional 
coherent quotients $\CO_{D_r}\twoheadrightarrow Q$ with $\chi(Q)=n$ and let 
\[ 
Z_r(q)=\sum_{n\geq 0} q^n \chi({\mathcal Hilb}(r,n)) 
\]
be the generating function of Euler numbers. Then Section 9.3.2 in \cite{Vertex_corner} shows that $Z_r(q)$ is equal to the vacuum character of the $W$-algebra $W_\kappa({\mathfrak gl}_r)$ up to a prefactor. 
 More precisely, for any level $k$, one has 
\[
\text{ch}[W_\kappa({\mathfrak gl}_r)] = q^{-\frac{c_{k, r}}{24}} Z_r(q) 
\]
with
\[
c_{k, r} = (r-1) - r(r^2-1)\left(\sqrt{k+r} - \frac{1}{\sqrt{k+r}} \right)^2.   
\]
Motivated by this observation, the main goal of the present work is to provide an explicit geometric construction of a $W$-algebra action on the localized equivariant Borel-Moore homology 
of the Hilbert scheme 
\[
{Hilb}(D_r) = \amalg_{n\geq 0}{Hilb}({r,n})
\] 
and prove that the resulting $W$-module is isomorphic to the vacuum module.
The equivariant structure is induced by the torus action
$(\IC^\times)^{\times 2}\times \IA^3 \to \IA^3$, 
\[
(t_1, t_2)\times (x_1,x_2,x_3) \mapsto (t_1x_1, t_2x_2, t_1^{-1}t_2^{-1}x_3)
~.
\]
Let ${\bf T}_0=\IC^\times \times \IC^\times$, let $R_0\simeq \IC[x,y]$ be the cohomology ring of the classifying space $B{\bf T}_0$ and let 
$K_0\simeq \IC(x,y)$ be its fraction field. For further reference let 
\[
{\bf V}_{K_0}^{(r)} = \bigoplus_{n\geq 0} H^{{\bf T}_0}({Hilb}(r,n)) 
\otimes_{R_0} K_0
\]
denote the localized equivariant Borel-Moore homology of the Hilbert scheme.
As in \cite{Cherednik_W,Quantum_Y} one first constructs a degenerate DAHA, or, equivalently, Yangian action via Hecke transformations. 

\subsection{The degenerate DAHA action}\label{degDAHAaction}

The degenerate DAHA is an
infinite dimensional associative algebra 
${\bf SH}^c$ constructed by Schiffmann and Vasserot \cite{Cherednik_W} 
as a degeneration of Cherednik's $GL_n$ double affine Hecke algebra. 
As briefly reviewed in Section \ref{SHsect}, this algebra is generated by the elements 
 $D_{1,l}, D_{0, l}, D_{-1,l}$, $l\geq 0$, 
 over a polynomial ring $\IC(\kappa)[{\bf c}_0, {\bf c}_1, \ldots]$ where 
 $\kappa$ is a formal parameter and ${\bf c}_l$, $l\geq 0$ are central elements. The defining relations are written explicitly in equations 
\eqref{eq:SHA} -- \eqref{eq:SHE}. As shown in \cite{DAHA_Yangian} it is in fact related by specialization
to the Yangian algebra of affine ${\mathfrak gl}_1$. A specialization 
$\IC(\kappa)[{\bf c}_0, {\bf c}_1, \ldots]\to K$ of 
this algebra has been shown in \cite[Thm. 3.2]{Cherednik_W} to act on the localized equivariant Borel-Moore homology 
\[
{\bf L}_{K_{r+2}}^{(r)} = \bigoplus_{n\geq 0} H^{{\bf T}_{r+2}}(\CM(r,n))_{K_{r+2}}
\]
of the 
moduli space of rank $r$ framed torsion free sheaves on the projective plane. Here ${\bf T}_{r+2} = 
(\IC^{\times})^{r+2}$ denotes the natural torus which acts on these moduli spaces and $K_{r+2}$ is the field of fractions of the cohomology ring of the classifying space, 
$H(B{\bf T}_{r+2})$.  More details 
are provided in Sections \ref{gentorus} and \ref{SHmodsect}. 

In the Hilbert scheme context the direct construction of Hecke transformations encounters significant technical problems. 
The natural Hecke correspondence in this case is the 
nested Hilbert scheme ${Hilb(r,n,n+1)}$ parametrizing flags of ideal sheaves 
$\CI_1\subset \CI_2$ on $D_r$ with $\chi(\CO_{D_r}/\CI_2)=n$ and 
$\chi(\CO_{D_r}/\CI_1) = n+1$. By construction there are natural projections 
$\rho_1 : {Hilb(r,n,n+1)} \to Hilb(r,n)$ and $\rho_2: {Hilb(r,n,n+1)}\to Hilb(r, n+1)$ and $\rho_2$, and $\rho_2$ can be shown to be proper. The main technical difficulty resides in the fact that $\rho_1$ is not a locally complete intersection morphism, hence one cannot construct a refined Gysin pullback 
$\rho_1^!$. This problem has been encountered in a similar context in
\cite{AGT_surfaces,COHA_quivers_I,COHA_surface,Categ_COHA,KHall} where it was solved using various techniques. The construction of \cite{COHA_quivers_I}
uses a factorization of $\rho_1$ into simpler morphisms which is specific to 
moduli spaces of Nakajima quiver varieties. The constructions of \cite{AGT_surfaces,COHA_surface,Categ_COHA,KHall}
use virtual pullbacks,  as constructed in \cite{Virt_pullbacks} or some derived variant. A common feature of all these cases is that their constructions take place in an abelian category of homological dimension two, which is an essential condition.
Such an approach does not apply to the present case since the relative obstruction theory of $\rho_1$ is not perfect of amplitude $[-1,\ 0]$. 
This reflects the fact that the present construction takes place in a 
category of homological dimension three rather than two as in loc. cit. 

In order to circumvent this obstacle, the strategy used in this paper employs an equivariant embedding of the Hilbert scheme
$Hilb(r,n)$ in a smooth ambient space which yields an injection for localized Borel-Moore homology. 
This is obtained in Propositions \ref{frquivpropA} and \ref{hilbhiggspropA} which construct a closed embedding 
$f: Hilb(r,n) \hookrightarrow \CM(r,n)$, in the
 moduli space of rank $r$ framed sheaves on the projective plane with second Chern class $n$. Moreover, as shown in Section \ref{CYtorus}, there is an injective 
group homomorphism ${\bf T}_0\hookrightarrow {\bf T}_{r+2}$, 
\be\label{eq:injtorus}
(t_1, t_2) \mapsto (t_1, t_2, t_1^{a-1}t_2^{a-1}), \qquad 
1\leq a\leq r
\ee
such that $f$ is equivariant. 
Then an important technical result in the current approach states that any connected component of the ${\bf T}_0$-fixed locus 
in $\CM(r,n)$ which intersects $Hilb(r,n)$ nontrivially must be:
\begin{itemize} 
\item[$(a)$] contained in $Hilb(r,n)$, and
\item[$(b)$] zero dimensional i.e. a 
single closed point.  
\end{itemize}
This is proven in Section \ref{CYtorus}, Corollary\ref{CYfixedcorB} 
and Lemma \ref{CYfixedlemmD}.  
Then
the pushforward map  for localized equivariant Borel-Moore homology is injective and yields an identification 
\be\label{eq:Hsub}
H^{{\bf T}_0}(Hilb(r,n))_{K_0}\simeq \bigoplus_{\alpha \in Hilb(r,n)^{{\bf T}_0}} K_0[\alpha] \subset H^{{\bf T}_0}(\CM(r,n))_{K_0}~.
\ee
For future reference let 
\[
{\bf L}^{(r)}_{K_0} = \bigoplus_{n\geq 0}  H^{{\bf T}_0}(\CM(r,n))_{K_0}~.
\]
Using these results,  Hecke 
transformations for the Hilbert scheme are constructed in Section \ref{Heckehilb} by analogy to \cite[Section 3.2]{Cherednik_W}. 
Let  $\pi_i: \CA(r,n)\times \CA(r,n+1) \to \CA(r,n+i-1)$, $1\leq i \leq 2$ 
denote the canonical  projections and let $\tau_{n,n+1}$ 
denote the
 universal line bundle on the correspondence variety $\CA(r,n,n+1)$. Let 
$\gamma:\CA(r,n,n+1) \to  \CA(r,n)\times \CA(r,n+1)$ denote the canonical closed embedding. Then Lemmas \ref{HeckelemmB} and \ref{HeckelemmB} show that the transformations 
\[
h^+_l(x) = \pi_{2*}(\gamma_*(e_{{\bf T}_0} (\tau_{n,n+1})^l) \cdot \pi_1^*(x)), \qquad 
h^-_l(x) = \pi_{1*}(\gamma_*(e_{{\bf T}_0} (\tau_{n-1,n})^l) \cdot \pi_2^*(x)) 
\]
are well defined for any $l\geq 0$ and factor through 
${\bf V}^{(r)}_{K_0}\subset {\bf L}^{(r)}_{K_0} $. Therefore, by a slight abuse of notation, one obtains linear transformations $h^{\pm}_l\in 
{\rm End}_{K_0}({\bf V}^{(r)}_{K_0})$. In addition, one defines the diagonal 
transformations $h^0_l\in {\rm End}_{K_0}({\bf V}^{(r)}_{K_0})$ using the 
plethystic powers of the universal vector bundle $\CV_{r,n}$ on $\CM(r,n)$. 
Namely, denoting by $\CL_k$, $1\leq k\leq n$, the virtual equivariant Chern roots of $\CV_{r,n}$, let 
\[ 
h_{l}^0 = \sum_{k=1}^n e_{{\bf T}_0}(\CL_k)^l \qquad l\geq 0.  
\]
Finally, in order to state the first main result of the present paper 
let ${\bf SH}^{(r)}_{K_0}$ be the specialization of the degenerate DAHA induced by 
the injective group homomorphism ${\bf T}_0\subset {\bf T}_{r+2}$. By construction, this is an algebra over   $K_0=\IC(x,y)$. Then one has: 

\begin{theo}\label{mainthmA} 
The map $\{D_{-1,l},\, D_{0,l},\, D_{1,l}\,|\, l \in \IZ,\, l\geq 0\}\to 
{\rm End}_{K_0}({\bf V}^{(r)}_{K_0})$ 
given by 
\[ 
D_{1,l}\mapsto x^{1-l}y h^+_l, \qquad D_{0,l} \mapsto x^{1-l}h_{l}^0, 
\qquad D_{-1,l} \mapsto (-1)^{r-1}x^{-l}h^{-}_l, \qquad l\geq 0 
\]
extends uniquely to an algebra homomorphism 
\[
{\bf SH}^{(r)}_{K_0} \to {\rm End}_{K_0}({\bf V}^{(r)}_{K_0}).
\] 
\end{theo}

Theorem \ref{mainthmA} is proven in Proposition \ref{SHaction} using Lemmas \ref{correspfixedlemmAB}, \ref{HeckelemmB}, \ref{HeckelemmC} and \ref{CYtruncationA}. The 
main strategy is to derive the above result from 
\cite[Thm. 3.2]{Cherednik_W} by a detailed equivariant fixed point analysis.

\subsection{The $W$-algebra action}\label{Walgact}
The next goal is to convert the ${\bf SH}^{(r)}_{K_0}$-action in Theorem 
\ref{mainthmA} into a $W$-algebra action. 
Let $W_\kappa({\mathfrak gl}_r)$ be the $W$-algebra for ${\mathfrak gl}_r$ at level $\kappa= -x^{-1}y$ considered as an algebra over the ground field $F=\IC(\kappa)$. A rigorous mathematical theory of $W$-algebras and their representations has been developed in \cite{Rep_th_W}. A brief introduction is provided for completeness in Section \ref{Wsect}. In particular there is a canonical vacuum module $\pi_0$ which admits a free field realization. 

In the present context the relation between the degenerate DAHA and the $W$-algebra is obtained by a careful specialization of 
\cite[Thm. 8.21]{Cherednik_W}, including the main steps in the proof.  This is carried out in detail in Appendix \ref{DDAHAW}. In particular,  Lemma \ref{CYSHW} shows that there is a surjective homomorphism of algebras 
\be\label{eq:SHtoWmapA} 
{\overline \Theta}_0^{(r)} : {\mathfrak U}({\bf SH}_{K_0}^{(r)}) 
\to \CU_0(W_\kappa({\mathfrak gl}_r)) 
\ee
where the domain is the current algebra of the degenerate DAHA and 
the target is the image of the current algebra ${\mathfrak U}(W_\kappa({\mathfrak gl}_r))$ in ${\rm End}(\pi_0)$. Corollary \ref{modequiv} states that the restriction 
\be\label{eq:SHtoWmapB} 
{ \Theta}_0^{(r)} : {\bf SH}_{K_0}^{(r)}
\to \CU_0(W_\kappa({\mathfrak gl}_r)) 
\ee
yields a categorical equivalence of admissible modules. 

As  in \cite[Definition 8.10]{Cherednik_W}, admissible modules are $\IZ$-graded modules with respect to a natural $\IZ$-grading on the degenerate DAHA
whose graded summands are trivial for sufficiently high degree.  
Moreover, by construction, the underlying $K_0$-vector space of any admissible 
${\bf SH}_{K_0}^{(r)}$-module is canonically identified with the underlying $K_0$-vector space of the corresponding $\CU_0(W_\kappa({\mathfrak gl}_r))$-module. This is briefly reviewed in Section \ref{degcompl}. 

In conclusion, one obtains a representation 
\be\label{eq:HilbWrep}
\pi_0^{(r)} : \CU_0(W_\kappa({\mathfrak gl}_r)) \to {\rm End}_{K_0} 
({\bf V}^{(r)}_{K_0})~.
\ee
Then second main result of this paper is: 
\begin{theo}\label{mainthmB} 
For any $r\geq 1$ the representation $\pi_0^{(r)}$ is isomorphic to the vacuum representation of the $W$-algebra.
\end{theo}

The proof of Theorem \ref{mainthmB} is analogous to the proof of \cite[Thm. 8.21]{Cherednik_W}. The required structure results are proven in 
the present context in Lemmas \ref{SHvaclemm} and \ref{SHmodA}. 

\subsection{Further remarks and open directions} 
In order to conclude this section, it may be helpful to add a few comments on the relation of the present work to \cite{Vertex_COHA}, as well 
as mention a few possible open directions. 

One of the main results of \cite{Vertex_COHA}, is a construction of a certain vertex algebra action on the dual to equivariant vanishing cycle cohomology of the moduli spaces of spiked instantons constructed in \cite{BPS_CFT,Spiked_Inst}. Spiked instantons admit a presentation 
in terms of framed quiver representations of a triply framed 
generalization of the ADHM quiver labelled by framing vectors $(r_1,r_2,r_3)\in (\IZ_{\geq 0})^{\times 3}$. It was shown in \cite{Vertex_COHA} that the dual of the compactly supported equivariant vanishing cycle cohomology of
moduli spaces of stable framed  representations with 
fixed $(r_1,r_2,r_3)$ carries an action of the vertex algebra $V_{r_1,r_2,r_3}$ introduced in \cite{Vertex_corner}. Moreover, the resulting module is identified with a Verma module with the highest weight depending on the equivariant parameters. In particular for $(r_1,r_2,r_3) = (r,0,0)$ the moduli space reduces to the standard moduli space of stable ADHM data, and the action reduces to the $W$-action constructed in \cite{Cherednik_W,Quantum_Y}. 

The present paper generalizes the results of \cite{Cherednik_W,Quantum_Y} in a different direction, using framed quiver representations associated to geometric objects as explained above.
In particular, as shown in Section \ref{frhiggsplane}, spectral correspondence leads to the new framed quiver with potential shown in diagram \eqref{eq:quiver}. The Hilbert scheme ${\mathcal Hilb}_n(D_r)$  is then isomorphic to a closed subscheme $\CQ_0(r,n)$ of the moduli space $\CQ(r,n)$ of framed quiver representations constructed in Proposition \ref{frquivpropC}. 
Moreover, this closed embedding yields an isomorphism in localized equivariant 
Borel-Moore homology and Theorem \ref{mainthmB} identifies the latter with a vacuum $W$-module.

From the point of view of the underlying Calabi-Yau geometry, the moduli 
space $\CQ(r,n)$ is in fact more natural than the Hilbert scheme, since it admits a 
global presentation as the critical locus of a polynomial potential. As such, it is endowed with an equivariant sheaf of vanishing cycles. Then one is naturally led to the following conjecture: 
\begin{conj}\label{conj} 
Let $H_{{\bf T}_0}^{\sf van}(\CQ(r,n))$ denote the dual to the compactly supported equivariant vanishing cycle cohomology of the framed quiver moduli space $\CQ(r,n)$. Let $H_{{\bf T}_0}^{\sf van}(\CQ(r,n))_{K_0}$ denote its localization at $(0)$. Then there is an explicit 
${\bf SH}_{K_0}^{(r)}$ 
action on $H_{{\bf T}_0}^{\sf van}(\CQ(r,n))_{K_0}$ constructed via Hecke correspondences for vanishing cycle cohomology by analogy to \cite{Vertex_COHA}.
\end{conj} 

The main difficulty in proving Conjecture \ref{conj} resides in the technical difficulties involved in working with sheaves of vanishing cycles. 

As was pointed out in \cite{Vertex_COHA}, a natural question for further study is whether similar algebraic structures 
can be associated to more general configurations of divisors in Calabi-Yau threefolds. It was explicitly conjectured in loc. cit. that this should be possible at least for divisors of the form $x_1^{r_1}x_2^{r_2}x_3^{r_3} =0$ in $\IA^3$. At the current stage, the construction of a quiver with potential associated to such configurations is an open problem.  At the same time the results of \cite{AGT_surfaces} point to another possible generalization associated to compact nonreduced divisors in Calabi-Yau threefolds.

{\it Acknowledgements.}  We would like to thank Davide Gaiotto, Nikita Nekrasov, Andrei Okounkov, Miroslav Rapcak, Francesco Sala, Olivier Schiffmann, Yaping Yang, Gufang Zhao and Yu Zhao for very helpful discussions. 
The work of W.-Y.C. was partly supported by MOST grant  107-2115-M-002-009-MY2. T. C. is supported by NSERC $\#$RES0020460.
The work of 
D.-E.D was partially supported by NSF grants DMS-1501612 and DMS-1802410. 
The work of Y.S. was partially supported by an NSF grant and Munson-Simu award of KSU.

\section{Hilbert schemes and framed Higgs sheaves}\label{hilbhiggs}  

Let $S$ be a smooth complex projective surface, let $L$ 
be a line bundle on $S$ and let $X$ be the total space of the line bundle $L$. Let $\pi:X\to S$ denote the canonical projection and $\zeta \in H^0(X, \pi^*L)$ denote the tautological section. The zero locus of $\zeta$ is the image of the zero section $S \to L$, which will be denoted by $S_1$. More generally, 
for any positive integer $r\geq 1$ let $S_r$ denote the nonreduced divisor $\zeta^r =0$, and let $\iota_r: S_r \hookrightarrow X$ denote the canonical closed embedding into $X$.
Note also that the projection $\pi:X\to S$ yields by restriction a projection map 
$\pi_r: S_r \to S$. Moreover, the zero section $\sigma:S \hookrightarrow X$ factors through a closed embedding $\sigma_r:S \to S_r$. 

As proven in \cite[Proposition 2.2]{TT17a}, the abelian category 
$Coh_c(X)$ of coherent $\CO_X$-modules with compact support 
is equivalent to the category of Higgs sheaves on $S$ with coefficients in $L$. 
A Higgs sheaf on $S$ with coefficients in $L$ is defined as a pair $(E, \Phi)$ where $E$ is a coherent $\CO_S$-module and $\Phi:E \to E\otimes L$ is a morphism of $\CO_S$-modules. Such pairs form naturally an abelian category $Higgs(S,L)$, where the morphisms are defined as morphisms $f: E\to E'$ of $
\CO_S$ modules 
such that $\Phi'\circ f = (f \otimes {\bf 1}_L)\circ \Phi$. 
Then Proposition 2.2 of loc. cit.  
proves that there is an equivalence of abelian categories  
\be\label{eq:cateqA} 
Coh_c(X) {\buildrel \sim \over \longto}  Higgs(S,L)
\ee
where $Coh_c(X)$ denotes the abelian category of coherent $\CO_X$-modules with compact support. This equivalence associates to any such sheaf $F$ its direct image $E=\pi_*F$, while the 
Higgs field $\Phi:E\to E\otimes L$ is the direct image, $\Phi=\pi_*(\zeta_F)$ of the canonical morphism 
$\zeta_F= {\bf 1}_F \otimes \zeta : F \to F \otimes L$. 

Now let $\Delta \subset S$ be a smooth connected effective divisor 
on $S$, and let $\Delta_r = \pi_r^{-1}(\Delta)\subset S_r$ be its inverse image in $S_r$. For any positive integer $r\geq 1$ let $D_r$ be the complement of $\Delta_r$ in $S_r$ and let 
$Hilb_n(D_r)$ be the Hilbert scheme of zero dimensional quotients 
\[ 
\CO_{S_r} \twoheadrightarrow Q
\]
with $\chi(Q)=n$ such that the support of $Q$ is contained in $D_r$.  
For each such quotient let $\CI_{S_r} = {\rm Ker}(\CO_{S_r} \twoheadrightarrow Q)$. Clearly the extension of $\CI_{S_r}$ by zero to $X$ is a pure dimension two torsion 
sheaf with compact support. 
Using the equivalence \eqref{eq:cateqA}, the goal of this section is to prove that the Hilbert scheme $Hilb_n(D_r)$ is isomorphic to a moduli space of framed Higgs sheaf quotients on $S$.

\subsection{Framed sheaves}\label{frsheaves} 
This section summarizes the main results on framed sheaves needed in this paper following \cite{Moduli_framed}.  Let $F$ be a fixed locally free $\CO_\Delta$-module. A framed torsion free sheaf on $S$ with respect to $(\Delta,F)$ is a pair $(E, \xi)$ where $E$ is a torsion free 
sheaf on $S$, and $\xi: E\otimes \CO_{\Delta} \to F$ is an isomorphism of 
$\CO_\Delta$-modules. For any such a pair $(E,\xi)$, let 
$f_\xi: E \to F$ denote the morphism of 
$\CO_S$-modules obtained by composing $\xi: E\otimes \CO_\Delta \to F$ 
with the natural projection $E\to E\otimes \CO_\Delta$. 
An isomorphism of framed sheaves $(E_1, \xi_1)$, $(E_2, \xi_2)$  is an isomorphism 
of sheaves $\phi:E_1\to E_2$ such that 
\[
f_{\xi_2} \circ \phi =f_{\xi_1}. 
\]
Moreover, given a parameter scheme $T$ over $\IC$, let $F_T$ denote the pull-back of $F$ to $\Delta \times T$. Then a flat family of rank $r\geq 1$ torsion free sheaves on $S$ 
is a coherent $\CO_{S\times T}$-module $E_T$, flat over $T$, and an isomorphism 
$\xi_T: E_T \otimes \CO_{\Delta\times T} {\buildrel \sim \over \longto} F_T$.

Next suppose the following additional conditions are satisfied: 
\begin{itemize} 
\item[(1)] $\Delta$ is nef and $\Delta^2\geq 0$ in the intersection ring of $S$, and 
\item[$(2)$] $F$ is a good framing sheaf as defined in \cite[Definition 2.4]{Moduli_framed}
satisfying in addition the vanishing condition
${\rm Hom}_\Delta(F, F\otimes 
\CO_\Delta(-k\Delta))=0$ for all $k\in \IZ$, $k\geq 1$. 
\end{itemize}
Then Theorem 3.1 of \cite{Moduli_framed} proves that there is a smooth quasi-projective fine moduli space $\CM(r,\beta, n)$ of framed torsion free sheaves with topological invariants 
\[ 
{\rm rk}(E)=r, \qquad c_1(E) =\beta, \qquad c_2(E) = n. 
\]
In particular, under the above conditions, any framed sheaf $(E,\xi)$ has trivial automorphism group. 

\subsection{Framed Higgs sheaves} 
In the above framework, a framed Higgs sheaf will be a triple $(E,\xi,\Phi)$ where $(E,\xi)$ is a framed torsion free sheaf and $\Phi:E\to E\otimes L$ a Higgs field 
satisfying a certain framing 
condition along $\Delta$. The framing condition is naturally inferred from the following special case.  

\begin{lemm}\label{frHiggslemmB}
Let $(O_r,\Lambda_r)$ be the Higgs sheaf associated to the torsion 
$\CO_X$-module $\CO_{S_r}$ via correspondence \eqref{eq:cateqA}. Then there is a direct sum decomposition
\be\label{eq:frHiggsA}
O_r\simeq \bigoplus_{a=0}^{r-1} L^{-a}  
\ee 
identifying the off-diagonal components of the Higgs field to the canonical isomorphisms \be\label{eq:frHiggsC}
\Lambda_{a+1,a} : L^{-a} {\buildrel \sim \over\longto} L^{-(a+1)} \otimes L.
\ee
Moreover, all other components of $\Lambda_r$ are identically zero. 
\end{lemm} 

{\it Proof}. 
By construction, $O_r=\pi_{*}\CO_{S_r}$. Equation \eqref{eq:frHiggsA} will be proven below by induction on $r \geq 1$. The case $r=1$ is obvious. Suppose  equation \eqref{eq:frHiggsA} holds 
for $(O_r, \Lambda_r)$ and note the canonical exact sequence
\be\label{eq:canexseqA} 
0 \to \CO_{rS_1}(-S_1) \to \CO_{(r+1)S_1} {\buildrel f\over \longto} \CO_{S_1} \to 0 
\ee 
Note also
the canonical isomorphism $\CO_X(S_1) \simeq \pi^*L$ determined by the tautological section $\zeta: \CO_X \to \pi^*L$. 
Then applying $\pi_*$ to the exact sequence \eqref{eq:canexseqA}, one obtains the exact sequences of $\CO_S$-modules
\be\label{eq:canexseqB} 
0\to O_{r}\otimes L^{-1} \to O_{r+1} \xlongrightarrow{\pi_*f} \CO_S \to 0~. 
\ee
Let ${\bar X}$ be the total space of the projective bundle 
$\IP_{S}(\CO_S\oplus L)$, which contains $X$ as an open subscheme. 
Let ${\bar \pi}: {\bar X} \to S$ denote the canonical projection. 
Then there is a commutative diagram of morphisms of $\CO_{\bar X}$-modules 
\[
\xymatrix{
& \CO_{\bar X} \ar[dr]^-{p} \ar[dl]_-{q} & \\
\CO_{(r+1)S_1} \ar[rr]^{f} & & \CO_{S_1} \\}
\]
where all arrows are canonical projections. Applying ${\bar\pi}_*$ yields 
a commutative diagram of morphisms of $\CO_S$-modules
\[
\xymatrix{
& \CO_{S} \ar[dr]^-{{\bar \pi}_*p} \ar[dl]_-{{\bar \pi}_*q} & \\
O_{r+1} \ar[rr]^{{\bar \pi}_*f} & & \CO_{S}~,\\}
\]
 where ${\bar \pi}_*p= {\bf 1}_{\CO_S}$. Therefore ${\bar \pi}_*q$ determines a splitting of 
 the exact sequence \eqref{eq:canexseqB}. This proves  the inductive step. 

The Higgs field $\Lambda_r$ is the direct image of the canonical map 
\[ 
\zeta \otimes {\bf 1}_{\CO_{S_r}} : \CO_{S_r} \to \CO_{S_r}(S_1). 
\]
As observed above, the the tautological section $\zeta: \CO_X \to \pi^*L$
determines a canonical isomorphism $\CO_X(S_1) \simeq \pi^*L$. 
Therefore, equation \eqref{eq:frHiggsC} follows immediately.  

\hfill $\Box$

The framing data for Higgs sheaves on $S$ will be specified by the Higgs sheaf $(F_r,\Psi_r)$ obtained by restricting $(O_r, \Lambda_r)$ to $\Delta \subset S$. Therefore, a framed Higgs sheaf of rank $r$ on $S$ will be a framed rank $r$ torsion sheaf $(E, \xi)$ and a Higgs field $\Phi: E\to E\otimes L$ such that the following diagram is commutative 
\be\label{eq:frHiggsD} 
\xymatrix{ 
E \ar[r]^-{\Phi}\ar[d]_-{f_\xi} & E\otimes L \ar[d]^-{f_\xi\otimes q} \\
F_r \ar[r]^-{\Psi_r} &  F_r \otimes L \otimes \CO_\Delta. 
\\}
\ee
For completeness, recall that the morphism $f_\xi$ in the above diagram is the composition of $\xi: E\otimes \CO_\Delta \to F_r$ with the natural projection 
$E\twoheadrightarrow E\otimes \CO_\Delta$. Similarly $q:L \twoheadrightarrow L\otimes \CO_\Delta$ is the natural projection. 
Naturally, an isomorphism of framed Higgs sheaves $(E, \xi, \Phi){\buildrel \sim \over \longto} (E', \xi', \Phi')$ is an isomorphism 
$f: (E, \xi) {\buildrel \sim \over \longto} (E', \xi')$ of framed sheaves which is at the same time an isomorphism 
$f:(E, \Phi) {\buildrel \sim \over \longto} (E', \Phi')$ of Higgs sheaves. The definition of flat families is also natural, hence the details will be omitted. 

The next result will show that the Hilbert scheme $Hilb_n(D_r)$ is isomorphic to a moduli space of  Higgs sheaf quotients. 
Let 
\be\label{eq:exseqXA} 
0\to \CI_{S_r} \to \CO_{S_r} \to Q\to 0
\ee
be an exact sequence of $\CO_{S_r}$ modules, where $Q$ is 
zero dimensional, supported in the complement  $S_r\setminus \Delta_r$. By extension by zero, this yields an exact sequence of $\CO_X$-modules.  Let $(E, \Phi)$ be the associated Higgs sheaf via correspondence \eqref{eq:cateqA}.  Taking the direct image of the exact sequence \eqref{eq:exseqXA} 
via $\pi:X\to S$ yields an exact sequence  
\be\label{eq:higgsseqA}
0 \to (E, \Phi) \to (O_r, \Lambda_r) \to (G,\Upsilon)\to 0
\ee
in the abelian category $Higgs(S,L)$. Here 
$(G, \Upsilon)$ is a zero dimensional Higgs sheaf on $S$ supported in the complement  $S\setminus \Delta$ such that $\chi(G) = \chi(Q)=n$.

For any $r,n \in \IZ$, $r\geq 1$ let $\CQ uot_n(O_r, \Lambda_r)$ be the {\it Quot}-scheme parametrizing quotients 
\be\label{eq:higgsquotA}
 (O_r, \Lambda_r) \twoheadrightarrow (G,\Upsilon)
 \ee
 in the abelian category $\CH(S,L)$, where $(G, \Upsilon)$ is 
a zero dimensional Higgs sheaf on $S$ supported in $S\setminus \Delta$ such that $\chi(G) = \chi(Q)=n$. Then the equivalence \eqref{eq:cateqA} yields 
\begin{prop}\label{hilbhiggspropA}
There is an isomorphism $Hilb_n(D_r) {\buildrel \sim \over \longto} \CQ uot_n(r,n)$ mapping a quotient 
$\CO_{S_r}\twoheadrightarrow Q$ to the corresponding quotient 
$(O_r, \Lambda_r) \twoheadrightarrow (G,\Upsilon)$. 
\end{prop} 

{\it Proof.} This follows by a routine verification for flat families. 

\hfill $\Box$

\begin{rema}\label{kerrem} 
Note that the kernel of the epimorphism \eqref{eq:higgsquotA} is naturally a framed Higgs sheaf $(E,\Phi)$ on $S$
with topological  invariants 
\[ 
{\rm rk}(E)=r, \qquad c_1(E)= 0, \qquad \langle c_2(E), [S]\rangle = n.
\]
Let $\CH(r,n)$ denote the moduli stack of all such framed Higgs sheaves.
Then Proposition \ref{hilbhiggspropA} 
yields a stack morphism 
\be\label{eq:hilbmorphism} 
h:Hilb_n(D_r)\to \CH(r,n). 
\ee
In the next section it will be shown that this is in fact a closed embedding of schemes in the case where $S$ is the complex projective space, $L = \CO_S$ and $\Delta \subset S$ is a projective line. 
\end{rema}

\section{Framed Higgs sheaves on the projective plane}\label{frhiggsplane}  

In this section $S=\IP^2$ 
with homogeneous coordinates $[z_0, z_1, z_2]$, the line bundle $L$ is trivial, $L=\CO_S$, and 
$\Delta \subset S$ is the projective line $z_0=0$.
Then the Higgs bundle $(O_r, \Lambda_r)$ obtained in Lemma 
\ref{frHiggslemmB} has underlying vector bundle $O_r =\IC^r\otimes \CO_S$. The Higgs field $\Lambda_r:O_r\to O_r$ is given by 
\[
\Lambda_r = A_r \otimes {\bf 1}_{\CO_S}
\]
where $A_r\in M_r(\IC)$ is the lower triangular  $r\times r$ regular nilpotent Jordan block. As in the previous section, $\CH(r,n)$ denotes the resulting moduli stack of framed Higgs sheaves on $S$ with rank $r\geq 1$ and second Chern number $n$. The framing condition implies that the first Chern class must vanish. Let $\CM(r,n)$ denote the moduli space of framed sheaves on $S$, which is smooth, quasi-projective. Then note that there is a natural forgetful morphism $\jmath: \CH(r,n)\to \CM(r,n)$ forgetting the Higgs field. In this section we will prove that the moduli stack $\CH(r,n)$ in fact is isomorphic to a moduli scheme of framed quiver representations and the morphism $\jmath$ is a closed embedding.

\subsection{ADHM data for framed Higgs 
sheaves}\label{ADHMsection} 
Consider the following quiver $\CQ$ 
\be
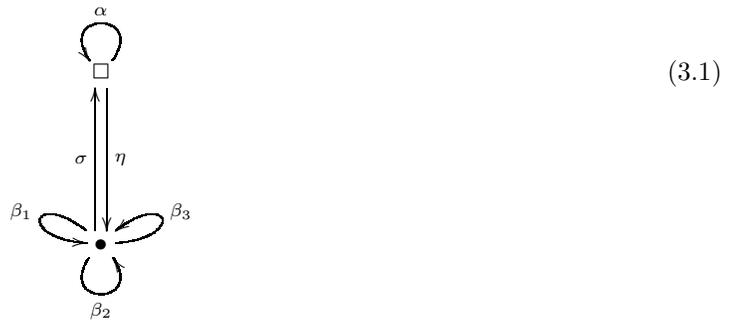
\label{eq:quiver}
\xymatrix{ 
& \Box \ar@(ur,ul)_-{\alpha}  \ar@<0.5ex>[dd]^-{\eta} & \\
& & \\
& \bullet \ar@(ul,l)_{\beta_1} \ar@(dl,dr)_{\beta_2} 
\ar@(r,ur)_{\beta_3} \ar@<0.5ex>[uu]^-{\sigma} \\}
\ee
with potential 
\[
\CW=\beta_3\circ(\beta_1\circ \beta_2-\beta_2\circ \beta_1) 
+ \beta_3 \circ \eta\circ \sigma -  \eta\circ\alpha\circ \sigma
\]
As usual, a representation  of the above quiver with potential is defined as a pair of vector spaces $(V,V_\infty)$ and linear maps 
\[
A: V_\infty \to V_\infty, \qquad B_i: V\to V,\ \ 1\leq i \leq 3,  \qquad I: V_\infty\to V, \qquad 
J:V\to V_\infty
\]
satisfying the relations derived from the potential function
\be\label{eq:potfct}
W = {\rm Tr}_V(B_3[B_1, B_2] + B_3IJ - IAJ).
\ee
The dimension vector of such a representation is the pair of integers $(r,n)$, where 
\[
n={\rm dim}(V) \qquad r={\rm dim}(V_\infty). 
\]
An isomorphism between two such representations 
$\varrho, \varrho'$ is a pair of vector space isomorphisms $\phi:V {\buildrel \sim \over \longto} V'$, $\psi: V_{\infty} {\buildrel \sim \over \longto} V'_{\infty}$ 
intertwining between the linear maps belonging to the two representations. 

\begin{defi}\label{framedQreps} 
A representation  of dimension vector $(r,n)$ will be called framed if  the following conditions are satisfied: 
\begin{itemize}
\item[$(Fr.1)$] The vector space $V_{\infty}$ is a fixed vector space of dimension $r$, 
equipped with a fixed basis ${\bf e}_{a}$, $1\leq a\leq r$. 
\item[$(Fr.2)$] The map
$A:V_{\infty}\to V_{\infty}$ is fixed and given by the regular nilpotent endomorhism
\be\label{eq:regnilpA} 
A_r ({\bf e}_a) = {\bf e}_{a+1}, \qquad 1\leq a \leq r,
\ee
 where ${\bf e}_{r+1}=0$. In particular if $0\leq r \leq 1$, the map $A$ is identically zero. 
\item[$(Fr.3)$] The representation  satisfies the relations 
derived from the potential function 
\be\label{eq:potentialC}
W^{fr} = {\rm Tr}_V \bigg( B_3[B_1,B_2] + B_3IJ -I A_r J \bigg),
 \ee
 where the map $A_r$ is fixed as in 
 $(Fr.1)$, $(Fr.2)$ above. 
 The resulting relations are:
 \be\label{eq:relationsA}
\bal
 & [B_1,B_2]+IJ =0,\qquad JB_3-A_rJ=0, \qquad B_3I-IA_r=0.\\
 & [B_3,B_1]=0, \qquad \  \qquad [B_3, B_2]=0\\
\eal
\ee
\end{itemize} 
Furthermore, 
a framed representation  as defined above 
will be called {\it cyclic} if  it satisfies the additional condition:
\begin{itemize} 
\item[$(Fr.4)$] There is no proper non-zero linear subspace $0\subset V'\subset V$ preserved by $B_1,B_2, B_3$ and at the same time containing the image of $I:V_{\infty}\to V$.
\end{itemize}
\end{defi}
Finally, in the formulation of the moduli problem, two framed representations $\varrho, \varrho'$ will be said to be isomorphic if and only if they are related by an isomorphism of the form 
$(\phi, {\bf 1}_{V_{\infty}})$. It is clear that such isomorphisms preserve the fixed framing data.

At this point it is helpful to note the relation between the quiver \eqref{eq:quiver} and the standard ADHM quiver. The latter is defined by the following diagram
\[
\xymatrix{ 
& \Box  \ar@<0.5ex>[dd]^-{\eta} & \\
& & \\
& \bullet \ar@(ul,dl)_{\beta_1} \ar@(ur,dr)^-{\beta_2} 
 \ar@<0.5ex>[uu]^-{\sigma} \\}
\]
subject to the quadratic relation 
\[
\beta_1\circ \beta_2-\beta_2\circ \beta_1
+  \eta\circ \sigma. 
\]
By analogy with $(Fr.1)$--$(Fr.4)$ above, 
a fixed isomorphism $V_\infty \simeq \IC^r$ will be fixed by choosing a basis, where $V_\infty$ is the vector space assigned to the node $\Box$. Then framed representations of the ADHM quiver are defined by 
data 
$\alpha=(V, V_\infty, B_1, B_2,I,J)$ where $B_i \in {\rm End}(V)$, $1\leq i\leq 2$, $I \in {\rm Hom}(V_\infty, V)$ and $J \in {\rm Hom}(V, V_\infty)$ 
are linear maps corresponding to $\beta_i, \eta, \sigma$ respectively. Hence they
satisfy the quadratic relation $[B_1,B_2]+IJ=0$. Isomorphisms of framed representations are required to preserve the identification $V_\infty \simeq \IC^r$. Moreover, recall that $\alpha$ is cyclic, or stable, if and only if there is no proper nontrivial linear subspace $V'\subset V$ preserved by $B_1, B_2$ and at the same time containing the image of $I$. As shown for example in Section 3.1 of \cite{Hilb_lect}, there is a smooth quasi-projective fine moduli space $\CA(r,n)$ of stable framed ADHM representations of fixed dimension vector $(n,r)$. 

Now note that 
 for each framed representation $\varrho=(V,V_\infty,B_1,B_2, B_3, I,J)$ defined in $(Fr.1)$--$(Fr.4)$ above, the data $\alpha=(V, V_\infty, B_1, B_2,I,J)$ is a framed representation of the ADHM quiver. A priori $\alpha$ need not be cyclic even if $\varrho$ is. The next result will show that in fact the two cyclicity conditions are in fact compatible.

\begin{lemm}\label{stablemma}
Let $\varrho$ be a framed representation of $(\CQ,\CW)$ satisfying conditions 
$(Fr.1)$ - $(Fr.3)$. Let $\alpha$ be the underlying ADHM representation. Then
$\varrho$ is cyclic as defined in $(Fr.4)$ if and only is $\alpha$ is cyclic as an 
ADHM representation. 
\end{lemm} 

{\it Proof}. The inverse implication is clear. In order to prove the direct implication,  
suppose $\varrho$ is cyclic, and suppose $V'\subset V$ is a proper nonzero linear subspace preserved by $B_1, B_2$  and at the same time containing the image of $I$. Let $v_a = I({\bf e}_a)$, $1\leq a \leq r$. Since $\varrho$ is cyclic, $V$ is generated by elements of the form 
\[
m(B_1,B_2,B_3)v_a 
\]
where $m(B_1,B_2,B_3)$ are monomials in ${\rm End}(V)$. Hence $V'$ will be generated by certain linear combinations of such elements. Then relations \eqref{eq:relationsA} imply that $V'$ is also preserved by $B_3$. 

\hfill $\Box$

The construction of the moduli space of framed cyclic representations 
$(V,B_1,B_2,B_3,I,J)$ with fixed dimension vector is very similar to the GIT construction of the moduli space of ADHM quiver representations presented in Section 3.1 of \cite{Hilb_lect}. In particular, proceeding by analogy with loc. cit., the cyclicity condition $(Fr.4)$ is equivalent to a GIT stability condition, and the following holds. 
\begin{prop}\label{frquivpropO} 
For any $n, r \in \IZ$, $n,r \geq 1$ there is a quasi-projective fine moduli space $\CQ(r,n)$ of framed cyclic representations $(V,B_1,B_2,B_3,I,J)$ of dimension vector $(r,n)$. In particular any such representation has trivial automorphism group. Furthermore there is a natural forgetful morphism 
$\imath : \CQ(r,n)\to \CA(r,n)$ to the moduli space of stable framed ADHM representations obtained by omitting the map $B_3:V\to V$. 
\end{prop} 

The connection to framed Higgs sheaves is provided by: 

\begin{prop}\label{frquivpropB} 
There is a commutative diagram of morphisms of moduli stacks
\be\label{eq:modspdiag} 
\xymatrix{ 
\CH(r,n) \ar[r]^-{\jmath} \ar[d]^-{\simeq} & \CM(r,n) 
\ar[d]^-{\simeq} \\ \CQ(r,n)\ar[r]^-{\imath} & \CA(r,n)\\ }
\ee
where the vertical arrows are isomorphisms. 
\end{prop} 

{\it Proof}. This follows from the construction of the left vertical isomorphism in diagram \eqref{eq:modspdiag} given in \cite[Thm. 2.1]{Hilb_lect}, which is based on the Beilinson spectral sequence. Then note that the latter is functorial with respect to morphisms of sheaves, and 
all the  steps carried out in Section 2.1 of loc. cit. are naturally compatible with morphisms of framed sheaves. 

\hfill $\Box$.

\subsection{Embedding in the ADHM moduli space}\label{embedding}

The main technical result of this section states that the forgetful 
morphism $\imath : \CQ(r,n)\to \CA(r,n)$ is a closed embedding. 
This will be shown in several steps using the criterion proved in 
\cite[Lemma 4]{Schematic_HN}. In order to formulate the required conditions, note that the map $\imath$ naturally preserves residual fields. That is, if $\varrho_K$ is a point of $\CN(r,n)$ 
with residual field $K$, then the residual field of the point 
$\alpha_K=\imath(\varrho_K)$ is canonically isomorphic to  $K$. 
This implies that there exists an induced map on Zariski tangent spaces 
\be\label{eq:tangentmapA} 
\iota_*: T_{\varrho_K}\CQ(r,n) \to T_{\iota(\varrho_K)}\CA(r,n).
\ee
Then,  as shown in \cite[Lemma 4]{Schematic_HN}, it suffices to prove that: 
\begin{itemize}
\item[(1)] $\imath$ is universally injective i.e. for any point
$\alpha_K$ of $\CA(r,n)$ with residual field $K$ there is at most one point $\varrho_K$ of $\CN(r,n)$, also with residual field $K$, such that $\imath(\varrho_K)=\alpha_K$. 
\item[(2)] $\imath$ is proper.
\item[(3)] For any point $\varrho_K$ of $\CQ(r,n)$  the induced map on Zariski tangent spaces \eqref{eq:tangentmapA} 
is injective. 
\end{itemize} 
This will be carried out in three separate lemmas. Since $\imath$ is the identity map for $r=1$, it will be assumed $r\geq 2$ below. 

\begin{lemm}\label{frquivlemmA} 
$\imath$ is universally injective.
\end{lemm} 

{\it Proof}.  
This follows easily from the cyclicity condition $(Fr.4)$. Let $\varrho_K$ be an arbitrary framed cyclic representation defined over a $\IC$-field $K$ 
of characteristic zero. Then its underlying vector space $V_K$ is the $K$-linear span of 
\[
m(B_{1,K}, B_{2,K}) I_K
\]
where $m(B_{1,K}, B_{2,K})$ is an arbitrary monomial in the endomorphism ring ${\rm End}_K(V_K)$.  This uniquely determines $B_{3,K}$ via the relations 
\be\label{eq:relationsB}
[B_{3,K}, B_{i,K}]=0,\ 1\leq i\leq 2,  \qquad B_{3,K}I_K = I_K\, (A_r\otimes {\bf 1}_{K}). 
\ee

\hfill $\Box$. 

\begin{lemm}\label{frquivlemmB} 
$\imath$ is proper. 
\end{lemm} 

{\it Proof.} 
This will be proven using the valuative criterion for properness. Let $R$ be a discrete valuation ring over $\IC$, let ${\mathfrak m}\subset R$ denote the unique maximal ideal, and $K$ its field of fractions. Using the categorical equivalence between 
$R$-modules and sheaves on  ${\rm Spec}(R)$, flat families of stable representations of the ADHM quiver parametrized by ${\rm Spec}(R)$ are in one-to-one correspondence to representations 
 $\alpha_R=(V_R, B_{1,R}, B_{2,R}, I_R, J_R)$ the ADHM quiver over $R$ satisfying the following conditions.  
\begin{itemize} 
\item $V_R$ is a free $R$-module, and $B_{1,R}, B_{2,R} \in {\rm End}_R(V_R)$, $I_R \in {\rm Hom}_R(R^{\oplus r}, V_R)$, 
$J_R\in {\rm Hom}_R(V_R, R^{\oplus r})$ are morphisms of $R$ modules such that 
\be\label{eq:relationsC}
[B_{1,R}, B_{2,R}]+I_RJ_R=0.
\ee 
\item For any prime ideal ${\mathfrak p}\subset R$ the induced ADHM representation $\alpha_{k({\mathfrak p})} = \alpha_R \otimes_R R/{\mathfrak p}$ over the residual field $k({\mathfrak p})$ is stable. 
\end{itemize}
For simplicity let $\alpha = (V, B_1, B_2, I,J)$ denote 
$\alpha_{k({\mathfrak m})}$, which is a 
complex ADHM representation. Let also $I_{R,a}: R \to V_R$, $1\leq a\leq r$ 
denote the components of $I_{R}$. Analogous notation will be used for the induced representations $\alpha_K, \alpha$. 

Since $V$ is a finite dimensional  complex vector space, it admits a finite set of generators 
\[
m_{i,a}(B_{1}, B_{2})I_a, \qquad 1\leq i \leq N_a, \quad 1\leq a\leq r,
\]
for some positive integers $N_a\geq 1$. Therefore there is an exact sequence of $\IC$-vector spaces 
\[
0 \to Y {\buildrel f\over \longto}  \oplus_{a=1}^r \oplus_{i=1}^{N_a} \IC\langle {\bf u}_{i,a}\rangle  {\buildrel p\over \longto} V\to 0
\]
where $p$ is the canonical projection. By Nakayama's lemma the 
elements 
 \[
m_{i,a}(B_{1,R}, B_{2,R})I_R, \qquad 1\leq i \leq N_a , \quad 1\leq a\leq r,
\]
generate $V_R$ as an $R$-module. In particular, the elements 
\[
m_{i,a}(B_{1,K}, B_{2,K})I_K, \qquad 1\leq i \leq N_a , \quad 1\leq a\leq r,
\]
generate $V_K$ as a $K$-vector space. Moreover there is an exact sequence of finitely generated $R$-modules 
\[ 
0 \to Y_R {\buildrel f_R\over \longto}  \oplus_{a=1}^r \oplus_{i=1}^{N_a} R\langle {\bf u}_{i,a}\rangle {\buildrel p_R\over \longto} V_R\to 0
\]
as well as an exact sequence of $K$-vector spaces 
\[
0 \to Y_K{\buildrel f_K\over \longto}  \oplus_{a=1}^r \oplus_{i=1}^{N_a} K\langle {\bf u}_{i,a}\rangle {\buildrel p_K\over \longto} V_K\to 0
\]
where $p_R,p_K$ are again the canonical projections. Since $V_R$ is free, it follows that $Y_R$ is isomorphic to a free module as well. 

Suppose $B_{3,K}: V_K \to V_K$ is a $K$-linear map satisfying 
relations \eqref{eq:relationsB}. In detail, this is equivalent to 
\[ 
B_{3,K}(m_{i,a}(B_{1,K},B_{2,K})I_{K,a}) = m_{i,a}(B_{1,K},B_{2,K})I_{K,a+1}
\]
for all $1\leq a\leq r-1$, $1\leq i \leq N_a$ and 
\[ 
{B}_{3,K}(m_{i,r}(B_{1,K},B_{2,K})I_{K,r})) =0
\]
 for all $1\leq i \leq N_r$. 
Let 
\[
{\overline B}_{3,K}: 
\oplus_{a=1}^r \oplus_{i=1}^{N_a} K\langle {\bf u}_{i,a}\rangle \to \oplus_{a=1}^r \oplus_{i=1}^{N_a} K\langle {\bf u}_{i,a}\rangle 
\]
be the linear map determined on basis elements by 
\[ 
{\overline B}_{3,K}({\bf u}_{i,a}) = {\bf u}_{i,a+1}
\]
for all $1\leq a\leq r-1$, $1\leq i \leq N_a$ and 
\[ 
{\overline B}_{3,K}({\bf u}_{i,r}) =0
\]
for all $1\leq i \leq N_r$. By construction there is a commutative diagram 
\be\label{eq:Kdiagram}
\xymatrix{ 
\oplus_{a=1}^r \oplus_{i=1}^{N_a} K\langle {\bf u}_{i,a}\rangle
\ar[rr]^-{{\overline B}_{3,K}} \ar[d]^-{p_K} & & 
\oplus_{a=1}^r \oplus_{i=1}^{N_a} K\langle {\bf u}_{i,a}\rangle
\ar[d]^-{ p_K} \\
V_K \ar[rr]^-{B_{3,K}} & & V_K. \\
}
\ee
In particular, 
\be\label{eq:Kvanishing}
p_K \circ {\overline B}_{3,K} \circ f_K =0.
\ee
Clearly, ${\overline B}_{3,K}$ extends to an $R$-module morphism 
\[
{\overline B}_{3,R}:\oplus_{a=1}^r \oplus_{i=1}^{N_a} R\langle {\bf u}_{i,a}\rangle \to \oplus_{a=1}^r \oplus_{i=1}^{N_a} R\langle {\bf u}_{i,a}\rangle 
\]
such that
\[ 
B_{3,R}(m_{i,a}(B_{1,R},B_{2,R})I_{K,a}) = m_{i,a}(B_{1,R},B_{2,R})I_{R,a+1}
\]
for all $1\leq a\leq r-1$, $1\leq i \leq N_a$ and 
\[ 
{B}_{3,R}(m_{i,r}(B_{r,R},B_{r,R})I_{R,r})) =0
\]
 for all $1\leq i \leq N_r$. 
Furthermore, relation \eqref{eq:Kvanishing} implies 
\[
p_R \circ {\overline B}_{3,R} \circ f_R =0.
\]
Hence $p_R \circ {\overline B}_{3,R}$ yields a morphism of $R$-modules $B_{3,R} : V_R \to V_R$. By construction, this also satisfies the relations 
\[ 
B_{3,R}(m_{i,a}(B_{1,R},B_{2,R})I_{R,a}) = m_{i,a}(B_{1,R},B_{2,R})I_{R,a+1}
\]
for all $1\leq a\leq r-1$, $1\leq i \leq N_a$ and 
\[ 
{B}_{3,R}(m_{i,r}(B_{1,R},B_{2,R})I_{R,r})) =0
\]
 for all $1\leq i \leq N_r$. Therefore this is the required extension. It is also unique by Lemma \ref{frquivlemmA}. 

\hfill $\Box$ 

Finally, the third part of the proof consists of: 
\begin{lemm}\label{frquivlemmC} 
For any point $\varrho_K$ of $\CQ(r,n)$ the induced map on Zariski tangent spaces $\imath_*: T_{\varrho_K}\CQ(r,n) \to T_{\imath(\varrho_K)}\CA(r,n)$  is injective.
\end{lemm} 

{\it Proof}. 
The Zariski tangent spaces for moduli of quiver representations are canonically determined by linearizing the specified relations. 
In the present case, $T_{\varrho_K}\CQ(r,n)$ is isomorphic to the
middle cohomology of the following complex of amplitude $[0,\ 2)$:
\[
\xymatrix{ 
0\ar[d] \\
{\rm End}(V) \ar[d]^-{g_0} \\
{\rm End}(V)^{\oplus 3}\oplus {\rm Hom}(V_{\infty}, V) 
\oplus {\rm Hom}_K(V, V_{\infty}) \ar[d]^-{g_1} \\
{\rm End}(V)^{\oplus 3}\oplus {\rm Hom}(V_{\infty}, V) 
\oplus {\rm Hom}(V, V_{\infty})
}
\]
where 
\[
g_0(\psi) = \left([\psi, B_i],\ \psi I, J\psi\right),  
\]
\[
\bal
& g_1(\epsilon_i, \eta, \delta) = \\
& \left([B_1,\epsilon_3],\ [B_2, \epsilon_3],\ 
[\epsilon_1, B_2]+[B_1,\epsilon_2]+I\delta + \eta J, \ 
B_3\epsilon +\epsilon_3I-\epsilon A_r, \ \delta B_3+J\epsilon_3-A_r\delta \right) \\
\eal 
\]
for $1\leq i\leq 3$.
For simplicity, the subscript $K$ in the notation used for the 
data of the representation $\varrho_K$ has been suppressed. This convention will be employed only in the proof of the current lemma. 

At the same time 
$T_{\imath(\varrho)}\CA(r,n)$ is isomorphic to the
middle cohomology of the following complex of amplitude $[0,\ 2)$:
\[
\xymatrix{ 
0\ar[d] \\
{\rm End}(V) \ar[d]^-{f_0} \\
{\rm End}(V)^{\oplus 2}\oplus {\rm Hom}(V_{\infty}, V) 
\oplus {\rm Hom}(V, V_{\infty}) \ar[d]^-{f_1} \\
{\rm End}(V)
}
\]
where 
\[
f_0(\psi) = \left([\psi, B_i],\ \psi I, J\psi\right), 
\]
\[
\bal
& f_1(\epsilon_i, \eta, \delta) = 
[\epsilon_1, B_2]+[B_1,\epsilon_2]+I\delta + \eta J\\
\eal 
\]
for $1\leq i \leq 2$. 
Note that there is a natural degree zero map of complexes determined by the obvious term-by-term canonical injections. 
Moreover, 
the cyclicity condition implies that $f_0, g_0$ are injective by a straightforward argument. 

Next note that the induced map on degree $1$ cochains is also injective. 
Suppose a degree 1 cochain $(\epsilon_i, \epsilon, \delta)$ of the first complex maps to $0$ in the second complex. This implies 
$\epsilon =0$ and $\eta=0$. Therefore the condition 
$g_1(\epsilon_i, \epsilon, \delta) =  0$ yields 
\[
[B_i, \epsilon_3]=0,\quad 1\leq i\leq 2, \qquad \epsilon_3 I =0, \qquad J\epsilon_3 =0 
\]
Then the cyclicity condition implies that $\epsilon_3$ must be identically zero. 
Since $f_0, g_0$ are injective, this implies that the induced map in degree 1 cohomology is also injective. 

\hfill $\Box$

In conclusion, using \cite[Lemma 4]{Schematic_HN}, Lemmas \ref{frquivlemmA}, \ref{frquivlemmB}, 
\ref{frquivlemmC} prove: 

\begin{prop}\label{frquivpropA}
The morphism $\imath: \CQ(r,n) \to \CA(r,n)$ is a closed embedding. 
\end{prop} 

\subsection{The Hilbert scheme as a quiver moduli space}\label{hilbnonred} 

In order to simplify the notation, let $Hilb(r,n)= Hilb_n(D_r)$ since 
$D_r$ will always be a nonreduced plane in the following. 
The main goal of this section is to show that the natural forgetful morphism $h: Hilb(r,n)\to \CH(r,n)$ obtained from Proposition \ref{hilbhiggspropA} is also a  closed embedding. More precisely using the isomorphism 
$\CH(r,n) {\buildrel \sim \over \longto} \CQ(r,n)$ constructed in Proposition  \ref{frquivpropB}, it will be shown that $h$ yields an isomorphism onto the closed subscheme of $\CQ_0(r,n)\subset \CQ(r,n)$ parametrizing representations 
$\varrho$ with $J=0$.  
The main idea of the proof is to show that the 
framed Higgs sheaf corresponding to a framed cyclic representation $\varrho=(V,B_1,B_2,B_3, I, J)$
is a Higgs sheaf quotient as in Proposition \ref{hilbhiggspropA} if and only if $J=0$.
This will require some intermediate steps. In the next two lemmas the ground field will be an arbitrary $\IC$-field $K$, although the index $K$ will be suppressed in order to simplify the notation. 

Let $\alpha = (V,B_1,B_2, I,J)$ be the underlying ADHM representation of $\varrho$. As shown in \cite[Section 2.1]{Hilb_lect}, the framed sheaf corresponding to $\alpha$ is the middle cohomology sheaf of the monad complex 
\[ 
0\to V\otimes \CO_{S}(-\Delta) \xlongrightarrow{f_{-1}}
\begin{array}{c} V\otimes \CO_S\\ \oplus \\ V\otimes \CO_S\\
\oplus \\
V_\infty\otimes \CO_S \end{array} \xlongrightarrow{f_0} 
V\otimes \CO_S(\Delta) \to 0 
\]
where the terms have degrees $-1,0,1$ and the differentials are 
given by 
\[
f_{-1}=\left(\begin{array}{c} z_0 B_1- z_1 \\ z_0 B_2- z_2 \\
z_0 J \end{array}\right) \qquad 
f_0=\left(-z_0B_2+z_2,\ z_0B_1-z_1,\\ z_0 I\right).
\]
This complex will be denoted by $\CC_\alpha$ in the following. 
As proven in \cite[Lemma 2.7]{Hilb_lect}, it has trivial cohomology in degrees $-1,1$. These statements are proven in loc. cit for ground field $\IC$, but the 
generalization to an arbitrary ground field $K$ is straightforward. 

Now let $\CB_\alpha$ denote the complex 
\[
0\to V\otimes \CO_{S}(-\Delta) \xlongrightarrow{g_{-1}} 
\begin{array}{c} V\otimes \CO_S\\ \oplus \\ V\otimes \CO_S\\
\end{array}  \xlongrightarrow{g_{0}} 
V\otimes \CO_S(\Delta) \to 0 
\]
where, again, the terms have degrees $-1,0,1$ and the differentials are 
given by 
\[
g_{-1}=\left(\begin{array}{c} z_0 B_1- z_1 \\ z_0 B_2- z_2  \end{array}\right) \qquad 
g_0=\left(-z_0B_2+z_2,\ z_0B_1-z_1\right).
\]
If $J=0$ there is an exact sequence of complexes 
\[
0\to \CB_\alpha \to \CC_\alpha \to V_\infty\otimes \CO_S\to 0
\]
where the map $\CB_\alpha \to \CC_\alpha$ is the natural injection in each degree. 
It will be shown next that $\CB_\alpha$ has trivial cohomology in degrees $-1,0$. 

\begin{lemm}\label{quivhilbA} 
Suppose $\alpha=(V,B_1,B_2,I,J)$ is a stable ADHM representation. 
Then the complex $\CB_\alpha$ has trivial cohomology in degrees $-1,0$
while its degree 1 cohomology sheaf has zero dimensional support contained in $S\setminus \Delta$. 
\end{lemm} 

{\it Proof}. First note that the restriction of $\CB_\alpha$ to 
$\Delta$ is exact since $z_1|_\Delta, z_2|_\Delta$ do not have any common zeroes. This implies that 
$$H^1(\CB_\alpha)|_\Delta\simeq H^1\big(\CB_\alpha|_\Delta\big)$$ 
is zero, hence $H^1(\CB_\alpha)$ must be a torsion sheaf on $S$ supported in the complement $S\setminus \Delta$. Since $S\setminus \Delta \simeq \IA^2_K$, this implies that 
$H^1(\CB_\alpha)$ must be a zero dimensional sheaf, which further implies that 
\be\label{eq:cohrankA}
{\rm rk}\, H^{-1}(\CB_\alpha) = {\rm rk}\, H^0(\CB_\alpha). 
\ee
Next note that $g_{-1}$ is injective by the same argument as in \cite[Lemma 2.7(1)]{Hilb_lect}. Namely, suppose $g_{-1}$ is not injective. Then its kernel must be a nonzero torsion-free sheaf and there is an exact sequence of $\CO_S$-modules
\[
0 \to {\rm Ker}(g_{-1}) \to V\otimes \CO_S \to {\rm Im}(g_{-1}) \to 0.
\]
Moreover, ${\rm Im}(g_{-1}) \subset V\otimes \CO_S^{\oplus 2}$ is a nonzero subsheaf since $g_{-1}$ is not identically zero. Therefore ${\rm Im}(g_{-1})$ is a also a nonzero torsion-free sheaf, hence locally free on the complement of a zero dimensional subscheme $Y\subset S$. This implies that ${\rm Ker}(g_{-1})$ is locally free at any $K$-point $p\in S\setminus Y$ and the fiber ${\rm Ker}(g_{-1})|_p$ is the kernel of the restriction
$g_{-1}|_p$. As in \cite[Lemma 2.7(1)]{Hilb_lect}, this further implies  that ${\rm Ker}(g_{-1})|_p$ is a simultaneous eigenspace 
of $B_1,B_2$ with eigenvalues $\lambda_1, \lambda_2$ determined by 
\[
z_0|_p \lambda_i = z_i|_p, \qquad 1\leq i \leq 2. 
\]
Since $B_1, B_2$ are fixed, this condition can only hold for finitely many 
$K$-points $p$, leading to a contradiction. 
In conclusion $g_{-1}$ is injective.  

In order to prove exactness in the degree zero, suppose 
the middle cohomology sheaf $H^0(\CB_\alpha)$ is not zero. 
Then the snake lemma yields an exact sequence 
\[
0\to H^0(\CB_\alpha) \to {\rm Coker}(g_{-1}) \to {\rm Im}(g_0) \to 0.
\]
At the same time, since $g_{-1}$ is injective, ${\rm Coker}(g_{-1})$ has a two term locally free resolution, which implies that it is torsion-free. Hence $H^0(\CB_\alpha)$ must be a nonzero torsion-free sheaf. However, since $H^{-1}(\CB_\alpha)=0$,  equation \eqref{eq:cohrankA} implies that ${\rm rk}\, H^0(\CB_\alpha)=0$, leading to a contradiction.
In conclusion, $H^0(\CB_\alpha)$ is also zero.  

\hfill $\Box$

\begin{lemm}\label{quivhilbB} 
Let $\varrho=(V,B_{1},B_{2},B_{3}, I,J)$ be a cyclic framed representation
of $(\CQ,\CW)$. Then the framed Higgs sheaf corresponding to $\varrho$ via Proposition \ref{frquivpropA} is a Higgs sheaf quotient as in Proposition \ref{hilbhiggspropA} if and only if $J=0$. 
\end{lemm} 

{\it Proof.} 
$(\Rightarrow)$ Suppose $J=0$.
As above, the framed sheaf $E_\alpha$ corresponding to the ADHM representation $\alpha=(V,B_1,B_2,I,J)$ is the middle cohomology sheaf of the monad complex $\CC_\alpha$. Since $J=0$, there is an exact sequence of complexes 
\[ 
0 \to \CB_\alpha \to \CC_\alpha \to V_\infty \otimes \CO_S \to 0 
\]
where the last term from the left is regarded as a complex supported in degree 0. Using Lemma \ref{quivhilbA}, this yields an exact sequence
\[
0\to E_\alpha \to V_\infty \otimes \CO_S \to G_\alpha \to 0 
\]
where $G_\alpha=H^1(\CB_\alpha)$ is a zero dimensional $\CO_S$-module. 
Using the monad construction, the map $B_3:V\to V$ yields morphisms $\Phi: E_\alpha\to E_\alpha$, respectively $\Upsilon: G_\alpha\to G_\alpha$ which are naturally compatible with the differentials of the above complex. Moreover, by construction, $\Phi$ also satisfies the framing condition \eqref{eq:frHiggsD}.  Therefore one obtains indeed an exact sequence of framed Higgs sheaves 
\[
0 \to (E_\alpha, \Phi) \to (O_r, \Lambda_r) \to (G_\alpha, \Upsilon) \to 0. 
\]

$(\Leftarrow)$ The above steps are reversible using the functoriality of the Beilinson spectral sequence.

\hfill $\Box$. 

To conclude, let $\CQ_0(r,n)$ denote the closed subscheme of the moduli scheme 
$\CQ(r,n)$ parameterizing cyclic framed representations $\varrho$ with $J=0$.
Let $\CH_0(r,n)$ denote the corresponding closed subscheme of 
the framed Higgs moduli space $\CH(r,n)$ under the isomorphism 
$\CQ(r,n) \simeq \CH(r,n)$. Proposition \ref{hilbhiggspropA} and Lemma \ref{quivhilbB} show that $h$ yields a bijection between the set of $K$-points of 
the Hilbert scheme and the set of $K$-points of $\CH_0(r,n)$ for any field $K$ over $\IC$.   In fact, a stronger statement holds: 

\begin{prop}\label{hilbhiggspropB}
The morphism $h: Hilb(r,n)\to \CH(r,n)$ constructed below 
Proposition \ref{hilbhiggspropA} factors through an isomorphism 
$h_0:Hilb(r,n){\buildrel \sim \over \longto} \CH_0(r,n)$. In particular, $h$ is a closed embedding.\end{prop}

{\it Proof}. It suffices to generalize Lemma \ref{quivhilbB} to flat families. 
This is 
straightforward since the monad construction works in the relative setting. The details are very similar to those in Sections 7.1 and 7.2 of \cite{ADHM_sheaves}, hence will be omitted. 

\hfill $\Box$

\section{Torus actions and fixed loci}\label{torusaction}

Summarizing the previous results, there is a commutative diagram of scheme morphisms 
\be\label{eq:modspdiag}
\xymatrix{
Hilb(r,n) \ar[r]^-{h} \ar[d]^-{\simeq} & \CH(r,n) \ar[r]^-{\jmath} \ar[d]^-{\simeq}  & \CM(r,n) \ar[d]^-{\simeq} \\
\CQ_0(r,n) \ar[r]^-{q} & \CQ(r,n) \ar[r]^-{\imath} & \CA(r,n) \\}
\ee
where all horizontal arrows are closed embeddings. The top row consists of geometric moduli spaces, while the bottom row consists of the framed quiver moduli spaces obtained from the Beilinson spectral sequence. 
The goal of this section is to construct certain torus actions on all moduli spaces in the above diagram such that all arrows are equivariant, and prove certain properties of the fixed loci.

\subsection{Generic torus action}\label{gentorus} 

First recall the natural action on the moduli space of framed sheaves of the $(r+2)$-dimensional torus, namely, ${\bf T}_{r+2} \times \CM(r,n)\to \CM(r,n)$ with ${\bf T}_{r+2}= \IC^\times \times \IC^\times \times (\IC^{\times})^{\times r}$. 
Namely, for any $(t_1,t_2, u_1, \dots, u_r)\in {\bf T}_{r+2}$ let 
$\eta(t_1, t_2): \IP^2 \to \IP^2$ be the morphism induced by 
\[
(t_1, t_2) \times [z_0,z_1,z_2]\mapsto [z_0, t_1z_1,t_2z_2]
\]
and let $u$ denote the diagonal $r\times r$ matrix with diagonal elements 
$u_1, \ldots, u_r$. Then the ${\bf T}_{r+2}$-action on $\CM(r,n)$ is given by 
\[
(t_1,t_2, u) \times (E,\xi) \mapsto (\eta(t_1,t_2)^{-1})^*(E,
u\xi).  
\]
This action does not preserve the moduli space of framed Higgs sheaves. 

Let ${\bf T}:={\bf T}_3$ denote the torus $\IC^\times \times \IC^\times \times \IC^\times$ with coordinates $(t_1,t_2,t_3)$.
The threefold $X$ in Section \ref{hilbhiggs} is canonically isomorphic to $\IP^2\times \IA^1$. Hence there is a three-dimensional torus action on the Hilbert scheme determined by the geometric action 
${\bf T}\times X \to X$,  
\[
(t_1,t_2,t_3) \times ([z_0,z_1,z_2], y) \mapsto ([z_0, t_1z_1,t_2z_2],t_3y)
\]
There is also a natural ${\bf T}$-action on framed sheaf moduli, as well as framed Higgs sheaf moduli, making all the morphisms in the above diagram equivariant. This is obtained by restricting the ${\bf T}_{r+2}$-action on 
$\CM(r,n)$ to a ${\bf T}$-action via the injective group morphism ${\bf T} \hookrightarrow {\bf T}_{r+2}$, 
\be\label{eq:Tsubtorus}
(t_1,t_2,t_3) \mapsto (t_1, t_2, t_3^{1-a}), \qquad 1\leq a \leq r.
\ee
It is straightforward to check that the resulting ${\bf T}$-action on 
$\CM(r,n)$ preserves the closed subscheme $\CH(r,n)$. 

The corresponding ${\bf T}$-actions on the framed quiver moduli spaces $\CN(r,n)$ and $\CA(r,n)$ are then easily obtained from the monad construction. 
As in the previous section, the components of $I, J$ with respect to the fixed basis of $V_\infty$ will be denoted by $I_a: \IC \to V$, $J_a: V \to \IC$, 
$1\leq a \leq r$. Then the potential function is rewritten as  
\be\label{eq:potentialC}
W^{fr} = {\rm Tr}_V \big( B_3[B_1,B_2] \big)+ 
\sum_{a=1}^r \big(J_aB_3I_a - J_{a-1}I_a \big),
 \ee
where by convention $J_{0}=0$. The ${\bf T}$-action on $\CN(r,n)$ will be given by 
\be\label{eq:toractA}  
(t_i)\times (B_i, I_a, J_a) \mapsto 
(t_i B_i, t_3^{a-1} I_a, t_1t_2t_3^{1-a} J_a), \qquad 
1\leq i \leq 3, \quad 1\leq a \leq r. 
\ee
while the ${\bf T}$-action on $\CA(r,n)$ will be given by the same 
expression with $B_3$ is omitted. 

The next goal is to determine the {\bf T}-fixed points in the framed quiver moduli space $\CQ(r,n)$ up to $GL(n, \IC)$ gauge transformations. The 
fixed locus $\CA(r,n)^{{\bf T}_{r+2}}$ is finite and in one-to-one correspondence to $r$-partitions $\mu=(\mu_1, \ldots, \mu_r)$ of $n$, as shown for example \cite[Prop. 2.9]{instcountA}. 
Geometrically, the fixed points in the moduli space of framed sheaves $\CM(r,n)$ parametrize framed sheaves of the form 
\[ 
E_\mu \simeq \bigoplus_{a=1}^r I_{Z_{\mu_a}} 
\]
where $I_{Z_a}$ is the ideal sheaf of the ${\bf T}_{r+2}$-invariant zero dimensional subscheme $Z_{\mu_a}$ parameterized by the partition $\mu_a$. The framing $\xi_\mu$ is the natural framing determined by the injections $I_{Z_{\mu_a}}\subset \CO_S$, $1\leq a\leq r$. 

Partitions will be identified to Young diagrams using the conventions in Section 0.1 of \cite{Cherednik_W}. 
A partition $\nu=(\nu^1\geq \cdots 
\geq \nu^{l})$ will be identified with the set of integral points 
\[
\{(i,j) \in \IZ\times \IZ\, |\, 1\leq i \leq l,\ 1\leq j \leq \nu^i\}
\]
consisting of $l$ columns of heights $\nu^i$, $1\leq i \leq l$. As usual, 
such a set is also canonically identified with a collection of boxes. 
Then note: 

\begin{prop}\label{frquivpropC}
$(i)$ The ${\bf T}$-fixed locus $\CA(r,n)^{\bf T}$ coincides with $\CA(r,n)^{{\bf T}_{r+2}}$. 

$(ii)$ A {\bf T}-fixed point $\alpha_\mu \in \CA(r,n)^{\bf T}$ belongs to 
$\CQ(r,n)^{\bf T}$ if and only if the $r$-partition $\mu$ is nested i.e. 
\be\label{eq:nestedpart}
\mu_r\subseteq \mu_{r-1} \subseteq \cdots \subseteq \mu_1 
\ee
as Young diagrams. Hence the {\bf T}-fixed locus $\CQ(r,n)^{\bf T}$ is finite 
and in one-to-one correspondence to nested $r$-partitions. 

$(iii)$ 
The natural closed embedding $\eta:\CQ_0(r,n) \hookrightarrow \CQ(r,n)$ yields an isomorphism of ${\bf T}$-fixed loci. Hence the {\bf T}-fixed locus $\CQ_0(r,n)^{\bf T}$ is also finite 
and in one-to-one correspondence to nested $r$-partitions. 
\end{prop} 

{\it Proof.} Statement 
$(i)$ follows from the observation that ${\bf T}\subset {\bf T}_{r+2}$ is sufficiently generic. In particular, since $(t_1, t_2, t_3)$ are independent parameters, any ${\bf T}$-fixed framed sheaf still has to split as a direct sum of equivariant ideal sheaves. 

$(ii)$ Suppose $(E,\xi, \Phi)$ is a ${\bf T}$-fixed framed Higgs sheaf. Then 
$(E,\xi)$ is a {\bf T}-fixed framed Higgs sheaf, hence it is isomorphic to 
a framed sheaf of the form $(E_\mu, \xi_\mu)$. Moreover, the {\bf T}-fixed condition implies that the only nonzero components of the Higgs field are injections $I_{Z_{a}} \hookrightarrow I_{Z_{a+1}}$, $1\leq a \leq r$, where $Z_{r+1}$ is the empty subscheme. Using the snake lemma, each such injection is equivalent to an equivariant surjective morphism $\CO_{Z_{a}}\twoheadrightarrow 
\CO_{Z_{a+1}}$, which yields the nesting condition \eqref{eq:nestedpart}. Clearly, the converse, also holds, since the equivariant projections $\CO_{Z_{a}}\twoheadrightarrow 
\CO_{Z_{a+1}}$
are uniquely determined by the inclusions $\mu_{a+1} \subseteq \mu_a$. 

$(iii)$ Using Proposition \ref{hilbhiggspropB}, it suffices to note that all ${\bf T}$-fixed ADHM data have $J=0$. 

\hfill $\Box$ 

Note that nested ideal sheaves on surfaces occur through localization in a 
similar context \cite{GSY17b,TT17a}.

\subsection{Calabi-Yau specialization}\label{CYtorus}
The Calabi-Yau torus is by definition the two-dimensional subtorus ${\bf T}_0\hookrightarrow {\bf T}$ defined by
\[
(t_1,t_2) \mapsto (t_1, t_2, t_1^{-1}t_2^{-1}). 
\]
Geometrically this is the subtorus of ${\bf T}$ which 
preserves the natural holomorphic three-form on $\IA^2\times \IA^1 \subset \IP^2 \times \IA^1$, where $\IA^2 \subset \IP^2$ is the complement of $\Delta = \{z_0=0\}$. The goal of this section is to analyze the behavior of the fixed loci under 
this specialization.

First note the following:
\begin{rema}\label{CYtorrem}
$(i)$ 
For rank $r=1$ the ${\bf T}_0$-action on the moduli space $\CA(1,n)$, $n\geq 1$, coincides with the standard two dimensional torus action induced by the scaling action on $\IA^2$.  In particular 
the ${\bf T}_0$-fixed locus in $\CA(1,n)$ is a finite set of closed points in one-to-one correspondence with partitions $\mu$ 
of $n$. 

$(ii)$ 
For $r\geq 2$ there is a natural action of the quotient torus 
${\bf S}={\bf T}/{\bf T}_0\simeq \IC^\times$ on the fixed locus $\CA(r,n)^{\bf T}$ given by 
\be\label{eq:quotaction} 
z \times (V, B_1, B_2, I_a, J_a) \mapsto 
(V, B_1, B_2, z^{a-1}I_a, z^{1-a} J_a), \qquad 1\leq a \leq r. 
\ee
Clearly, the fixed locus of the {\bf S}-action on $\CA(r,n)^{{\bf T}_0}$ coincides with the fixed locus $\CA(r,n)^{\bf T}$. The ${\bf S}$-action on 
the ${\bf T}_0$-fixed locus will be called residual torus action. 
\end{rema}

The main goal of this section is a detailed analysis of ${\bf T}_0$-fixed locus in $\CA(r,n)$. In order to fix ideas, note the following basic facts 
on flat families of ${\bf T}_0$-fixed loci. 

Let $Z$ be an arbitrary  parameter scheme over $\IC$. A flat family of stable framed ADHM quiver representations parametrized by $Z$ is defined by a pair $(\alpha, \zeta)$ where
\begin{itemize}
\item $\alpha =(V, V_\Box, B_1, B_2, I, J)$ is a locally free ADHM quiver sheaf on $Z$ with ${\rm rk}(V)=n$, 
\item $\zeta: V_\Box {\buildrel \sim \over \longto} \CO_Z^{\oplus r}$ is an isomorphism of $\CO_Z$-modules, and
\item the restriction of the pair $(\alpha, \zeta)$ to any point in $Z$ 
is a stable framed representation of the ADHM quiver over the corresponding residual field.
\end{itemize} 
For each $1\leq a \leq r$ let  $I_a, J_a$ denote the components of the morphisms $I\circ \zeta^{-1}: \CO_Z^{\oplus r} \to V$ and $\zeta\circ J : V \to \CO_Z^{\oplus r}$ respectively. Let also $GL(V)$ denote the principal 
$GL(n,\IC)$-bundle associated to $V$ on $Z$. 

A flat family of ${\bf T}_0$-fixed framed stable ADHM quiver representations is defined by the data $(\alpha, \zeta, \eta)$ where $(\alpha, \zeta)$ are as above, and $\eta: {\bf T}_0\times Z \to GL(V)$ is a morphism of groups over $Z$ such that 
\be\label{eq:CYfixedA} 
\bal 
t_i B_i & = \eta(t_1,t_2)B_i \eta(t_1,t_2)^{-1} \\
t_1^{1-a}t_2^{1-a} I_a &  = \eta(t_1,t_2)^{-1}I_a \\
t_1^{a} t_2^{a} J_a & = J_a \eta(t_1, t_2).\\
\eal 
\ee
for any morphisms $t_1, t_2:Z \to {\bf T}_0$. 
In particular, $\eta$ determines 
a ${\bf T}_0$-equivariant structure on the underlying locally free $\CO_Z$-module $V$. Let
\[
V = \bigoplus_{(i,j)\in \IZ^2} V(i,j) 
\]
denote the resulting character decomposition.
Then conditions \eqref{eq:CYfixedA} imply that 
\be\label{eq:CYfixedB} 
{\rm Im}(I_a) \subseteq V(a-1,a-1), \qquad 
J_a\big|_{V(i,j)}=0\quad {\rm for\ all}\quad (i,j)\neq (a,a), 
\ee
for all $1\leq a \leq r$. Moreover, only the following components  
\be\label{eq:CYfixedC}
B_1(i,j): V(i,j) \to V(i-1,j), \qquad  B_2(i,j):V(i,j)\to V(i,j-1)
\ee
of $B_1, B_2$ 
are allowed to be nonzero. All other components must vanish identically. 
Hence the ADHM relation reduces to 
\be\label{eq:CYfixedD}
B_1(a,a-1) B_2(a,a) - B_2(a-1,a) B_1(a,a) + I_a J_a =0 
\ee
for all $1\leq a\leq r$.

The first structure result is the following. 
\begin{lemm}\label{complemm} 
Let $Y$ be a connected component of the fixed locus $\CA(r,n)^{{\bf T}_0}$. 
Let $y\in Y$ be an arbitrary closed point and let $f_y: \IA^1\setminus \{0\}\to Y$ be the ${\bf S}$-orbit through $y$. Then $f_y$ extends uniquely to a morphism ${\bar f}_y : \IA^1\to Y$. In particular the ${\bf S}$-fixed locus $Y^{\bf S}$ contains at least one point. 
\end{lemm} 

{\it Proof}. If the ${\bf S}$-action on $Y$ is trivial, $Y$ is a connected component of the fixed locus $\CA(r,n)^{\bf T}$, which is finite. Hence the claim is obvious. 

Suppose the ${\bf S}$-action on $Y$ is not trivial. 
If $y \in Y^{\bf S}$, which is again finite, the claim follows. 
Therefore it suffices to consider the case where  $y$ is not ${\bf S}$-fixed. Then the given orbit is nontrivial and 
the claim is proven by analogy to \cite[Thm. 3.7]{Lect_inst}. 
The moduli space $\CA(r, n)$ is a GIT quotient $\CR(r,n)//GL(n,\IC)$ where
\[ 
\CR(r,n)\subset {\rm End}(\IC^n)^{\oplus 2}\oplus {\rm Hom}(\IC^n, \IC^r) 
\oplus {\rm Hom}(\IC^r, \IC^n) 
\]
is the zero locus $[B_1,B_2]+IJ=0$. By construction, there is a projective morphism $\pi:\CA(r,n)\to \CA_{0}(r,n)$ 
to the affine algebro-geometric quotient, which is the spectrum of the 
ring of $GL(n,\IC)$-invariant polynomials on $\CR(r,n)$. 
As in the proof of \cite[Thm. 3.7]{Lect_inst}, the latter is generated by the following types of functions 
\[
{\rm Tr}_{\IC^n}(m(B_1,B_2)), \qquad J_b m(B_1,B_2) I_a, \qquad 1\leq a,b\leq r,
\]
where $m(B_1,B_2)$ denotes an arbitrary  monomial in ${\rm End}(\IC^n)$.
The  functions of the first type have weights 
\[
({\rm deg}_{B_1}m(B_1,B_2), \ {\rm deg}_{B_2}m(B_1,B_2),\ 0) 
\]
under the action of ${\bf T}=(\IC^{\times})^{\times 3}$. 
The functions of the second type have weights 
\[
({\rm deg}_{B_1}m(B_1,B_2)+1, \ {\rm deg}_{B_2}m(B_1,B_2)+1,\ a-b). 
\]
If all the above invariant functions have trivial restriction to $Y$, it follows that $Y$ is a closed subscheme of $\pi^{-1}(0)$, hence it is projective. Then the claim follows. 

Suppose this is not the case. 
Then note that the  fixed locus conditions \eqref{eq:CYfixedB} and \eqref{eq:CYfixedC} imply that $f_y^*\pi^*{\rm Tr}_{\IC^n}(m(B_1,B_2))=0$. Therefore the exists a pair $(a,b)$, $1\leq a,b\leq r$ such that 
$f_y^*\pi^*(J_b(m(B_1,B_2)I_a)$ is nonzero. Then the fixed locus conditions 
imply that 
\[
{\rm deg}_{B_1}m(B_1,B_2) = {\rm deg}_{B_2}m(B_1,B_2) = a-b-1.
\]
and the ${\bf S}$-weight of $f_y^*\pi^*(J_b(m(B_1,B_2)I_a)=0$ is $a-b \geq 1$. This further implies that the morphism 
$\pi\circ f_y: \IA^1 \setminus\{0\} \to \CA_0(r,n)$ extends uniquely to $\IA^1$. Since $\pi$ is projective, $f_y$ can be also extended to a morphism 
$\IA^1 \to \CA(r,n)$. Since $Y$ is a closed subscheme of $\CA(r,n)$, the claim follows. 

\hfill $\Box$

In order to formulate the next result note that, omitting the rigidifying 
isomorphism $\zeta: V_\Box \to \CO_Z^{\oplus r}$, coherent ADHM quiver sheaves on 
$Z$ form an abelian category $\CA_Z$. Let also ${\bf e}_a$, $1\leq a \leq r$, be the standard basis vectors in  
$\IC^r$ and let 
\be\label{eq:inftyfiltrA} 
0\subset W_1 \subset \cdots W_2\subset \cdots \subset W_r  
\ee
be the filtration defined by 
\[
W_a = \IC \langle {\bf e}_1, \ldots, {\bf e}_a\rangle, \qquad 1\leq a\leq r,
\]
all inclusions being canonical. 
Then one has: 

\begin{lemm}\label{CYfixedlemmA} 
Let $\alpha$ be a flat family of ${\bf T}_0$-fixed stable framed ADHM quiver representations with dimension vector $(r,n)$ parametrized by a connected scheme $Z$. 
Suppose furthermore that $J$ is identically zero. 
Let 
\be\label{eq:inftyfiltrB}
V_\Box^{(a)} = \zeta^{-1}\left(W_a \otimes \CO_Z\right), \qquad 1\leq a \leq r, 
\ee
be the filtration induced by \eqref{eq:inftyfiltrA}. 
Then there exist an unique $r$-partition $\mu$ of $n$ 
 and a filtration $\alpha_\bullet$ 
\be\label{eq:CYfixedfiltrA} 
0 \subset \alpha_{r} \subset \cdots \subset \alpha_1=\alpha 
\ee
in the abelian category $\CA_Z$ such that 

$(i)$ The restriction of $\alpha_\bullet$ to the node $\Box$ 
coincides with the filtration \eqref{eq:inftyfiltrB}, and 

$(ii)$ Each succesive quotient ${\bar \alpha}_a$, $1\leq a\leq r$, 
is isomorphic to the stable framed ADHM representation corresponding to  
${\bf T}_0$-fixed point $\alpha_{\mu_a}\in \CA(1, |\mu_a|)^{{\bf T}_0}$.
\end{lemm} 

{\it Proof}. 
As observed above, $V$ has a ${\bf T}_0$-character decomposition 
\[
V = \bigoplus_{(i,j)} V(i,j)
\]
such that the ADHM data satisfy conditions \eqref{eq:CYfixedB} and 
\eqref{eq:CYfixedC}. Moreover, stability implies that $V(i,j)$ is identically zero if $i\geq r$ or $j\geq r$. 
Let 
\be\label{eq:filtrB}
0\subset V^{(1)} \subset \cdots\subset V^{(r)}=V
\ee
be the filtration defined by 
\[ 
V^{(a)} = \bigoplus_{\substack{i\leq a-1\\ j \leq a-1}} V(i,j) 
\]
where all the inclusions are canonical. 
Then conditions \eqref{eq:CYfixedB}, \eqref{eq:CYfixedC} imply that 
\be\label{eq:fitrC} 
B_i(V^{(a)})\subseteq V^{(a)}, \qquad 1\leq a\leq r, \qquad 1\leq i \leq 2,
\ee
\be\label{eq:filtrC}
I(V^{(a)}_\Box) \subseteq V^{(a)}, \qquad 1\leq a \leq r,
\ee
In addition, since $J$ is assumed identically zero, the ADHM relation restricts to 
\[
B_1(a,a-1) B_2(a,a) - B_2(a-1,a) B_1(a,a) =0. 
\]
Let ${\overline V}_a= V^{(a)}/V^{(a-1)}$, respectively 
${\overline V}_{\Box,a}= V^{(a)}_\Box/V^{(a-1)}_\Box\simeq \IC\langle {\bf e}_a\rangle$, 
$1\leq a\leq r$ where $V^{(0)}=0$ and $V_\Box^{(0)}=0$.   Equations \eqref{eq:fitrC} imply that there are naturally induced maps 
\[ 
{\overline B}_{i,a} :{\overline V}_a\to {\overline V}_a, \qquad 
{\overline I}_a : {\overline V}_{\Box,a} \to {\overline V}_a, \qquad 
1\leq a \leq r,\qquad 1\leq i \leq 2
\]
such that $[{\overline B}_{1,a}, {\overline B}_{2,a}]=0$. At the same time, the ADHM stability condition on $\alpha$ 
implies that the data ${\overline \alpha}_a = ({\overline V}_a, {\overline B}_{i,a}, {\overline I}_a)$ is a framed stable rank one ADHM representation over $Z$
for each $1\leq a \leq r$.  
Finally, by construction, each ${\overline \alpha}_a$ is a ${\bf T}_0$-fixed point in the moduli space of rank one ADHM representations
$\CA(1,n_a)^{{\bf T}_0}$ where $n_a= {\rm dim}({\overline V}_a)$. 
Since $Z$ is connected and $\CA(1,n_a)^{{\bf T}_0}$ is a finite set of closed points indexed by partitions of $n_a$, the claim follows. 

\hfill $\Box$

In order to formulate a useful consequence of Lemmas \ref{complemm} and \ref{CYfixedlemmA} 
recall the projective morphism to the affine geometric quotient, 
$\pi:\CA(r,n)\to \CA_0(r,n)$, used in the proof 
of Lemma \ref{complemm}. Then one has: 
\begin{coro}\label{CYfixedcorAA}
Let $Y\subset \CA(r,n)^{{\bf T}_0}$ be a connected component of the ${\bf T}_0$-fixed locus such that ${\bf S}$-action on $Y$ is nontrivial. Then $Y$ is not projective over $\IC$.  
\end{coro}

{\it Proof}. A priori $Y$ is smooth quasi-projective. 
Suppose $Y$ is projective. Then the fixed locus $Y^{\bf S}$ is nonempty
and consists of finitely many points in $\CA(r,n)^{\bf T}$. Moreover, all {\bf T}-fixed points are mapped to $0$ by $\pi$. Since $\CA_0(r,n)$ is affine and $Y$ is connected, it follows that $Y$ must be contained as a closed subscheme in the fiber $\pi^{-1}(0)$. 

Let $\alpha_Y$ denote the restriction of the universal family of stable framed ADHM data to $Y$. 
Since $\pi(Y)=\{0\}$ all the generators of polynomial ring $\Gamma(\CA_0(r,n))$, in particular
\[
J_b m(B_1, B_2) I_a, \qquad 1\leq a, b \leq r, 
\]
have trivial pull-back to $Y$. 
Then the ADHM 
stability condition implies that the family $\alpha_Y$ has $J=0$. Therefore, as shown in Lemma \ref{CYfixedlemmA}, there is a filtration 
of the form 
\eqref{eq:CYfixedfiltrA}. 
However, given the action \eqref{eq:quotaction} of ${\bf S}$, a point $y\in Y$ is fixed by $S$ if and only if the induced filtration on $\alpha_Y|_y$ is split. 
Therefore $Y^{\bf S}$ consists of the unique closed point $\alpha_\mu \in 
\CA(r,n)^{\bf T}$, where $\mu$ is the $r$-partition determined by $\alpha_Y$ as in Lemma \ref{CYfixedlemmA}. 
This contradicts the assumption that $Y$ is smooth projective and the ${\bf S}$-action on $Y$ is nontrivial. 

\hfill $\Box$

Next, 
a ${\bf T}$-fixed point $\alpha_\mu\in \CA(r,n)^{\bf T}$ will be called ${\bf T}_0$-isolated if and only if $\{\alpha_\mu\}$ is a zero dimensional connected component of the fixed locus 
$\CA(r,n)^{{\bf T}_0}$. Then it will be shown below that $\alpha_\mu$ is 
${\bf T}_0$-isolated for all nested partitions $\mu$. 
The proof will use the ${\bf T}_{r+2}$-character decomposition of the tangent space to the fixed point ${\alpha_\mu}\in \CA(r,n)^{{\bf T}_{r+2}}$. As in Section 3.2 of \cite{Cherednik_W}, let $q,t,\chi_a: {\bf T}_{r+2}\to \IC$ be the characters defined by 
\be\label{eq:torchar}
q(t_1, t_2,u) = t_1^{-1}, \qquad 
t(t_1, t_2,u) = t_2^{-1}, \qquad 
\chi_a(t_1, t_2, u) = u_a^{-1}, \quad 1\leq a \leq r.
\ee

Given a partition $\nu \subset \IZ^2$, for any box $s= (i(s), j(s))\in \IZ^2$ one defines 
\begin{itemize} 
\item $\ell_\nu(s) = \nu_{i(s)} - j(s)-1$
\item $a_\nu(s) = \nu^{t}_{j(s)}  -i(s) -1$ 
\end{itemize} 
where $\nu_{i}$ denotes the number of boxes on the $i$-th column 
of $\nu$ and $\nu^t_j$ denotes the number of boxes on  the $j$-th 
column of $\nu^t$, which is the same as the number of boxes on the $j$-th row of $\nu$. Note that the above numbers are negative if $s\notin \nu$.

Then, as shown in \cite[Thm. 2.11]{instcountA}, the explicit formula for the character decomposition of $T_\mu$ is 
\be\label{eq:equivtangentA} 
{\rm ch}_{{\bf T}_{r+2}}\, T_\mu = \sum_{b,c=1}^r \sum_{s\in \mu_b} \chi_b\chi_c^{-1} t^{\ell_{\mu_c}(s)}q^{-{a}_{\mu_b}(s)-1} 
+ \sum_{b,c=1}^r \sum_{s\in \mu_c} \chi_b\chi_c^{-1} t^{-\ell_{\mu_b}(s)-1}q^{a_{\mu_c}(s)}.
\ee 

Abusing notation, below  let $q,t$ denote  the restrictions of the characters $q,t$ to ${\bf T}$. Let also $\sigma:{\bf T}\to \IC$ denote the character $\sigma(t_1,t_2,t_3) = t_1t_2t_3$. For any $r$-partition $\mu$ of $n$, let 
\be\label{eq:Tcharmu}
\ch_{\bf T}(T_\mu) = \sum_{i,j,k\in \IZ} c_{i,j,k}(\mu) q^i t^j \sigma^k.
\ee
be the ${\bf T}$-character decomposition of $T_\mu$. 
Moreover, let ${\sf S}_1(\mu)$ denote the set of triples $(b,c,s)$ defined by the following conditions:
\be\label{eq:Scharone} 
1\leq c<b \leq r,\qquad 
s\in \mu_b\setminus \mu_c,\qquad b-c +\ell_{\mu_c}(s) =0,\qquad 
 b-c -a_{\mu_b(s)}-1=0.
 \ee
Let also ${\sf S}_2({\mu})$ denote the set of triples $(b,c,s)$ defined by:
\be\label{eq:Schartwo} 
1\leq b<c\leq r,\qquad 
s\in \mu_c\setminus \mu_b,\qquad b-c -\ell_{\mu_b}(s)-1 =0,\qquad b-c+a_{\mu_c(s)}=0.
\ee
Then the following holds:
\begin{lemm}\label{CYfixedlemmB} 
Let $\mu$ be an $r$-partition of $n$. Then, using the notation in \eqref{eq:Tcharmu}, there is an identity 
\be\label{eq:equivtangentC} 
\bal
\sum_{k\in \IZ} c_{0,0,k}(\mu) \sigma^k  = 
  \sum_{(b,c,s)\in {\sf S}_1(\mu)} 
\sigma^{b-c} \  + \
\sum_{(b,c,s)\in {\sf S}_2(\mu)} \sigma^{b-c}.
\eal
\ee
in the character ring of ${\bf T}$. 
\end{lemm}

{\it Proof}. 
Specializing \eqref{eq:equivtangentA} to ${\bf T}\subset {\bf T}_{r+2}$, one obtains 
\be\label{eq:equivtangentB} 
{\rm ch}_{\bf T}\, T_\mu = \sum_{b,c=1}^r \sum_{s\in \mu_b} \sigma^{b-c} t^{b-c+\ell_{\mu_c}(s)}q^{b-c-{a}_{\mu_b}(s)-1} 
+ \sum_{b,c=1}^r \sum_{s\in \mu_c} \sigma^{b-c} t^{b-c-\ell_{\mu_b}(s)-1}q^{b-c+a_{\mu_c}(s)}.
\ee 
Therefore such eigenvectors are in one-to-one correspondence to triples 
$(b,c,s)$, $1\leq b,c \leq r$, satisfying one of the following two conditions:
\be\label{eq:equivtangentD}
s\in \mu_b,\qquad b-c +\ell_{\mu_c}(s) =0, \qquad b-c -a_{\mu_b}(s)-1=0,
\ee
\be\label{eq:equivtangentE}
s\in \mu_c, \qquad b-c -\ell_{\mu_b}(s)-1 =0, \qquad b-c +a_{\mu_c}(s)=0.
\ee

For any triple $(b,c,s)$ as in \eqref{eq:equivtangentD} one has $a_{\mu_b}(s) \geq 0$ since $s\in \mu_b$. Hence $b-c \geq 1$ from the second equation in \eqref{eq:equivtangentD}, and $\ell_{\mu_c}(s) = c-b \leq -1$ from the first equation in \eqref{eq:equivtangentD}. In particular, 
$s\in \mu_b \setminus \mu_c$, which must be necessarily non-empty. This yields 
conditions \eqref{eq:Scharone}. 

For any triple $(b,c,s)$ as in \eqref{eq:equivtangentE},  
one has $a_{\mu_c}(s) \geq 0$ since $s\in \mu_c$. Hence $b-c\leq 0$ from the  second equation in \eqref{eq:equivtangentE}. If $b-c=0$, then $\mu_b=\mu_c$ and $\ell_{\mu_b}(s) \geq 0$. This contradicts the
first equation in  \eqref{eq:equivtangentE}. Therefore $ b-c\leq -1$, and the 
 first equation in  \eqref{eq:equivtangentE} implies $\ell_{\mu_b}(s) 
 = b-c-1\leq -2$. Hence $\mu_c\setminus \mu_b$ must be nonempty and $s\in \mu_c \setminus \mu_b$. 
  
In conclusion, equation \eqref{eq:equivtangentC} follows from \eqref{eq:equivtangentB}. 
 
 \hfill $\Box$

In order to formulate the next result, note that an $r$-partition $\mu$ of $n$ will be said to be contained, $\mu\subset \lambda$, into an $r$-partition of $n+1$ if there exists 
$1\leq a \leq r$ such that $\mu_b = \lambda_b$ for all $1\leq b\leq r$, $b\neq a$, while $\mu_a\subset \lambda_a$ as Young diagrams. In particular, $\lambda_a \setminus 
\mu_a$ consists of a single box. Then one has the following consequence of 
Lemma \ref{CYfixedlemmB}

\begin{coro}\label{CYfixedcorB}
$(i)$ Let $\mu$ be a nested $r$-partition of $n$ i.e. 
\[
\mu_r \subseteq \mu_{r-1} \subseteq \cdots \subseteq\mu_1.
\]
Then identity \eqref{eq:equivtangentC} reduces to 
\be\label{eq:CYnestedA} 
\sum_{k\in \IZ} c_{0,0,k}(\mu) \sigma^k  = 0.
\ee 

$(ii)$ Let $\lambda$ be a nested $r$-partition of $n$ and let $\mu \supset \lambda$
be an $r$-partition of $n+1$ such that $\mu$ is not nested. Then 
identity \eqref{eq:equivtangentC} reduces to 
\be\label{eq:CYnestedB} 
\sum_{k\in \IZ} c_{0,0,k}(\mu) \sigma^k  = \sigma.
\ee

$(iii)$ Let $\lambda$ be a nested $r$-partition of $n+1$ and let $\mu\subset \lambda$
be an $r$-partition of $n$ such that $\mu$ is not nested. 
Then identity \eqref{eq:equivtangentC} reduces to 
\be\label{eq:CYnestedA} 
\sum_{k\in \IZ} c_{0,0,k}(\mu) \sigma^k  = \sigma.
\ee

Furthermore, in each case suppose $Y$ is the unique connected component of 
the ${\bf T}_0$-fixed locus containing the {\bf T}-fixed point $\alpha_\mu$. Then $Y=\{\alpha_\mu\}$ in case $(i)$, while in cases 
$(ii)$ and $(iii)$ one has an isomorphism $Y\simeq \IA^1$ mapping $\alpha_\mu$ to $0\in \IA^1$. 
\end{coro} 

{\it Proof}. $(i)$ Since $\mu$ is nested, one has ${\sf S}_1(\mu) = {\sf S}_2(\mu) =\emptyset$. 

$(ii)$ Under the stated assumptions, there exists exactly one index $2\leq b \leq r$ such that $\lambda_c=\mu_c$ for all $1\leq c \leq r$, $c\neq b$ while $\lambda_b \setminus \mu_b$ consists of a single box $s$. Since $\mu$ is nested one has ${\sf S}_1(\lambda)=\{(b,b-1,s)\}$. At the same time note that for all triples in ${\sf S}_2(\mu)$ one has $\ell_{\mu_b}\leq -2$. 
Since $\mu$ is nested, this implies ${\sf S}_2(\mu)=\emptyset$. 

The proof of $(iii)$ is analogous to that of $(ii)$. 

The last statement then follows from Corollary \ref{CYfixedcorAA}.

\hfill $\Box$

Finally, note the following. 
\begin{lemm}\label{CYfixedlemmD} 
The ${\bf T}_0$-fixed locus $Hilb(r,n)^{{\bf T}_0}$ coincides with 
the ${\bf T}$-fixed locus $Hilb(r,n)^{{\bf T}}$. Moreover, any ${\bf T}_0$-connected component $Y$ intersecting $Hilb(r,n)$ nontrivially must be 
a single ${\bf T}_0$-isolated point belonging to $Hilb(r,n)$. 
\end{lemm}

{\it Proof.} Clearly, there is a closed embedding $Hilb(r,n)^{{\bf T}}\subset Hilb(r,n)^{{\bf T}_0}$.
Let $Z$ be a nonempty connected component of $Hilb(r,n)^{{\bf T}_0}$. 
Since $Hilb(r,n)$ is a closed subscheme of $\CA(r,n)$, there exists a unique connected component $Y$ of $\CA(r,n)^{{\bf T}_0}$ such that $Z$ is 
the scheme theoretic intersection $Y \times_{\CA(r,n)}Hilb(r,n)$. In particular $Z$ is a closed subscheme of $Y$. 
Moreover, $Z$ is also naturally preserved by the residual {\bf S}-action. 

If the ${\bf S}$-action on $Z$ is trivial, then $Z$ is a finite set of 
${\bf T}$-fixed closed points in $\CA(r,n)$ which also belong to $Hilb(r,n)$.  Then Corollary \ref{CYfixedcorB} shows that each such point is ${\bf T}_0$-isolated. Hence $Z$ coincides with $Y$, which is assumed connected. Hence $Y$ must reduce to a single {\bf T}-fixed point belonging to $Hilb(r,n)$. 

Suppose the ${\bf S}$-action on $Z$ in not trivial.  Since $Z$ is closed in $Y$,
Lemma \ref{complemm} implies that $Z$ contains at least one  {\bf S}-fixed point. Then the same argument as in the previous paragraph shows that $Y$ must be a single ${\bf T}$-fixed point 
belonging to $Hilb(r,n)$. 

\hfill $\Box$

\section{Hecke transformations}\label{convsect}

The goal of this section is to review the construction of Hecke transformations used in \cite{Cherednik_W}, at the same time proving certain 
properties of fixed points in the correspondence variety. 

\subsection{The ADHM correspondence variety}\label{correspvar}

As in Section 3.3 of \cite{Cherednik_W}, let $\CA(r,n,n+1)\subset \CA(r,n)
\times \CA(r,n+1)$ denote the Hecke correspondence parameterizing elementary modifications of framed torsion sheaves on $\IP^2$ supported at a single closed point $p\in \IA^2$. Let $\CA_c(r,n,n+1)\subset \CA(r,n,n+1)$ denote the closed subvariety where $p=0$. 
According to 
\cite[Proposition 3.1]{Cherednik_W}, the following holds.
\begin{prop}\label{corresprop}
$(1)$ 
The correspondence variety $\CA(r,n,n+1)$ is a smooth quasi-projective variety of dimension
$2rn+r+1$. 

$(2)$ The natural morphism $\CA(r,n,n+1)\to \CA(r,n) \times \CA(r,n+1)$ is a closed embedding, and the restriction of the projection $\pi_2: \CA(r,n) \times \CA(r,n+1)\to \CA(r, n+1)$ to $\CA(r,n,n+1)$ is proper. 

$(3)$ The restriction of the projection $\pi_1: \CA(r,n) \times \CA(r,n+1)\to \CA(r, n)$ to $\CA_c(r,n,n+1)$ is also proper. 

$(4)$ The correspondence variety is preserved by the 
${\bf T}_{r+2}$-action on the product and the fixed locus 
$\CA(r, n, n+1)^{{\bf T}_{r+2}}$ is a finite set of closed points in one-to-one correspondence with pairs of $r$-partitions $(\mu, \lambda)$ such that 
\be\label{eq:convfixedA}
|\mu|=n, \qquad |\lambda| = n+1, \qquad \mu\subset \lambda.
\ee
Written more explicitly, the last condition in \eqref{eq:convfixedA} 
states that there exists $1\leq b \leq r$ such that  $\mu_a=\lambda_a$ for all $1\leq a\leq r$, $a\neq b$, while $\mu_b \subset \lambda_b$ and $\lambda_b \setminus \mu_b$ consists of a single box $s$. 
\end{prop}

\begin{rema}\label{Tcorresp} 
Note that Proposition \ref{frquivpropC}.$i$ implies that the fixed locus 
$\CA(r,n,n+1)^{\bf T}$ coincides with the ${\bf T}_{r+2}$ fixed locus. 
Therefore {\bf T}-equivariant Hecke transformations are defined in complete analogy to \cite{Cherednik_W}. This is briefly reviewed below. 
\end{rema}

Let $\gamma : 
\CA(r,n,n+1)\hookrightarrow \CA(r,n)\times \CA(r,n+1)$ 
denote the natural closed embedding, which is clearly ${\bf T}$-equivariant. 
As in Section 2.2 of \cite{Cherednik_W}, any
equivariant Borel-Moore homology class $z\in H^{{\bf T}}(\CA(r,n,n+1))$ 
determines a Hecke transformation, 
\be\label{eq:HeckemapA}
u_{z,n}^{+}:H^{{\bf T}}(\CA(r,n))_{K} \to H^{{\bf T}}(\CA(r, n+1))_{K}, \qquad 
u_{z,n}^{+}(x)=\pi_{2*}((\gamma_*z)\cdot \pi_1^*x).
\ee
The pull-back and push-forward maps in equivariant Borel-Moore homology are well defined since $\gamma$ is a closed embedding and 
the restriction of $\pi_2$ to $\CA(r,n,n+1)$ is proper. Although $\pi_1$ is 
not proper, its restriction to the ${\bf T}$-fixed locus is. Therefore, using the localization theorem one can define downward Hecke transformations on equivariant Borel-Moore homology 
\be\label{eq:HeckemapB}
u_{z,n}^{-}:H^{{\bf T}}(\CA(r,n+1))_{} \to H^{{\bf T}}(\CA(r, n))_{K}, \qquad 
u_{z,n}^{-}(x)=\pi_{1*}((\gamma_*z)\cdot \pi_2^*x).
\ee
Moreover, the following explicit formulas hold by a straightforward application of the localization theorem for the correspondence variety.
These formulas have been used for example in Appendix C of \cite{Cherednik_W}, hence the proof will be omitted. 
\begin{lemm}\label{HeckelemmA} 
$(i)$ For any $r$-partition $\mu$, 
\be\label{eq:HeckemapAA}
u_{z,n}^+([\alpha_\mu]) = \sum_{\substack{\lambda \in \calP_{r,n+1}\\  
\lambda \supset \mu\\ }}
e_{{\bf T}}(T_\mu) e_{{\bf T}}(T_{\mu,\lambda})^{-1}
z_{\mu, \lambda}[\alpha_\lambda],
\ee
where $z_{\nu, \lambda}$ denotes the restriction of $z$ to the 
fixed point $(\alpha_\nu, \alpha_\lambda)$. 

$(ii)$ In the opposite direction, 
\be\label{eq:HeckemapBB}
u_{z,n+1}^-([\alpha_\lambda]) = \sum_{\substack{\mu\in \calP_{r,n}\\ \mu\subset \lambda \\ }}
e_{{\bf T}}(T_\lambda) e_{{\bf T}}(T_{\mu,\lambda})^{-1}
z_{\mu, \lambda}[\alpha_\mu].
\ee
\end{lemm}

In order to construct analogous transformations for the Hilbert scheme 
$Hilb(r,n)\subset \CA(r,n)$ one needs first some structure results for the fixed locus $\CA(r,n)^{{\bf T}_0}$, which are proven in the next section.

\subsection{The fixed locus of the Calabi-Yau 
torus action}\label{correspfixed} 

For any fixed point $\alpha_\mu \in \CA(r,n)^{{\bf T}_{r+2}}$ let $V_\mu$ denote the 
underlying vector space, which carries a linear ${\bf T}_{r+2}$-action such that 
\[ 
{\rm ch}_{{\bf T}_{r+2}}\, V_\mu = \sum_{b=1}^r \chi_b^{-1} \sum_{s\in \mu_b} q^{i(s)-1} t^{j(s)-1}. 
\]
As in \eqref{eq:torchar}, the characters $q,t, \chi_a$, $1\leq a \leq r$,  are defined by 
\[ 
q(t_1, t_2, u) = t_1^{-1}, \qquad 
t(t_1, t_2, u) =t_2^{-1}, \qquad 
\chi_a(t_1, t_2, u) = u_a^{-1}.
\]
Using the notation in Section 3.4 of \cite{Cherednik_W}, set $\tau_\mu = {\rm ch}_{{\bf T}_{r+2}}\, V_\mu$ and $w= \sum_{b=1}^r \chi_b^{-1}$ and let 
 $T_{\mu}$ denote the tangent space to $\CA(r,n)$ at the fixed point $\alpha_\mu$. 
As shown in  \cite[Thm. 2.11]{instcountA}, the character of the equivariant tangent space $T_{\mu}$ is given by 
\be\label{eq:equivtangentG} 
{\rm ch}_{{\bf T}_{r+2}}\, T_{\mu} = -(1-q^{-1})(1-t^{-1}){\tau}_\mu \tau_\mu^{\vee} + \tau_\mu w^{\vee} + q^{-1}t^{-1} \tau_\mu^{\vee}w~.
\ee

For any pair $(\mu,\lambda)$ with $|\mu|=n$, $|\lambda|=n+1$ and $\mu \subset \lambda$ let $N_{\mu, \lambda}$ denote the fiber of the equivariant normal bundle to  $\CA(r,n,n+1)$ in the product $\CA(r,n)\times 
\CA(r, n+1)$ at the fixed point $(\alpha_\mu, \alpha_\lambda)$.  Then the character of the ${\bf T}_{r+2}$-action on $N_{\mu, \lambda}$ is given by equation (3.11) in \cite{Cherednik_W}, which reads
\be\label{eq:equivnormalA}
{\rm ch}_{{\bf T}_{r+2}}\, N_{\mu, \lambda} = -(1-q^{-1})(1-t^{-1})
{\tau}_\mu\tau_\lambda^\vee + \tau_\mu w^{\vee} + q^{-1}t^{-1} \tau_\lambda^{\vee} w
- q^{-1} t^{-1}. 
 \ee
 
 As observed in Remark  \ref{Tcorresp}, the fixed locus in the correspondence variety remains unchanged under specialization to ${\bf T}\subset {\bf T}_{r+2}$.  Moreover, the character decomposition of the tangent space 
to a fixed point $(\alpha_\mu, \alpha_\lambda)$ is obtained by straightforward specialization. As in Lemma \ref{CYfixedlemmB},  let $\sigma: {\bf T}\to \IC^\times$ denote the character $\sigma(t_1,t_2, t_3) = t_1t_2t_3$. Then  
\[
\chi_b|_{\bf T} = (qt\sigma)^{b-1}, \qquad 1\leq b \leq r.
\]
Let $T_{\mu,\lambda}$ denote the ${\bf T}$-equivariant tangent space 
to the fixed point $(\alpha_\mu, \alpha_\lambda)\in \CA(r,n,n+1)^{\bf T}$. 
Let 
\[
{\rm ch}_{\bf T}( T_{\mu, \lambda}) = \sum_{i,j,k\in \IZ} c_{i,j,k} q^i t^j \sigma^k 
\]
be the ${\bf T}$-character decomposition of $T_{\mu, \lambda}$.  
As in Section \ref{CYtorus}, a  fixed point $(\alpha_\mu, \alpha_\lambda)\in \CA(r,n,n+1)^{\bf T}$ will be called 
${\bf T}_0$-isolated if and only if it is a connected component of the 
${\bf T}_0$-fixed locus. 
Then one has: 

\begin{lemm}\label{correspfixedlemmAB}
$(i)$ Suppose at least one $r$-partition of the pair $(\mu, \lambda)$ is nested. Then $c_{0,0,k}=0$ for all $k \geq 0$. In particular the fixed 
point $(\alpha_\mu, \alpha_\nu)$ is  ${\bf T}_0$-isolated. 

$(ii)$ Suppose $\mu$ and $\lambda$ are not nested, $n\geq 1$, and there is a nested $r$-partition $\nu$ of $n-1$ such that $\nu \subset \mu$. Then 
$c_{0,0,1}=1$ and $c_{0,0,k}=0$ for all $k\in \IZ$, $k\neq 1$. 

\end{lemm}

{\it Proof.} 
Note that 
\[ 
{\rm ch}_{{\bf T}_{r+2}}\, T_{\mu, \lambda}  = 
{\rm ch}_{{\bf T}_{r+2}}\, T_\mu + {\rm ch}_{{\bf T}_{r+2}}\, T_\lambda - 
{\rm ch}_{{\bf T}_{r+2}}\, N_{\mu,\lambda}.
\]
Using equations, \eqref{eq:equivtangentG} and \eqref{eq:equivnormalA} 
one has
\be \label{notrivial1}
\bal
&
{\rm ch}_{{\bf T}_{r+2}}\, T_\mu + {\rm ch}_{{\bf T}_{r+2}}\, T_\lambda - 
{\rm ch}_{{\bf T}_{r+2}}\, N_{\mu,\lambda}
 = \\
& -(1-q^{-1})(1-t^{-1}){\tau}_\mu \tau_\mu^{\vee} + \tau_\mu w^{\vee} + q^{-1}t^{-1} \tau_\mu^{\vee}w \\
& -(1-q^{-1})(1-t^{-1}){\tau}_\lambda \tau_\lambda^{\vee} + \tau_\lambda w^{\vee} + q^{-1}t^{-1} \tau_\lambda^{\vee}w \\
& + (1-q^{-1})(1-t^{-1})
{\tau}_\mu\tau_\lambda^\vee - \tau_\mu w^{\vee} - q^{-1}t^{-1} \tau_\lambda^{\vee} w + q^{-1} t^{-1} = \\
&  -(1-q^{-1})(1-t^{-1})(\tau_\lambda - \tau_\mu) (\tau_\lambda- \tau_\mu)^{\vee} - (1-q^{-1})(1-t^{-1}) \tau_\lambda \tau_\mu^{\vee}\\
&\ + \tau_\lambda w^{\vee} + q^{-1}t^{-1} \tau_\mu^{\vee}w + q^{-1} t^{-1} \ .
\eal
\ee
By assumption $\tau_\lambda - \tau_\mu$ is a one dimensional representation of ${\bf T}_{r+2}$ since $\mu \subset \lambda$, and $\lambda \setminus \mu$ consists of a single box. Hence the 
the right hand side of \eqref{notrivial1} reduces to 
\be \label{notrivial2}
\bal
&
{\rm ch}_{{\bf T}_{r+2}}\, T_\mu + {\rm ch}_{{\bf T}_{r+2}}\, T_\lambda - 
{\rm ch}_{{\bf T}_{r+2}}\, N_{\mu,\lambda}
= \\
& -1+q^{-1}+t^{-1} - (1-q^{-1})(1-t^{-1}) \tau_\lambda \tau_\mu^{\vee} + \tau_\lambda w^{\vee} + q^{-1}t^{-1} \tau_\mu^{\vee}w.\\
\eal
\ee
Let 
\[
E_{\mu, \lambda} = -(1-q^{-1})(1-t^{-1}) \tau_\lambda \tau_\mu^{\vee} + \tau_\lambda w^{\vee} + q^{-1}t^{-1} \tau_\mu^{\vee}w.
\] 
By analogy with
\cite[Theorem 2.11]{instcountA} this expression is given by 
\be \label{notrivial3}
\bal
& E_{\mu, \lambda} = \sum_{b,c=1}^r \sum_{s\in \lambda_b} \chi_b\chi_c^{-1}
t^{\ell_{\lambda_c}(s)}q^{-a_{\mu_b}(s)-1} 
+ \sum_{b,c=1}^r \sum_{s\in \mu_c} \chi_{b}\chi_c^{-1} t^{-\ell_{\mu_b}(s)-1}q^{a_{\lambda_c}(s)}\ .
\eal
\ee  
The ${\bf T}$-specialization of $E_{\mu, \lambda}$ is
\be\label{notrivial4}
E_{\mu, \lambda}|_{\bf T}=\sum_{b,c=1}^r \sum_{s\in \lambda_b}\sigma^{b-c}t^{b-c+\ell_{\lambda_c}(s)}q^{b-c-a_{\mu_b}(s)-1} 
+ \sum_{b,c=1}^r \sum_{s\in \mu_c} \sigma^{b-c} t^{b-c-\ell_{\mu_b}(s)-1}q^{b-c+a_{\lambda_c}(s)}
\ee
Let $1\leq d\leq r$ be such that $\lambda_d \setminus \mu_d=\{u\}$ while 
$\lambda_b = \mu_b$ for $b\neq d$. The terms with $b=d$ and $s=u$ in the
in the right hand side of the above equation are 
\[
\sum_{c=1}^r \sigma^{d-c} t^{d-c+\ell_{\lambda_c}(u)}q^{d-c-a_{\mu_d}(s)-1}
\]
where $a_{\mu_d}(u) =-1$. Therefore the only term that specializes to $1$ 
when $\sigma=1$ corresponds to $c=d$.  This implies that all terms in  \eqref{notrivial4} specializing to $1$ as $\sigma=1$ are obtained as follows. 

Let ${\sf S}_1(\mu, \lambda)$ denote the set of triples $(b,c,s)$ with 
$1\leq b, c \leq r$ and $s\in \lambda_b$, $s\neq u$ if $b=d$, such that 
\be\label{eq:vancondA} 
b-c+\ell_{\lambda_c}(s)=0, \qquad b-c-a_{\mu_b}(s)-1=0.
\ee
Let ${\sf S}_2(\mu, \lambda)$ denote the set of triples $(b,c,s)$ with 
$1\leq b, c \leq r$ and $s\in \mu_c$ such that 
\be\label{eq:vancondA} 
b-c-\ell_{\mu_b}(s)-1=0, \qquad b-c+a_{\lambda_c}(s)=0.
\ee
Then the terms in  \eqref{notrivial4} which specialize to $1$ as  $\sigma=1$ are given by 
\be\label{eq:ztermsA}
1+ \sum_{(b,c,s)\in {\sf S}_1(\mu, \lambda)} \sigma^{b-c}t^{b-c+\ell_{\lambda_c}(s)}q^{b-c-a_{\mu_b}(s)-1} + \sum_{(b,c,s)\in {\sf S}_2(\mu, \lambda)} \sigma^{b-c}t^{b-c-\ell_{\mu_b}(s)-1}q^{b-c+a_{\lambda_c}(s)}.  
\ee

Now consider a triple $(b,c,s)\in {\sf S}_1(\mu, \lambda)$. 
Since $s\neq u$ for $b=d$, it follows that $s\in \mu_b$. This implies $a_{\mu_b}(s)\geq 0$, hence $b\geq c+1$. 
At the same time, $\ell_{\lambda_c}(s) = c-b \leq -1$, hence $s\notin \lambda_c$. In conclusion, $b\geq c+1$ and  
$s\in \mu_{b} \setminus \lambda_c$ for any triple $(b,c,s)\in {\sf S}_1(\mu, \lambda)\setminus \{(d,d,u)\}$. Since $\mu_b\subseteq \lambda_b$ and $\mu_c \subseteq\lambda_c$, 
if $\mu$ or $\lambda$ is nested, this leads to a contradiction. Hence in that case ${\sf S}_1(\mu, \lambda)=\emptyset$. 

Suppose $\mu, \lambda$ are not nested, and the conditions of Lemma \ref{correspfixedlemmAB}.ii are satisfied. 
Since $\nu$ is nested 
and $\nu \subset \mu$ there is a unique index $1\leq e \leq r-1$ and a unique 
box $v \in \mu_{e+1}\setminus \nu_{e+1}$ such that $v \notin \mu_e$. For all other indices $1\leq b \leq r-1$, $b \neq e$ one has $\mu_{b+1}\subseteq \mu_b$. 
Then the argument in the previous paragraph implies that  
\[ 
{\sf S}_1(\mu, \lambda)= \{(e+1, e, v)\}. 
\]

A similar analysis applies to the set ${\sf S}_2(\mu, \lambda)$. Let 
$(b,c,s) \in {\sf S}_2(\mu, \lambda)$. 
Since $\mu_c \subseteq \lambda_c$ and $s\in \mu_c$, it follows that $a_{\lambda_c}(s)\geq 0$, hence $c\geq b$. 
If $c=b$, in order for $s\in {\sf S}_2(\mu, \lambda)$ one must have 
\[
\ell_{\mu_b}(s) = -1, \qquad a_{\lambda_b}(s) =0 
\]
This leads to a contradiction since $\mu_b \subseteq \lambda_b$ by assumption. 
Therefore one must have $c\geq b+1$. This implies that $a_{\lambda_c}(s)\geq 1$ and $\ell_{\mu_b}(s) \leq -2$. The second condition implies that $s\notin \mu_b$. If $\mu$ or $\lambda$ is nested this leads to a contradiction. Hence 
${\sf S}_2(\mu, \lambda)=\emptyset$.

In order to finish the proof, suppose $\mu,\lambda$ are not nested and the
conditions of Lemma \ref{correspfixedlemmAB}.ii are satisfied. 
Hence $\mu = \nu \cup\{v\}$ with $\nu$ nested. Since
$s\in \mu_c\setminus \mu_b$ with $b \leq c-1$ one must have $(b,c,s) = (e,e+1,v)$. However in that case $\ell_{\mu_b}(v) =-1$, which leads again to a contradiction. In conclusion ${\sf S}_2(\mu, \lambda)=\emptyset$. 

\hfill $\Box$

To conclude, note the following consequence of Lemma 
\ref{correspfixedlemmAB}. 
\begin{coro}\label{correspfixedcorA} 
 $(i)$ Let $\mu$ be a nested $r$-partition of $n$. 
Suppose a connected component $Z$ of $\CA(r,n,n+1)^{{\bf T}_0}$ 
has nontrivial set theoretic intersection with the closed subvariety $$\{\alpha_\mu\} \times \CA(r,n+1)\subset \CA(r,n)\times \CA(r,n+1).$$ Then $Z$ must be a closed point $(\alpha_\mu, \alpha_\lambda)$ for some $r$-partition $\lambda$ of $n+1$ such that $\lambda \supset \mu$. 

$(ii)$ Let $\lambda$ be a nested $r$-partition of $n+1$. 
Suppose a connected component $Z$ of $\CA(r,n,n+1)^{{\bf T}_0}$ 
has nontrivial set theoretic intersection with the closed subvariety $$\CA(r,n)\times \{\alpha_\lambda\} \subset \CA(r,n)\times \CA(r,n+1).$$ Then $Z$ must be a closed point $(\alpha_\mu, \alpha_\lambda)$ for some $r$-partition $\mu$ of $n$ such that $\mu \subset \lambda$. 
\end{coro}

 {\it Proof.} 
 $(i)$ Corollary \ref{CYfixedcorB}.$i$ shows that $\alpha_\mu$ is an isolated point of the fixed locus $\CA(r,n)^{{\bf T}_0}$. Therefore, $Z$ 
 must be contained as a closed subvariety  in the fiber $\pi_1^{-1}(\alpha_\mu)$. 
 Moreover, since the fixed locus $(\IA^1)^{{\bf T}_0}=\{0\}$, 
 any connected component of $\CA(r,n,n+1)^{{\bf T}_0}$ is a closed subvariety of $\CA_c(r,n,n+1)$. Since $\CA_c(r,n,n+1)$ is proper over 
 $\CA(r,n)$ by Proposition \ref{corresprop}.$3$, this implies that 
 $Z$ is proper over $\IC$, hence projective. Therefore the induced residual 
 ${\bf S}$-action on $Z$ has at least one fixed point $(\alpha_\mu, \alpha_\lambda)$ with $\mu \subset \lambda$. Since $\mu$ is nested by assumption, Lemma \ref{correspfixedlemmAB}.$i$ shows that this fixed point is ${\bf T}_0$-isolated. Then the claim follows since $Z$ is assumed connected.

 The proof of $(ii)$ is analogous. $Z$ must be again projective since 
 $\CA(r,n,n+1)$ is proper over $\CA(r,n,n+1)$.
 
\hfill $\Box$

\subsection{Hecke transformations for the Hilbert scheme}\label{Heckehilb} 

The goal of this section is to construct analogues of the Hecke transformations \eqref{eq:HeckemapA} and \eqref{eq:HeckemapB} for the
${\bf T}_0$-equivariant Borel-Moore homology of the moduli space ${Hilb}(r,n)$ 
of framed cyclic representations of the quiver \eqref{eq:quiver}. 
As shown in Proposition \ref{frquivpropA}, for each pair $(r,n)$ there is a closed embedding ${Hilb}(r,n) \subset \CA(r,n)$. 
Let ${{Hilb}(r,n,n+1)}$ denote the scheme theoretic intersection
of ${Hilb}(r,n) \times {Hilb}(r, n+1)$ and $\CA(r,n,n+1)$ in 
in $\CA(r,n) \times \CA(r,n+1)$. Hence one has a cartesian square 
\be\label{eq:HilbHeckecorresp}
\xymatrix{ 
{{Hilb}(r,n,n+1)}\ar[d] \ar[rr]^-{\kappa} & & \CA(r,n,n+1) \ar[d]^\gamma \\ 
{Hilb}(r,n)\times {Hilb}(r,n+1) \ar[rr]^-{} & & \CA(r,n) \times \CA(r,n+1). \\}
\ee
Let $\rho_1:{{Hilb}(r,n,n+1)} \to 
{Hilb}(r,n)$ and $\rho_2:{{Hilb}(r,n,n+1)}\to{Hilb}(r,n+1)$. The second is proper 
by base change since the projection $\pi_2\circ \gamma : \CA(r,n,n+1) \to \CA(r,n+1)$ is proper. However \ref{degDAHAaction}, 
one cannot define a refined Gysin pullback 
$\rho_1^!$ because $\rho_1$ is not l.c.i.  Furthermore, the construction of a virtual pullback as in 
 \cite{Virt_pullbacks} also fails because the relative obstruction theory 
of $\rho_1$ is not perfect of amplitude $[-1,\ 0]$. This precludes 
a straightfoward generalization of the Hecke transformations \eqref{eq:HeckemapA} and 
\eqref{eq:HeckemapB}.

One can use instead the embedding into the smooth ambient space $\CA(r,n)$. 
Lemma \ref{CYfixedlemmD} proves that the 
fixed locus ${Hilb}(r,n)^{{\bf T}_0}$ is finite and in one-to-one correspondence to nested $r$-partitions $\mu$ of $n$. Moreover, each fixed
${\bf T}_0$-fixed point $\alpha_\mu$ is isolated as a ${\bf T}_0$-fixed point in $\CA(r,n)$. Therefore the pushforward map for localized homology 
is injective and yields an identification 
\be\label{eq:Hilbsub}
 H^{{\bf T}_0}({Hilb}(r,n))_{K_0}\simeq \bigoplus_{\substack{\mu \in \calP(r,n)\\ \mu\, {\rm nested}}} K_0[\alpha_\mu] \subset 
  H^{{\bf T}_0}(\CA(r,n))_{K_0}~.
 \ee
 Here $K_0$ denotes the fraction field of the cohomology ring $H(B{\bf T}_0)$ and $[\alpha_\mu]=i_{\mu*}(1)$ for any 
nested $r$-partition $\mu$. 
Let also $\gamma : 
\CA(r,n,n+1)\hookrightarrow \CA(r,n)\times \CA(r,n+1)$ 
denote the natural closed embedding, which is clearly ${\bf T}_0$-equivariant. 
As in Section 2.2 of \cite{Cherednik_W}, any
equivariant Borel-Moore homology class $z\in H^{{\bf T}_0}(\CA(r,n,n+1))$ 
determines a Hecke transformation, 
\[
h_{z,n}^{+}:H^{{\bf T}_0}(\CA(r,n))_{K_0} \to H^{{\bf T}_0}(\CA(r, n+1))_{K_0}, \qquad 
h_{z,n}^{+}(x)=\pi_{2*}((\gamma_*z)\cdot \pi_1^*x).
\]
The pull-back and push-forward maps in equivariant Borel-Moore homology are well defined since $\gamma$ is a closed embedding and 
the restriction of $\pi_2$ to $\CA(r,n,n+1)$ is proper. 

As shown in Lemma \ref{correspfixedlemmAB}.i, for any pair of $r$-partitions $(\nu, \lambda)$ with $\nu$ nested, the closed point 
$(\alpha_\nu, \alpha_\lambda) \in \CA(r,n,n+1)$ is an isolated 
${\bf T}_0$-fixed point. Let $T_{\nu, \lambda}$ denote the tangent space to $\CA(r,n,n+1)$ at  $(\alpha_\nu, \alpha_\lambda)$. For any pair $(r,n)$ let 
$\calP(r,n)$ denote the set of $r$-partitions of $n$. 

\begin{lemm}\label{HeckelemmB} 
For any nested $r$-partition $\nu$, 
\be\label{eq:HeckemapC}
h_{z,n}^+([\alpha_\nu]) = \sum_{\substack{\lambda \in \calP_{r,n+1}\\  
\lambda\, {\rm nested}\\
\lambda \supset \nu\\ }}
e_{{\bf T}_0}(T_\nu) e_{{\bf T}_0}(T_{\nu,\lambda})^{-1}
z_{\nu, \lambda}[\alpha_\lambda],
\ee
where $z_{\nu, \lambda}$ denotes the restriction of $z$ to the 
fixed point $(\alpha_\nu, \alpha_\lambda)$. 
In particular,
the Hecke transformation \eqref{eq:HeckemapC} maps $H^{{\bf T}_0}({Hilb}(r,n))_{K_0}$ to $H^{{\bf T}_0}({Hilb}(r,n+1))_{K_0}$.
\end{lemm}

{\it Proof}. 
Note the cartesian diagram 
\[
\xymatrix{ 
\{\alpha_\nu\} \times \CA(r, n+1) \ar[r]^-{k_\nu} 
\ar[d]_-{\pi_{1,\nu}} & \CA(r,n)\times \CA(r,n+1) \ar[d]^-{\pi_1} \\
\{\alpha_\nu\} \ar[r]^-{i_\nu} & \CA(r,n) \\}
\]
where the horizontal maps are closed embeddings. 
As observed below equation 
$(2.8)$ in Section 2.1 of \cite{Cherednik_W}, this yields the base-change identity
\[ 
\pi_1^*[\alpha_\nu] = k_{\nu *}(\pi_{1,\nu}^*1) =  k_{\nu*}[\CA(r,n+1)]
\]
Since $k_\nu$ is a closed embedding of smooth varieties, 
there is an identity 
\[ 
(\gamma_*z)\cdot k_{\nu*}[\CA(r,n+1)] = k_{\nu*} k_{\nu}^* \gamma_*(z)
\]
in the intersection ring of $\CA(r,n)\times \CA(r,n+1)$. See for example Section 2.6.21 of \cite{Rep_th_geom}. 
For each connected component $Z$ of the fixed locus
$\CA(r,n,n+1)^{{\bf T}_0}$ let 
$q_Z:Z\hookrightarrow \CA(r,n,n+1)$ denote the natural closed embedding. Since the correspondence variety is smooth, the localization theorem yields 
\[
z= \sum_{Z} q_{Z*}\left(e_{{\bf T}_0}(\CV_Z)^{-1}q_Z^* z\right)
\]
where $\CV_Z$ is the normal bundle to $Z$ in $\CA(r,n,n+1)$. 
Then 
\be\label{eq:locHeckeC}
k_{\nu}^* \gamma_{*} z = 
\sum_{Z}  k_{\nu}^*\gamma_*q_{Z*}
\left( e_{{\bf T}_0}(\CV_Z)^{-1} q_Z^* z\right).
\ee
Clearly, if the set theoretic intersection of $Z$ with  $\alpha_\nu \times \CA(r,n+1)$ in $\CA(r,n)\times \CA(r, n+1)$ is empty, the corresponding term in the right hand side 
of \eqref{eq:locHeckeC} vanishes.
On the other hand, Corollary \ref{correspfixedcorA} shows that the connected components intersecting $\alpha_\nu \times \CA(r,n+1)$ nontrivially coincide with the finite set of isolated fixed points $(\alpha_\nu, \alpha_\lambda)$ where $\lambda \supset \nu$. For each such fixed point, 
let $q_{\nu, \lambda}:(\alpha_\nu, \alpha_\lambda)\hookrightarrow \CA(r,n+1,n)$ denote the natural closed embedding and let
\[
j_{\nu, \lambda}=\gamma \circ q_{\nu, \lambda}: \{(\alpha_\nu, \alpha_\lambda)\}\hookrightarrow \CA(r,n)\times \CA(r,n+1).
\]
Then 
\be\label{eq:locHeckeA}
k_{\nu}^* \gamma_{*} z =  \sum_{\substack{\lambda \in \calP_{r,n+1}\\ 
\lambda \supset \nu}} e_{{\bf T}_0}(T_{\nu,\lambda})^{-1}k_{\nu}^*(j_{\nu, \lambda})_* z_{\nu, \lambda},
\ee
where $T_{\nu, \lambda}$ is the tangent space to $\CA(r,n,n+1)$ at 
$(\alpha_\nu, \alpha_\lambda)$, and 
$z_{\nu,\lambda}=q_{\nu,\lambda}^*z$ denotes the restriction of $z$ to $(\alpha_\nu, \alpha_\lambda)$. 

Now let $i_Y : Y\hookrightarrow \CA(r, n+1)$ denote a connected component 
of the fixed locus $\CA(r,n+1)^{{\bf T}_0}$ and  let 
$k_{\nu, Y}= k_\nu \circ i_Y: \alpha_\nu \times Y \hookrightarrow 
\CA(r,n)\times \CA(r,n+1)$ denote the corresponding closed embedding. Let  also $N_Y$ denote the normal bundle to $Y$ in $\CA(r, n+1)$. 
The localization theorem yields 
\be\label{eq:locHeckeB}
k_{\nu}^* \gamma_{*} z = \sum_{Y}
i_{Y*} \left(e_{{\bf T}_0}(N_Y)^{-1}
 i_Y^*k_\nu^* \gamma_{*} z\right)
\ee
where 
$k_\nu^* \gamma_{*} z$ is given by the right hand side of equation 
\eqref{eq:locHeckeA}. 
Clearly, one has 
\[
i_Y^*k_\nu^*(j_{\nu, \lambda})_* z_{\nu, \lambda}= 
k_{\nu, Y}^*(j_{\nu, \lambda})_* z_{\nu, \lambda} =0 
\]
unless the component $Y$ contains the point $\alpha_\lambda$. 
If this is the case, Corollary \ref{CYfixedcorB} shows that 
\[
Y \simeq \left\{\begin{array}{ll} 
\{\alpha_\lambda\}, & {\rm for}\ \lambda\ {\rm nested} \\ & \\
\IA^1, &  {\rm for}\ \lambda\ {\rm not\ nested}.\\
\end{array}\right. 
\]
Moreover, in each case $\alpha_\lambda$ is the unique ${\bf S}$-fixed point in $Y$. 
For ease of exposition, such  a component will be denoted below by $Y_\lambda$ while their equivariant normal bundles in $\CA(r,n+1)$ will be denoted by $N_\lambda$. 

In conclusion, 
using equation \eqref{eq:locHeckeA},
 equation \eqref{eq:locHeckeB} reduces 
to 
\be\label{eq:locHeckeD} 
k_{\nu}^* \gamma_{*} z = \sum_{\substack{\lambda \in \calP_{r,n+1}\\ 
\lambda \supset \nu}}i_{\lambda*}\left(
e_{{\bf T}_0}(N_\lambda)^{-1} e_{{\bf T}_0}(T_{\nu,\lambda})^{-1}
i_\lambda^*k_{\nu}^*(j_{\nu, \lambda})_* z_{\nu, \lambda}\right).
\ee
Now let $f_\lambda: \{\alpha_\lambda \}\hookrightarrow 
Y_\lambda$ denote the natural closed embedding. Note that $\alpha_\nu \times Y_\lambda$ is a connected component of the ${\bf T}_0$ action on $\CA(r,n)\times \CA(r,n+1)$ and its equivariant 
normal bundle in the product is naturally isomorphic to 
\[
N_\lambda\oplus T_\nu\otimes \CO_{Y_\lambda}.
\]
Moreover, 
\[
j_{\nu, \lambda} = k_{\nu, \lambda} \circ f_\lambda, 
\]
hence  
\[
i_\lambda^*k_{\nu}^*(j_{\nu, \lambda})_* z_{\nu, \lambda}=
k_{\nu, \lambda}^*(k_{\nu, \lambda})_* f_{\lambda *} z_{\nu, \lambda} = 
e_{{\bf T}_0}(T_\nu) e_{{\bf T}_0}(N_\lambda) 
f_{\lambda*} z_{\nu, \lambda}. 
\]
Here $f_{\lambda*}$ is the pushforward map 
\[
f_{\lambda*}: H_0(\{\alpha_\lambda\}) \otimes_\IC 
K_{{\bf T}_0} \to H_0(Y_\lambda) \otimes_\IC 
K_{{\bf T}_0}.
\]
If $\lambda \supset \nu$ is not nested, $Y_\lambda\simeq \IA^1$, hence the degree zero Borel-Moore homology vanishes, 
$H_0(Y_\lambda) =0$. Therefore if this is the case, 
\[ 
f_{\lambda*} z_{\nu, \lambda}=0.
\]
If $\lambda\supset \nu$ is nested, $Y_\lambda$ coincides with $\alpha_\lambda$, hence 
\[
f_{\lambda*} z_{\nu, \lambda}= z_{\nu, \lambda}.
\]
Therefore equation \eqref{eq:locHeckeD} reduces to 
\[
k_{\nu}^* \gamma_{*} z = \sum_{\substack{\lambda \in \calP_{r,n+1}\\ \lambda\ {\rm nested}\\ \lambda \supset \nu\\ }}i_{\lambda*}\left(
e_{{\bf T}_0}(T_\nu) e_{{\bf T}_0}(T_{\nu,\lambda})^{-1}
z_{\nu, \lambda}\right).
\]
This yields equation \eqref{eq:HeckemapC}. 

\hfill $\Box$

In the opposite direction, let $[\alpha_\lambda]\in H^{{\bf T}_0}({Hilb}(r,n+1))_{K}$ be a basis element, where $\lambda$ is a
nested $r$-partition of $n+1$. Let $k_\lambda: \CA(r,n)\times \alpha_\lambda\hookrightarrow \CA(r,n)\times \CA(r,n_1)$ be the natural closed embedding.  Then 
\[
\pi_2^*[\alpha_\lambda] = k_{\lambda*} [\CA(r,n)\times \alpha_\lambda]
\]
and 
\[
(\gamma_* z) \cdot \pi_2^*[\alpha_\lambda] = k_{\lambda*}k_\lambda^*\gamma_*z. 
\]
The composition $\pi_1\circ k_{\lambda}: \CA(r,n)\times \alpha_\lambda \to \CA(r,n)$ is the identity, hence clearly, 
$\pi_{1*} ((\gamma_* z) \cdot \pi_2^*[\alpha_\lambda])$ is well defined.
This defines a Hecke transformation in the opposite direction, 
\[
h_{z,n+1}^-: H^{{\bf T}_0}({Hilb}(r,n+1))_{K} \to H^{{\bf T}_0}(\CA(r,n))_K.
\]
\begin{lemm}\label{HeckelemmC} 
For any nested $r$-partition $\lambda$, 
\be\label{eq:HeckemapE}
h_{z,n+1}^-([\alpha_\lambda]) = \sum_{\substack{\nu \in \calP_{r,n}\\ \nu\, {\rm nested}\\ \nu\subset \lambda \\ }}
e_{{\bf T}_0}(T_\lambda) e_{{\bf T}_0}(T_{\nu,\lambda})^{-1}
z_{\nu, \lambda}[\alpha_\nu].
\ee
In particular,
the Hecke transformation \eqref{eq:HeckemapE} maps $H^{{\bf T}_0}({Hilb}(r,n+1))_{K_0}$ to $H^{{\bf T}_0}({Hilb}(r,n))_{K_0}$.
\end{lemm} 

{\it Proof.} Completely analogous to the proof of Lemma \ref{HeckelemmB}. 

\hfill $\Box$

In conclusion, Lemmas \ref{HeckelemmB} and \ref{HeckelemmC} yield upward and downward Hecke transformations on the equivariant homology space  
\[
{\bf V}^{(r)}_{K_0} = \bigoplus_{n\geq 0} H^{{\bf T}_0}({Hilb}(r,n))_{K_0} \simeq \bigoplus_{n\geq 0}  \bigoplus_{\substack{\mu \subset \calP_{n,r}\\ \mu\, {\rm nested}}} K_{0}[\alpha_\mu]
\]
Abusing notation, they will be denoted by the same symbols,
$h^+_{z,n}$, $h_{z,n+1}^-$, the distinction being  clear from the context.

\section{Degenerate DAHA action}\label{DAHAsect}
 
The goal of this section is to prove that the transformations \eqref{eq:HeckemapC} and \eqref{eq:HeckemapE} yield a 
 degenerate DAHA action on ${\bf V}^{(r)}_{K_0}$.

\subsection{The degenerate DAHA of Schiffmann and Vasserot}\label{SHsect}  

This section is a brief review the construction of the family algebras 
${\bf {SH}^{\bf c}}$ introduced in \cite[Section 1]{Cherednik_W}, and further studied in \cite{DAHA_Yangian}. The construction employs a formal parameter 
$\kappa$ as well as an infinite set of formal parameters 
${\bf c}=({\bf c}_l)_{l\geq 0}$.
Using the notation of \cite[Section 1.5]{Cherednik_W}, for any $l\geq 0$ one defines 
\be\label{eq:SHnotation} 
\bal 
\xi & = 1-\kappa \\
G_0(s) & = -{\rm ln}(s) \\
G_l(s) & = (s^{-l}-1)/l, \quad l\neq 0\\
\varphi_l(s) & = \sum_{q=1,-\xi, -\kappa} s^l(G_l(1-qs) -G_l(1+qs))\\
\phi_l(s) & = s^lG_l(1+\xi s). \\
\eal
\ee
where $s$ is yet another formal parameter.  The right hand sides of the last two equations should be regarded as formal Laurent power series of $s$ by formally expanding the log functions. 
Then 
${\bf SH}^{\bf c}$ is generated by  $D_{-1,l}$, $D_{0,l}$, $D_{1,l}$, $l\in \IZ$, $l\geq 0$, satisfying the relations:
\be\label{eq:SHA}
\bal
&\ \, [D_{0,l}, D_{0,k}]  = 0, \qquad l,k\geq 1,\\
 &\ \, [D_{0,l}, D_{1,k}]  =  D_{1,l+k-1}, \qquad l\geq 1,\ 
 k\geq 0\\
 & [D_{0,l}, D_{-1,k}] = -D_{-1,l+k-1}, \qquad l\geq 1,\ 
 k\geq 0\\ 
 \eal 
 \ee
 \be\label{eq:SHB}
\bal  
& 3[D_{1,2},D_{1,1}] - [D_{1,3},D_{1,0}]+[D_{1,1},D_{1,0}] + 
\kappa(\kappa-1)(D_{1,0}^2 + [D_{1,1},D_{1,0}]) = 0 \\
& 3[D_{-1,2},D_{-1,1}] - [D_{-1,3},D_{-1,0}]+[D_{-1,1},D_{-1,0}] + 
\kappa(\kappa-1)(-D_{-1,0}^2 + [D_{-1,1},D_{-1,0}]) = 0 \\
\eal
\ee
\be\label{eq:SHC} 
\bal
 & [D_{-1,k}, D_{1,l}] =  E_{k+l}, \qquad l,k\geq 0.\\
 \eal
 \ee 
 \be\label{eq:SHD} 
[D_{1,0},[D_{1,0},D_{1,1}]]=0, \qquad [D_{-1,0}, [D_{-1,0},D_{-1,1}]] =0
\ee
The elements $E_l$ are expressed in terms of the generators 
 via the power series identity 
 \be\label{eq:SHE}
 1+\xi\sum_{l \geq 0} E_l s^{l+1} = {\rm exp}\big(\sum_{l\geq 0} (-1)^{l+1}
 {\bf c}_l \phi_l(s)\big) {\rm exp}\big(\sum_{l\geq 0} D_{0,l+1}
 \varphi_l(s)\big). 
 \ee
 The parameters ${\bf c}_l$, $l\geq 0$, are central. 
 Note that relations \eqref{eq:SHB} and \eqref{eq:SHC} were derived in 
 \cite{DAHA_Yangian}.
 Moreover, the following structure results were proven in Propositions (1.34) and (1.36) of \cite{Cherednik_W} respectively. 
 \begin{lemm}\label{SHgenlemm} 
 $(i)$ The algebra ${\bf SH}^{\bf c}$ is generated by the elements   
 ${\bf c}_l, D_{1,0}, D_{-1,0}, D_{0,2}$.
 
 $(ii)$ Any element of ${\bf SH}^{\bf c}$ can be written as a linear combination of monomials in the generators $D_{k,l}$ such that $D_{1,l}$, 
 $D_{0,l}$, $D_{-1,l}$, $l\geq 0$ appear exactly in this order from left to right. 
 \end{lemm} 
 

\subsection{Calabi-Yau DAHA}\label{CYDAHAsect}

In order to simplify the formulas let $K_{r+2} = K_{{\bf T}_{r+2}}$ denote the fraction field of the cohomology ring $H(B{\bf T}_{r+2})$. 
Then note that 
Theorem 3.2 in \cite{Cherednik_W} proves that a certain  specialization 
of ${\bf SH}^{\bf c}$ acts on the localized equivariant Borel-Moore homology 
\[
{\bf L}_{K_{r+2}}^{(r)} = \bigoplus_{n\geq 0}  
H^{{\bf T}_{r+2}}(\CA(r,n))_{K_{r+2}}
\]
via Hecke correspondences. The specialization
relates the formal parameters $\kappa, {\bf c}$ to the canonical 
generators $(x, y, e_1, \ldots, e_r)$ of the cohomology ring of the classifying space $B{\bf T}_{r+2}$ as shown below.  Let $\IC(\kappa)[{\bf c}] = \IC(\kappa)[{\bf c}_0, {\bf c}_1, \ldots]$  and let 
\be\label{eq:SHspecA} 
\IC(\kappa)[{\bf c}] \to K_{r+2}
\ee
be the algebra homomorphism 
mapping 
\[
\kappa\mapsto  - x^{-1}y, \qquad {\bf c}_l \mapsto p_l(\epsilon_1, \ldots, \epsilon_r),
\]
where $\epsilon_a= x^{-1}e_a$, $1\leq a\leq r$,  and
$p_l$, $l\geq 0$, are the symmetric power functions in $r$ variables.
Using the $\IC(\kappa)[{\bf c}]$-module structure on  $K_{{\bf T}_{r+2}}$ obtained from \eqref{eq:SHspecA}, let 
\[
{\bf SH}^{(r)}_{K_{r+2}}= {\bf SH}^{\bf c}\otimes_{\IC(\kappa)[{\bf c}]} K_{r+2}.
\]
Let $(x,y,z)$ be the canonical generators of the cohomology ring of $B{\bf T}$.   Then ${\bf T}$-specialization is defined by setting 
\be\label{eq:Tspec}
\epsilon_a = (1-a) x^{-1}z , \qquad 1\leq a\leq r. 
\ee
Finally, the  specialization to ${\bf T}_0\subset {\bf T}$ is obtained by further setting $z= -(x+y)$, which yields 
\be\label{eq:Tzerospec}
\epsilon_a = (a-1)\xi, \qquad 1\leq a \leq r, 
\ee
where $\xi = 1-\kappa$ as in \eqref{eq:SHnotation}. 

Note that both specializations are well defined since the right hand side of equation \eqref{eq:SHB} is a formal power series of ${\bf c}_l$, $l\geq 0$. 
In particular the ${\bf T}_0$-specialization of ${\bf SH}^{\bf c}$ will be denoted by ${\bf SH}^{(r)}_{K_0}$. 
  
\subsection{${\bf SH}^{(r)}_{K_0}$-module structure}\label{SHmodsect} 

This section will prove that the specialization ${\bf SH}^{(r)}_{K_0}$ acts on 
${\bf V}_{K_0}^{(r)}$ using the  Hecke transformations \eqref{eq:HeckemapC} 
and \eqref{eq:HeckemapE}. 
First recall the subalgebra ${\bf U}^{(r)}_{K_{r+2}}$ of the convolution algebra constructed in Section 3.5 of \cite{Cherednik_W}.
As in Section 3.4 of loc. cit., let $\tau_{n+1,n}$ denote the universal line bundle on the correspondence variety $\CA(r,n,n+1)$. Using the notation of Section \ref{correspvar}, the algebra ${\bf U}_{K_{r+2}}^{(r)}$ is generated by the Hecke transformations 
\be\label{eq:SHeckeA}
f_{1,l} = \prod_{n\geq 0} f_{1,l,n}, \qquad f_{1,l,n}(w)=
\pi_{2*}(\gamma_* c_1(\tau_{n+1,n})^l \cdot \pi_1^*(w)), \qquad 
l \geq 0, 
\ee
\[ 
f_{-1,l} = \prod_{n\geq 0} f_{-1,l,n}, \qquad f_{1,l,n}(w)=
\pi_{1*}(\gamma_* c_1(\tau_{n,n+1})^l \cdot \pi_2^*(w)) , \qquad l\geq 0, 
\]
\[
e_{0,l} = \prod_{n\geq 0} e_{0,l,n}, \qquad 
e_{0,l,n}(w) = c_l(E_n)\cdot w , \qquad l\geq 0,
\]
where $E_n$ is the universal rank $n$ vector bundle over 
$\CA(r,n)$. The Chern classes $c_l$, $l\geq 0$, in the above formulas are ${\bf T}_{r+2}$-equivariant. In addition one defines the diagonal operators 
\be\label{eq:SHeckeB} 
f_{0,l}([\alpha_\mu]) = \sum_{a=1}^r\sum_{s\in \mu_a} c_a(s)^l [\alpha_\mu], \qquad c_a(s) = i(s) x + j(s) y - e_a
\ee
which can be written as polynomial functions of the $e_{0,l}$, $l\in \IZ$.

Then Theorem 3.2 of \cite{Cherednik_W} proves that there is a unique isomorphism 
of algebras 
\be\label{eq:SHconvA}
{\bf SH}^{(r)}_{K_{r+2}} {\buildrel \sim \over \longto} 
{\bf U}^{(r)}_{K_{r+2}} 
\ee
mapping 
\be\label{eq:SHconvB}
D_{1,l}\mapsto x^{1-l}yf_{1,l}, \qquad 
D_{0,l} \mapsto x^{1-l}f_{0,l-1}, \qquad 
D_{-1,l}\mapsto (-1)^{r-1} x^{-l}f_{-1,l}.
\ee
As observed in Remark \ref{Tcorresp} the above construction admits a 
straightforward specialization to ${\bf T}\subset {\bf T}_{r+2}$. This results in an isomorphism 
\be\label{eq:SHconvB}
{\bf SH}^{(r)}_{K} {\buildrel \sim \over \longto} 
{\bf U}^{(r)}_K. 
\ee
Abusing notation, the ${\bf T}$-equivariant counterparts of the above Hecke transformations will be denoted by the same symbols. The distinction will be clear from the context. 
Note also that explicit formulas for the matrix elements of the generators in the above representation follow immediately from Lemma \ref{HeckelemmA}.

Next note the direct sum decomposition 
\be\label{eq:nestedsplit}
{\bf L}^{(r)}_{{\bf T}} \simeq 
{\bf V}_{K} \oplus {\bf V}_{K}^\perp~,
\ee
where 
\[
{\bf V}_K = \bigoplus_{\mu\, {\rm nested}} 
K[\alpha_\mu], \qquad {\bf V}_K^\perp = \bigoplus_{\mu\, {\rm not\ nested}} 
K[\alpha_\mu]~.
\]
Suppose $A$ is one of the $K$-linear transformations $f_{1,l,n}$, $f_{-1,l,n}$, $f_{0,l,n}$ which generate the convolution algebra. 
Let 
\[ 
A= \left[\begin{array}{cc} A_{11} & A_{12} \\ A_{21} & A_{22}\\ \end{array} \right]
\]
be the block form decomposition of $A$ with respect to the direct sum \eqref{eq:nestedsplit}.
Then define the $K$-linear transformations 
\be\label{eq:trHeckeA} 
g_{i,l,n}\in {\rm End}_K({\bf V}_K), \qquad i\in \{-1,0,1\}, \ l\in \IZ,\ n \in \IZ_{\geq 0} 
\ee
as 
\[
g_{i,l,n} = \left(f_{i,l,n}\right)_{11}. 
\]
Recall that $K$ is 
 the field of rational functions $\IC(x,y,z)$, where $(x,y,z)$ are 
 the canonical generators of the cohomology ring of $B{\bf T}$. 
Let $s= c_1(\sigma)$ where $\sigma(t_1,t_2,t_3) = t_1t_2t_3$ is the character used in Lemma \ref{correspfixedlemmAB}, and note that $s= -(x+y+z)$. At the same time, ${K}_0=\IC(x,y)$ and there is a canonical isomorphism of 
${K}$-vector spaces 
\be\label{eq:blockdecomp}
{\bf V}_{K} \simeq {\bf V}_{K_0} \otimes_{K_0} K
\ee
mapping 
\[
[\alpha_\mu] \mapsto [\alpha_\mu]\otimes 1
\]
for any nested $r$-partition $\mu$. Let 
\[
g_{i,l} = \prod_{n \geq 0} g_{i,l,n}. 
\]
Finally, recall that each {\bf T}-fixed point $\alpha_\mu\in \CA(r,n)^{\bf T}$, with $\mu$ nested, 
is ${\bf T}_0$-isolated  and the equivariant Euler class \[
e_{{\bf T}_0}(T_\mu) = e_{\bf T}(\alpha_\mu)|_{s=0}
\]
is nonzero. Similarly, as shown in Lemma \ref{correspfixedlemmAB}.$i$, for any nested partitions $\mu \subset \lambda$ the fixed point $(\alpha_\mu, \alpha_\lambda)\in \CA(r,n,n+1)^{\bf T}$ is ${\bf T}_0$-isolated and 
\[
 e_{{\bf T}_0}(T_{\mu, \lambda}) = e_{\bf T}(T_{\mu,\lambda})|_{s=0}
\]
is nonzero. Moreover, for any box $s\in \mu$ set 
\[ 
c_a^0(s) = (i(s)-a)x + (j(s)-a)y 
\]
where $1\leq a \leq r$ indicates that $s\in \mu_a$, 
and for any pair $\mu\subset \lambda$ with $|\lambda\setminus \mu|=1$ set 
\[ 
\tau_{\mu, \lambda}^0 = (i(s)-a)x+(j(s)-a)y
\]
where $\{s\}=\lambda\setminus \mu$. 

Then the following holds 
\begin{lemm}\label{CYtruncationA}
The matrix elements of the $K$-linear transformations $g_{i,l} \in {\rm End}_{K}({\bf V}_{K})$, $i\in \{-1,0,1\}, \ l\in \IZ$ with respect to the fixed point basis $\{[\alpha_\mu]\}$, with $\mu$ a nested $r$-partition, have well defined specializations $g_{i,l}|_{s=0}$ at $s=0$. In particular there exist unique $K_0$-linear transformations $g_{i,l}^0\in {\rm End}_{K_0}({\bf V}_{K_0})$
such that $g_{i,l}|_{s=0}=g_{i,l}^0\otimes {\bf 1}$. Moreover, the 
 explicit expressions of $g_{i,l}^0$ in the fixed point basis are given by 
\be\label{eq:HeckemapF} 
\bal 
g_{0,l+1}^0([\alpha_\mu]) & = \sum_{a=1}^r\sum_{s\in \mu} c_a^0(s)^l [\alpha_\mu]~, \\
g_{1,l}^0([\alpha_\mu]) & = \sum_{\substack{\lambda\ {\rm nested} \\\lambda\supset \mu\\ |\lambda| = |\mu|+1}} (\tau_{\mu, \lambda}^0)^l e_{{\bf T}_0}(T_\mu) e_{{\bf T}_0}(T_{\mu, \lambda})^{-1} 
[\alpha_\lambda]~, \\
g_{-1,l}^0([\alpha_\lambda]) & = \sum_{\substack{\mu\ {\rm nested} \\\mu \subset \lambda \\ |\lambda| = |\mu|+1}} (\tau_{\mu, \lambda}^0)^l e_{{\bf T}_0}(T_\lambda) e_{{\bf T}_0}(T_{\mu, \lambda})^{-1} 
[\alpha_\lambda]~. \\
\eal 
\ee
\end{lemm} 

{\it Proof}. Let $A$ be an element of $\{f_{0,l,n},\ f_{1,l,n}\ f_{-1,l,n}\}$.  
Consider the following cases: 

$(1)$ Suppose $A$ is one of the $f_{0,l,n}$. Then the above block decomposition is diagonal and the claim is obvious.  

$(2)$ Suppose $A$ is one of the $f_{1,l,n}$. Then Lemma \ref{CYfixedlemmB}.$i$ shows that for any nested $r$-partition $\mu$ of $n$ the equivariant Euler class $e_{{\bf T}}(\alpha_\mu)\in \IC(x,y,z)$ has well defined specialization at $s=0$, which is furthermore equal to $e_{{\bf T}_0}([\alpha_\mu])$. Moreover, for any pair of nested $r$-partitions $\mu, \lambda$ of $n,n+1$ respectively, with $\mu\subset \lambda$, Lemma 
\ref{correspfixedlemmAB}.$i$ shows that the equivariant Euler class $e_{{\bf T}}(T_{\mu,\lambda})\in \IC(x,y,z)$ has well defined specialization at $s=0$, which is furthermore equal to $e_{{\bf T}_0}(T_{\mu, \lambda})$. 
Then the claim follows from equations \eqref{eq:HeckemapC}. 

$(3)$ Suppose $A$ is one of the $f_{-1,l,n}$. This case is completely analogous to $(2)$. 

Equations \eqref{eq:HeckemapF} follow immediately by specialization from 
\eqref{eq:HeckemapAA} and \eqref{eq:HeckemapBB}. 

\hfill $\Box$

\begin{lemm}\label{SHaction}
The map
\be\label{eq:SHconvC}
D_{1,l}\mapsto x^{1-l}y g^0_{1,l}, \qquad 
D_{0,l} \mapsto x^{1-l} g^0_{0,l-1}, \qquad 
D_{-1,l}\mapsto (-1)^{r-1} x^{-l} g^0_{-1,l}.
\ee
extends uniquely to a homomorphism of $K_0$-algebras 
\be\label{eq:SHconvD} 
\rho_0^{(r)}:{\bf SH}^{(r)}_{K_0} \to {\rm End}({\bf V}^{(r)}_{K_0}). 
\ee
\end{lemm} 

{\it Proof}. The proof will proceed by truncating the relations satisfied by the generators \eqref{eq:SHconvB} to their $(1,1)$ blocks and specializing to $s=0$. 
Again, suppose $A$ is one of the transformations $f_{i,l,n}$, $i \in \{-1,0,1\}$, $l\in \IZ$, $n \in \IZ_{\geq 0}$ and consider the following cases. 

$(1)$ Suppose $A$ is one of the $f_{0,l,n}$. Then, clearly, $A_{12}=0$ and $A_{21}=0$. 

$(2)$ Suppose $A$ is one of the $f_{1,l,n}$. Then note that the matrix elements of 
\[
A_{12}: {\bf V}_K^\perp \to {\bf V}_K, \qquad A_{21}: {\bf V}_K
\to {\bf V}_K^\perp 
\]
are given by
\[
\left(A_{12}\right)_{\lambda,\mu} = e_{\bf T}(\alpha_\lambda) e_{\bf T}(\alpha_\mu, \alpha_\lambda)^{-1}, \qquad 
\left(A_{21}\right)_{\nu, \rho} = e_{\bf T}(\alpha_\nu) e_{\bf T}(\alpha_\nu, \alpha_\rho)^{-1}
\]
where $\mu, \nu$ are nested $r$-partitions, 
\[
|\lambda| = |\mu|+1, \qquad \mu \subset \lambda~,
\]
and 
\[
|\rho| = |\nu|+1, \qquad \nu \subset \rho~.
\]
Moreover, $\mu, \nu$ are nested while $\lambda, \rho$ are not. 
Then Lemma \ref{CYfixedlemmB}.$ii$ shows that 
\[
e_{\bf T}(T_\lambda) = s {\tilde e}(T_\lambda) 
\]
where ${\tilde e}(T_\lambda)$ has well defined specialization at $s=0$. At the same time 
Lemma \ref{correspfixedlemmAB}.$i$ shows that $e_{\bf T}(T_{\mu, \lambda})^{-1}$ has well defined specialization at $s=0$. 
In conclusion,
\[
A_{12} = s {\tilde A}_{12} 
\]
where ${\tilde A}_{12}$ has well defined specialization at $s=0$. 
Similarly, Lemmas \ref{CYfixedlemmB}.$i$ and \ref{correspfixedlemmAB}.$i$
imply that $A_{21}$ has well defined specialization at $s=0$.

$(3)$ Suppose $A$ is one of the $f_{-1,l,n}$. In complete analogy to $(2)$, 
Lemmas \ref{CYfixedlemmB}  and \ref{correspfixedlemmAB} 
imply that, again, 
\[
A_{12} = s {\tilde A}_{12} 
\]
where ${\tilde A}_{12}$ has well defined specialization at $s=0$. At the same time $A_{21}$ has well defined specialization at $s=0$.

The above observations imply that the claim holds for the quadratic relations  \eqref{eq:SHA} -- \eqref{eq:SHC}. 
 Suppose $A,B$ are two linear transformations among 
the $f_{i,l,n}$, $i \in \{-1,0,1\}$, $l\in \IZ$, $n \in \IZ_{\geq 0}$
such that the target of $B$ coincides with the domain of $A$. 
Using the block form decomposition \eqref{eq:blockdecomp}, the product 
$AB$ is written as
\be\label{eq:quadrelA}
AB = \left[\begin{array}{cc} 
A_{11}B_{11} + s {\tilde A}_{12}B_{21} & s A_{11} {\tilde B}_{12} + 
s {\tilde A}_{12}B_{22} \\ 
A_{21}B_{11} + A_{22} B_{21} & s A_{21} {\tilde B}_{12} + 
{ A}_{22}B_{22} \\
\end{array}\right]
\ee
Note that the product $A_{11}B_{11}$ has well defined specialization at $s=0$ by Lemma \ref{CYtruncationA}. Then, using $(1)$, $(2)$ and $(3)$ above it follows that the same holds for the component $(AB)_{11}$, and 
\[
(AB)_{11} |_{s=0} = \left(A_{11}|_{s=0} \right) \left(B_{11}|_{s=0} \right). 
\]
This proves that the transformations $g_{i,l}^0$ satisfy the quadratic relations \eqref{eq:SHA} -- \eqref{eq:SHC}. 

In order to prove the cubic relations \eqref{eq:SHD}, suppose $A,B,C$ are 
transformations of the form $f_{1,l,n+1}$, $f_{1,l,n}$, $f_{1,l,n-1}$ respectively, where $n\geq 1$. 
The $(1,1)$ block of the triple product $ABC$ reads 
\be\label{eq:tripleoneoneA}
A_{11}B_{11}C_{11} + s {\tilde A}_{12}B_{21}C_{11} + s A_{11}{\tilde B}_{12}C_{21}+ s {\tilde A}_{12}B_{22}C_{21}. 
\ee
Again, remarks $(1)$, $(2)$, $(3)$ above imply that 
\[
A_{11}B_{11}C_{11},\qquad {\tilde A}_{12}B_{21}C_{11}, \qquad A_{11}{\tilde B}_{12}C_{21}
\]
have well defined specializations at $s=0$, and 
\be\label{eq:triplespec}
(A_{11}B_{11}C_{11})|_{s=0} =  \left(A_{11}|_{s=0} \right) \left(B_{11}|_{s=0} \right) \left(C_{11}|_{s=0} \right).
\ee
Using Lemma \ref{correspfixedlemmAB}.$ii$, it will be shown below that the product ${\tilde A}_{12}B_{22}C_{21}$  also has well defined specialization at $s=0$. The matrix elements of 
the triple product ${A}_{12}B_{22}C_{21}$ are given by 
\[ 
\sum_{\nu\subset\mu\subset \lambda\subset \rho} {e_{\bf T}(T_\nu)\over 
e_{\bf T}(T_{\nu,\mu})} {e_{\bf T}(T_\mu)\over 
e_{\bf T}(T_{\mu,\lambda})} {e_{\bf T}(T_\lambda)\over 
e_{\bf T}(T_{\lambda, \rho})} 
\]
where $\nu, \rho$ are nested $r$-partitions of $n-1, n+2$ respectively and the sum is over all sequences of $r$-partitions $\nu\subset \mu \subset \lambda \subset \rho$ with $|\mu|=n$, $|\lambda|=n+1$, and $\mu, \lambda$ not nested. 
Lemma \ref{CYfixedlemmB}.$i$ and Lemma \ref{correspfixedlemmAB}.$i$ imply that 
\[
{e_{\bf T}(T_\nu)\over 
e_{\bf T}(T_{\nu,\mu})}
\]
has well defined specialization at $s=0$. 
Moreover, Lemma \ref{CYfixedlemmB}.$ii$ and Lemma \ref{correspfixedlemmAB}.$ii$ imply that 
\[
e_{\bf T}(T_\mu) = s\, {\tilde e}(T_\mu) \qquad 
{\rm{and} }\qquad
e_{\bf T}(T_{\mu, \lambda})^{-1} = s^{-1}\, {\tilde e}(T_{\mu, \lambda})^{-1}
\]
where ${\tilde e}(T_\mu)$, ${\tilde e}(T_{\mu, \lambda})^{-1}$ have well defined specializations at $s=0$. Therefore the same holds for 
\[
{e_{\bf T}(T_\mu)\over 
e_{\bf T}(T_{\mu,\lambda})}
\]
Similarly, Lemma \ref{CYfixedlemmB}.$ii$ and Lemma \ref{correspfixedlemmAB}.$i$ imply that 
\[
{e_{\bf T}(T_\lambda)\over 
e_{\bf T}(T_{\lambda, \rho})} =  {s\, {\tilde e}(T_\lambda)}
e_{\bf T}(T_{\lambda,\rho})^{-1}
\]
where ${\tilde e}(T_{\lambda}), e_{\bf T}(T_{\lambda,\rho})^{-1}$ also have well defined specialization at $s=0$. Since ${A}_{12}B_{22}C_{21}=s{\tilde A}_{12}B_{22}C_{21}$, 
it follows that, indeed, ${\tilde A}_{12}B_{22}C_{21}$ has well defined specialization at $s=0$ as claimed above. Then, using relation \eqref{eq:triplespec}
this implies that the transformations $g_{1,l}^0$ also satisfy the cubic relations \eqref{eq:SHD}.

\hfill $\Box$

{\it Proof of Theorem \ref{mainthmA}.} 
For any  $l\geq 0$
set $z=e_{{\bf T}_0}(\tau_{n,n+1})^l$ in Lemmas \ref{HeckelemmB} and 
\ref{HeckelemmC}. 
Let $h^+_{l}, h^-_{l}$ denote the resulting linear transformations in ${\rm End}({\bf V}^{(r)}_{K_0})$. Then equations \eqref{eq:HeckemapC},  
\eqref{eq:HeckemapE} and \eqref{eq:HeckemapF} show that 
\[
h^+_l = g^0_{1,l}, \qquad h^-_l = g^0_{-1,l} 
\]
for all $l\geq 0$. Moreover $h^0_l = g^0_{0,l}$ holds by construction for all $l\geq 0$. 
Therefore Theorem \ref{mainthmA} follows from Lemma
\ref{SHaction}.

\hfill $\Box$

\subsection{Some structure results}\label{SHstructure} 

The next goal is to prove some structure results for the action 
\eqref{eq:SHconvD} which are needed in the proof of Theorem 
\ref{mainthmB}. 
\begin{lemm}\label{SHfaitful} 
For $r=1$ the representation \eqref{eq:SHconvD} is faithful.
\end{lemm} 

{\it Proof.} 
 This follows from \cite[Prop.6.7]{Cherednik_W} since for $r=1$ the representations 
$\rho_0^{(r)}$ and $\rho^{(r)}$ are isomorphic. 

\hfill $\Box$ 


Next, one has to prove the analogues of Lemmas 8.33 and 8.34.a in \cite{Cherednik_W}.  Let $[\alpha_{\emptyset^r}]\in {\bf V}^{(r)}_{K_0}$ denote the element corresponding to the empty $r$-partition.

\begin{lemm}\label{SHvaclemm}
$[\alpha_{\emptyset^r}]$ is annihilated by all endomorphisms of the form 
$\rho_0^{(r)}(D_{0,l})$, $\rho_0^{(r)}(D_{-l,0})$, $l\geq 1$. 
\end{lemm}

{\it Proof}. This is analogous to Lemma 8.34.a in \cite{Cherednik_W}. The claim follows from equations \eqref{eq:HeckemapE}. 

\hfill $\Box$

 Lemma 8.33 in \cite{Cherednik_W} shows that 
\[
{\bf L}^{(r)}_{K_{r+2}} = \rho^{(r)}({\bf SH}^{(r)}_{K_r+2})([\alpha_{\emptyset^r}])
\]
The proof of loc. cit. is based on the following observation. 
Recall that $\calP_{r,n}$ denotes the set of $r$-partitions of $n$.
Let $\phi: \calP_{r,n}\to \big(\IZ^{2}\times \IZ^r\big)^{n}/\CS_n$ be the map defined by 
\[ 
\phi(\mu) = [(i(s),j(s),{\bf e}_a)] 
\]
where in the right hand side $s\in \mu_a$ and ${\bf e}_a$ are the canonical generators of $\IZ^r$ for $1\leq a\leq r$. Then $\phi$ is injective. 

In order to prove the analogue of Lemma 8.33 under Calabi-Yau specialization, one has to first prove the following analogous result.  
\begin{lemm}\label{injmap}
Let $\calP_{r,n}^{\subset} \subset \calP_{r,n}$ be the subset of nested 
$r$-partitions of $n$. Let 
\[
\phi_0: \calP_{r,n}^{\subset} \to (\IZ^2)^{\times n}/\CS_n
\]
be the map defined by 
\[ 
\phi_0(\mu) = [(i(s)-a, j(s)-a)]
\]
where in the right hand side $1\leq a \leq r$ and $s\in \mu_a$. Then 
$\phi_0$ is injective.
\end{lemm} 

{\it Proof}. Suppose $\mu=(\mu_a)$ and $\lambda=(\lambda_a)$, $1\leq a\leq r$,  are two nested $r$-partitions of $n$ such that $\phi_0(\mu)=\phi_0(\lambda)$. For each $1\leq a \leq r$ 
the  partition $\mu_a\subset \IZ^2$ can be written as a union of L-shaped sets in the plane 
\[
\bal
\mu_a  =\ & \cup_{k\in \IZ} \mu_{a,k},\\
\mu_{a,k} =\  & \{(i,j) \in \IZ^2\, |\, (i,j)\in \mu_a,\ 
i = k,\, \ j\geq k\}\ \cup \\
& \{ (i,j) \in \IZ^2\, |\, (i,j)\in \mu_a,\ 
i \geq k,\, \ j= k\}.\\
\eal
\]
Clearly, $\mu_{a,k}=\emptyset$ for $k\leq 0$ and $\mu_{a,k}\cap \mu_{a,l}=\emptyset$ for $k \neq l$. The translation of each set $\mu_{a,k}$ by $(-a,-a)$ will be denoted by $\mu_{a,k}-(a,a)$. Then note that 
\be\label{eq:inclA}
\mu_{a,k+l}-(l,l) \subseteq \mu_{a,k}  
\ee
for any $l \geq 0$ since 
since $\mu_a$ is a Young diagram. 
Moreover, for any $1\leq a, b\leq r$,
and any $k,l\in \IZ$  the subsets
$\mu_{a,k}-(a,a)$ and $\mu_{b,l}-(b,b)$ of $\IZ^2$ are disjoint unless $k-a = l-b$. 
For any $n \in \IZ$ let ${\sf S}_n(\mu)$ be the disjoint union 
\[
{\sf S}_n(\mu) = \coprod_{\substack{1\leq a \leq r,\, k \in \IZ\\ k-a =n}}
\ (\mu_{a,k}-(a,a))
\]
and let 
\[
\big|{\sf S}_n(\mu)\big| = \bigcup_{\substack{1\leq a \leq r,\, k \in \IZ\\ k-a =n}}
\ (\mu_{a,k}-(a,a))
\]
denote the set theoretic union as subsets of $\IZ^2$. 
As observed above, $|{\sf S}_n(\mu)|\cap |{\sf S}_m(\mu)|=\emptyset$  for $n\neq m$. Hence 
 \[
\coprod_{a=1}^r (\mu_a-(a,a)) = \bigcup_{n\in \IZ} {\sf S}_n(\mu),
\]
Since $\mu$ is nested, inclusion \eqref{eq:inclA} implies that 
there are inclusions 
\be\label{eq:inclB}
\mu_{a+1,k+a-r+1}\subseteq \mu_{a+1,k+a-r}\subseteq \mu_{a,k+a-r}, \qquad 1\leq a\leq r, \qquad k\geq 1, 
\ee
where by convention $\mu_{r+1,k+1}=\emptyset$. Therefore the following statement holds for each set ${\sf S}_{k-a}(\mu)$:

$(i)$ Each point 
\[
(i,j)\in \mu_{a,k+a-r}\setminus \mu_{a+1,k+a-r+1}
\]
occurs with multiplicity exactly $a$ in ${\sf S}_{k-a}(\mu)$. 

In complete analogy, one also has 
\[
\coprod_{a=1}^r (\lambda_a-(a,a)) = \bigcup_{m\in \IZ} {\sf S}_m(\lambda) 
\]
where the subsets ${\sf S}_m(\lambda)$ satisfy analogous properties.
In particular, one has:

$(ii)$ Each point 
\[
(i,j)\in \lambda_{a,k+a-r}\setminus \lambda_{a+1,k+a-r+1}
\]
occurs with multiplicity $a$ in ${\sf S}_{k-a}(\lambda)$.

In order to finish the proof  note that 
$|{\sf S}_n(\mu)|\cap |{\sf S}_m(\lambda)|=\emptyset$ for $n\neq m$. 
 Therefore $\phi_0(\mu)=\phi_0(\lambda)$ if and only if 
${\sf S}_n(\mu)={\sf S}_n(\lambda)$ for all $n \in \IZ$. If that is the case, properties 
$(i)$ and $(ii)$ above imply that 
\[
 \mu_{a,k+a-r}\setminus \mu_{a+1,k+a-r+1} = \lambda_{a,k+a-r}\setminus \lambda_{a+1,k+a-r+1}
 \]
 for all $1\leq a\leq r$ and all $k \in \IZ$. This implies that $\mu=\lambda$.

\hfill $\Box$ 

Now one has: 
 
\begin{lemm}\label{SHmodA} 
${\bf V}^{(r)}_{K_0} = 
\rho_0^{(r)}({\bf SH}^{(r)}_{K_0})([\alpha_{\emptyset^r}])$.
\end{lemm}

{\it Proof}. This is analogous to Lemma 8.33 in \cite{Cherednik_W}. The proof proceeds by induction on $n$ i.e. suppose $[\alpha]_\mu\in \rho_0^{(r)}({\bf SH}^{(r)}_{K_0})(|0\rangle)$ for all $r$-partitions $\mu$ of $n-1$. 
Recall the formula for $g_{1,l}^0$ in equation \eqref{eq:HeckemapE}. 
\[ 
{ g}^0_{1,l}([\alpha_\mu]) = \sum_{\substack{\lambda \in \calP_{r,n}\\ \lambda\, {\rm nested}\\ \lambda \supset \mu\\ }} (\tau^0_{\mu, \lambda})^l
e_{{\bf T}_0}(T_\mu) e_{{\bf T}_0}(T_{\mu,\lambda})^{-1}
[\alpha_\lambda].
\]
As in loc. cit., this implies that for any 
 $r$-partition $\lambda$ 
with $|\lambda|=n$ one can find a nested partition $\mu$, $|\mu|=n-1$ such that the coefficient of $[\alpha_\lambda]$ in $\rho_0^{(r)}(D_{1,l})([\alpha_\mu])$ is nonzero for some $l\geq 0$. 
Next note that 
\[
{ g}^0_{0,l+1}([\alpha_\lambda] )= \sum_{a=1}^r \sum_{s\in \lambda_a} 
\big((i(s)-a)x+(j(s)-a)y\big)^l[\alpha_\lambda] 
\]
from \eqref{eq:HeckemapE}. By analogy with Lemma 8.33, there is a map 
\[ 
\calP_{r,n}^{\subset} \to (K_0^2)^n/\CS_n 
\]
mapping 
\[
\lambda \mapsto [(i(a)-a)x+(j(a)-a)y)], \qquad s\in \lambda_a, \ 1\leq a\leq r.
\]
Lemma \ref{injmap} shows that this map is injective. Hence, as in loc. cit., the Hilbert nullstellensatz implies that there 
 there is a polynomial $f$ in the 
generators $D_{0,l}$ such that $f([\alpha_\lambda])=1$ and 
$f([\alpha_\rho])=0$ for any $r$-partition $\rho\neq \lambda$ of $n$.

\hfill $\Box$

\section{$W$-module structure}\label{Wmodsect}

This section provides a brief overview of $W$-algebras and concludes the proof of Theorem \ref{mainthmB}. 

\subsection{$W$-algebras}\label{Wsect} 
A succinct definition of the $W$-algebra $W_\kappa({\mathfrak gl}_r)$ and its free field realization via quantum Miura transform is presented in \cite[Section 8.4]{Cherednik_W}. The main points will be briefly 
summarized below.

In this section the ground field is $F=\IC(\kappa)$. 
Let $b_a$, $1\leq a \leq r$, be a basis of the standard Cartan algebra ${\mathfrak h}\subset {\mathfrak gl}_r$ and let $b^{(a)}$, $1\leq a \leq r$ be the dual basis. Let $b^{(a)}(z)$, $1\leq a\leq r$ be free boson fields such that the zero modes $b^{(a)}_0$ coincide with $b^{(a)}$ 
for all $1\leq a \leq r$ and their OPEs are given by 
\[ 
\partial_z b^{(a)}(z) \partial_w b^{(c)}(w) = -{\kappa^{-1}\over (z-w)^2} 
\] 
For any $h \in {\mathfrak h}^\vee$ let 
\[
h(z) = \sum_{a=1}^r \langle h, b_a\rangle b^{(a)}(z)
\]
where the angular brackets denote the canonical pairing ${\mathfrak h}^\vee \times {\mathfrak h}\to \IC(\kappa)$. Let $\pi_0$ be the Fock space of $b^{(a)}(z)$, $1\leq a\leq r$.

Now let $h^{(a)}$ be the fundamental weights of ${\mathfrak sl}_r$. The $W$-fields $W_d(z)\in {\rm End}(\pi_0)[[z^{-1}, z]]$, $0\leq d \leq r$ are defined by 
\[ 
-\kappa : \prod_{a=1}^r (Q\partial_z+h^{(a)}(z)):\ = \sum_{d=0}^r W_d(z)(Q\partial_z)^{r-d} 
\]
where $:\quad :$ indicates normal ordering
and 
\[
Q= -\kappa^{-1}\xi.
\]
This yields $W_0(z)=1$, $W_1(z)=0$, and some more complicated expressions for $W_d(z)$, $d\geq 2$. Abusing notation, for $d=1$ one sets 
\[ 
W_1(z) = \sum_{a=1}^r b^{(a)}(z)
\]
as opposed to $W_1(z)=0$. 
 Then the $W$-algebra $W_\kappa({\mathfrak gl}_r)$ is the vertex subalgebra of 
 $\pi_0$ generated by the Fourier modes of $W_d(z)$, $1\leq d \leq r$. 

Next note that any $\beta\in {\mathfrak h}$ determines a Verma module for the Heisenberg algebra ${\mathcal H}^{(r)}$ generated by the highest weight vector 
$|\beta \rangle$ satisfying 
\be\label{eq:HweightsA}
b^{(a)}_l |\beta\rangle = \delta_{l,0} \langle b^{(a)}, \beta \rangle |\beta\rangle, \qquad l \geq 0. 
\ee
Then there is a representation of $W_\kappa({\mathfrak gl}_r)$ on $\pi_\beta$ such that 
\be\label{eq:WweightsA}
W_{d,0}|\beta\rangle =w_d(\beta)|\beta\rangle, \qquad W_{d,l}|\beta\rangle =0, \quad l \geq 1, 
\ee
where 
\[ 
w_1(\beta) = \sum_{a=1}^r \langle b^{(a)}, \beta\rangle, \qquad 
w_d(\beta) = -\kappa \sum_{i_1< \cdots < i_d} \prod_{t=1}^d 
\big(\langle h^{(i_t)}, \beta\rangle + (d-t) \kappa^{-1} \xi \big) 
\]
In particular, for $\beta=0$ one obtains the vacuum $W$-module $\pi_0$.

\subsection{Module structure} 

Let ${\mathfrak U}(W_{\kappa}({\mathfrak gl}_r))$ denote the current algebra 
of the $W$-algebra, and let ${\mathcal U}_0(W_{\kappa}({\mathfrak gl}_r))$ denote its image in ${\rm End}(\pi_0)$
As shown in Corollary \ref{CYSHW}, which is analogous to Lemmas 8.22 and 8.24 in \cite{Cherednik_W}, there is an embedding of degreewise topological $K_0$-algebras
\be\label{eq:thetamapA}
\Theta_0^{(r)}: {\bf SH}^{(r)}_{K_0} \to {\mathcal U}_0(W_{\kappa}({\mathfrak gl}_r))
\ee
which lifts to a surjective morphism of degreewise topological $K_0$-algebras
\be\label{eq:thetamapB}
{\overline \Theta}_0^{(r)}: {\mathfrak U}( 
{\bf SH}^{(r)}_{K_0}) \to {\mathcal U}_0(W_{\kappa}({\mathfrak gl}_r)).
\ee
In particular, as observed in Corollary \ref{modequiv}, this yields  a representation 
\[
\pi_0^{(r)} : {\mathcal U}_0(W_{\kappa}({\mathfrak gl}_r)) \to 
{\rm End}\, ({\bf V}^{(r)}_{K_0}).
\]
As in \cite{Cherednik_W}, the morphism \eqref{eq:thetamapA} is obtained from a comparison the free field realizations of ${\bf SH}^{(r)}_{K_0}$ and respectively the $W$-algebra. The proof consists of a step by step ${\bf T}_0$-specialization of the proof given in \cite{Cherednik_W}. This is a straightforward, if somewhat tedious process, the details being provided in Appendix \ref{DDAHAW}.

Now one can conclude the proof of Theorem \ref{mainthmB} which states that 
$\pi_0^{(r)}$ is isomorphic to the vacuum representation of the $W$-algebra.

{\it Proof of Theorem \ref{mainthmB}.} Having proven the ${\bf T}_0$ variants  of Lemmas 8.33 and 8.34.a,  in \cite{Cherednik_W}, namely Lemmas \ref{SHmodA} and \ref{SHvaclemm}, 
the proof is now completely analogous to the proof of the first part of Theorem 8.32 in loc. cit. Using Corollary \ref{modequiv}, 
the above lemmas imply that ${\bf V}^{(r)}_{K_0}$ is Verma module for 
${\mathfrak U}_0(W_\kappa({\mathfrak gl}_r))$ with highest weight vector 
$[\alpha_{\emptyset^r}]$. Using the epimorhism \eqref{eq:thetamapB}, 
for each element $W_{d,0}$ of  ${\mathfrak U}(W_\kappa({\mathfrak gl}_r)$ there exists 
an element $W'_{d,0}$ in ${\mathfrak U}({\bf SH}^{(r)}_{K_0})$ mapped to $W_{d,0}$ in ${\mathfrak U}(W_\kappa({\mathfrak gl}_r)$.
Lemma \ref{freefieldvac} and Corollary \ref{Thetamap} 
imply that $W'_{d,0}$ acts in the same way on 
the vacua $[\alpha_{\emptyset^r}]$ and $[\alpha_{\emptyset}]^{\otimes r}$.
This  concludes the proof. 

\hfill $\Box$ 

\appendix

\section{From degenerate DAHA to $W$}\label{DDAHAW}

This main goal of this section is to show that the relation between degenerate DAHA modules and $W$-modules proven in \cite{Cherednik_W} also holds for 
their Calabi-Yau specializations. This includes
the construction of the free field representation of ${\bf SH}^{(r)}_{K_{r+2}}$ carried out in Section 8.5 
of \cite{Cherednik_W}, the algebra morphism obtained in Theorem 8.21 
of loc. cit. as well as the resulting categorical equivalence of admissible modules. The proof consists of a straightforward step-by-step verification that 
all intermediate steps in loc. cit. admit correct specialization under the inclusion ${\bf T}_0\subset {\bf T}_{r+2}$. The details are included here for completeness.

\subsection{Grading and order filtration}\label{gradingorder}

This is brief summary of Section 1.9 of \cite{Cherednik_W}. Using the notation of Section \ref{SHsect},  
let $D_{l,0}$, $l\geq \IZ$, be the elements of ${\bf SH}^{\bf c}$ 
defined recursively by
\[
[D_{1,1},D_{l,0}]=lD_{l+1,0}, \qquad [D_{-l,0},D_{-1,1}] = lD_{-l-1,0}, \qquad l\geq 0. 
\]
Moreover, for $l,k\geq 1$ set 
\[ 
D_{k,l}=[D_{0,l+1},D_{k,0}], \qquad D_{-k,l}=[D_{-k,0},D_{0,l+1}].
\]
As shown in Section 1.9 of \cite{Cherednik_W}, by construction 
there is an order filtration  
\[
\cdots \subset {\bf SH}^{\bf c}[\leq l] \subset {\bf SH}^{\bf c}[\leq l+1]
\subset \cdots 
{\bf SH}^{\bf c}
\]
where $l\geq 0$. 
 Proposition 1.38 in loc. cit shows that this filtration is completely determined by assigning the elements 
$D_{k,l}$, ${\bf c}_l$  order degrees $l$ and $0$ respectively.
Moreover, one has 
\[ 
{\bf SH}^{\bf c}[\leq l_1]\cdot {\bf SH}^{\bf c}[\leq l_2]\subseteq
{\bf SH}^{\bf c}[\leq l_1+l_2],
\]
hence the associated graded inherits an algebra structure.

In addition there is also a $\IZ$-grading such that $D_{l,0}$ has 
degree $l$ while $D_{0,l}$ has degree zero. This $\IZ$-grading is compatible with the above filtration.

\subsection{Coproduct}\label{coprodsect}

The topological tensor product ${\bf SH}^{\bf c} {\hat \otimes} 
{\bf SH}^{\bf c}$ is defined as 
\[
{\bf SH}^{\bf c} {\hat \otimes} 
{\bf SH}^{\bf c} = \bigoplus_{s\in \IZ} \varprojlim \bigg(
\bigoplus {\bf SH}^{\bf c}[s-t]\otimes {\bf SH}^{\bf c}[t]\bigg)/{\mathcal J}_N[s],
\]
\[
{\mathcal J}_N[s] = \bigoplus_{t\geq N} {\bf SH}^{\bf c}[s-t]\otimes {\bf SH}^{\bf c}[t].
\]
Then Theorem 7.9 in \cite{Cherednik_W} proves that there exists a 
coproduct $\Delta: {\bf SH}^{\bf c} \to {\bf SH}^{\bf c} {\hat \otimes}
{\bf SH}^{\bf c}$ 
which is uniquely determined by the formulas 
\be\label{eq:coprodA}
\bal 
\Delta({\bf c}_l)  & =\delta ({\bf c}_l), \qquad l \geq 0 \\
\Delta(D_{l,0}) & = \delta(D_{l,0}), \qquad l\neq 0\\
\Delta(D_{0,1}) & = \delta(D_{0,1}), \\ 
\Delta(D_{0,2}) & = \delta(D_{0,2}) + \xi \sum_{l\geq 1} 
l \kappa^{l-1} D_{-l,0} \otimes D_{l,0}, \\
\Delta(D_{1,1}) & = \delta(D_{1,1}) + \xi{\bf c}_0\otimes 
D_{1,0},\\
\Delta(D_{-1,1}) & = \delta(D_{-1,1}) + \xi
D_{-1,0}\otimes {\bf c}_0.\\
\eal
\ee
Here $\delta$ is the standard diagonal map. 

The main application of the coproduct resides in the construction of the free field representation 
\be\label{eq:freefieldA}
\rho^{(1^r)}: {\bf SH}^{(r)}_{K_{r+2}} \to 
{\rm End}\big({\bf L}^{(1)}_{K_{r+2}}\big)^{\otimes r}
\ee
in Section 8.5 of \cite{Cherednik_W}. 
As shown in loc. cit.,  Proposition 8.5, the coproduct determines 
naturally an injective algebra homomorphism 
\be\label{eq:coprodB} 
\Delta^{(1^r)} : {\bf SH}^{(r)}_{K_{r+2}} \to 
\bigg({\bf SH}^{(1)}_{K_{r+2}}\bigg)^{{\hat \otimes} r}.  
\ee
Moreover the case $r=1$ of Theorem 3.2 in loc. cit shows that there is a faithful representation 
\be\label{eq:rankonerepA}
{\bf SH}^{(1)}_{K_{r+2}}\to {\rm End}\, ({\bf L}^{(1)}_{K_{r+2}})
\ee
As observed in Corollary 8.7 of loc. cit, this yields a faithful representation \eqref{eq:freefieldA}.

The important point for the present purposes is the following: 
\begin{lemm}\label{CYfreefield} 
There is a faithful representation 
\be\label{eq:freefieldB}
\rho_0^{(1^r)}: {\bf SH}^{(r)}_{K_0} \to 
{\rm End}\, \big({\bf V}^{(1)}_{K_0}\big)^{\otimes r}.
\ee
\end{lemm} 

{\it Proof}. 
Clearly, using formulas \eqref{eq:coprodA} the construction of the injective algebra homomorphism \eqref{eq:coprodB} 
specializes immediately to ${\bf SH}^{(r)}_{K_0}$.  
Moreover, fact for $r=1$, the factor $({\IC}^{\times})^{\times r}$ acts trivially on $\CA(1,n)$ and 
the quotient ${\bf T}_{r+2}/ ({\IC}^{\times})^{\times r}$ is isomorphic to ${\bf T}_0$. Therefore the specialization of the rank $r=1$ case of Theorem 3.2 to ${\bf SH}^{(r)}_{K_0}$ is also immediate. This implies that the analogue of Corollary 8.7 in 
loc. cit also holds for the ${\bf T}_0$-specialization. Therefore 
Lemma \ref{CYfreefield} holds. 

\hfill $\Box$

To conclude this section, by analogy to Lemma 8.34.b in
\cite{Cherednik_W}, one has: 
\begin{lemm}\label{freefieldvac}
Let $[\alpha_{\emptyset}]^{\otimes r}\in ({\bf V}^{(1)}_{K_0})^{\otimes r}$ denote the vacuum vector.  Then 
$[\alpha_{\emptyset}]^{\otimes r}$ is annihilated by all endomorphisms of the form 
$\rho^{(1^r)}_0(D_{0,l})$, $\rho_0^{(1^r)}(D_{-l,0})$, $l\geq 1$. 
\end{lemm}

{\it Proof.} This follows from the $r=1$ case of Lemma \ref{SHvaclemm} using \cite[Lemma 7.11]{Cherednik_W}.

\hfill $\Box$

\subsection{Degreewise completion}\label{degcompl} 

As explained for example in Appendix A of \cite{Rep_th_W},  any $\IZ$-graded algebra as above has a natural degreewise linear topology defined by the decreasing sequence 
\[ 
{\mathcal J}_N = \bigoplus_{s\in \IZ} {\mathcal J}_N[s], \qquad 
{\mathcal J}_N[s] = \sum_{t\geq N} {\bf SH}^{(r)}_{K_{r+2}}[s-t] 
 {\bf SH}^{(r)}_{K_{{\bf T}_{r+2}}}[t].  
\]
This means that for any degree $s$ element $v\in {\bf SH}^{(r)}_{K_{r+2}}[s]$ the subsets $\{ v + {\mathcal J}_N[s]\}$ form a fundamental system of open neighborhoods of $v$. 

The standard degreewise completion of ${\bf SH}^{(r)}_{K_{r+2}}$ is defined by 
\be\label{eq:degcomplA}
{\mathfrak U}{\bf SH}^{(r)}_{K_{r+2}}=
\bigoplus_{s\in \IZ} {\mathfrak U}{\bf SH}^{(r)}_{K_{r+2}}[s], \qquad 
{\mathfrak U}{\bf SH}^{(r)}_{K_{r+2}}[s]=
\varprojlim {\bf SH}^{(r)}_{K_{r+2}}[s]/{\mathcal J}_N[s].
\ee

Clearly the grading, order filtration and degreewise completion remain well defined under  ${\bf T}_0$ specialization. 

Now, Definition 8.10 in \cite[Section 8.6]{Cherednik_W} introduces the notion 
of admissible ${\bf SH}^{(r)}_{K_{r+2}}$-module. This is a $\IZ$-graded module 
$M = \oplus_{s\in \IZ} M[s]$ 
such that $M[s]=0$ for sufficiently large $s$. This definition readily extends to modules over the degreewise completion \eqref{eq:degcomplA}. 
Then Proposition 8.11 in loc. cit. proves that 
 
$(1)$ The faithful representation \eqref{eq:freefieldA} extends to a  faithful representation of ${\mathfrak U}{\bf SH}^{(r)}_{K_{r+2}}$ on 
$\big({\bf L}^{(1)}_{K_0}\big)^{\otimes r}$. 

$(2)$ The canonical map ${\bf SH}^{(r)}_{K_{r+2}}\to {\mathfrak U}{\bf SH}^{(r)}_{K_{r+2}}$ is an embedding of degreewise topological algebras. 

Both statements are consequences on Corollary 8.7 in loc. cit, hence, in the view of Lemma \ref{CYfreefield}, they remain valid under ${\bf T}_0$-specialization. 

\subsection{From ${\bf SH}^{\bf c}$ to $W$}\label{SHtoWsect} 

Now let ${\mathfrak U}(W_\kappa({\mathfrak gl}_r))$ be the current algebra associated to the $W$-algebra. It was shown in \cite{Rep_th_W} that this is a graded degreewise complete topological $F$-algebra as defined in Section \ref{degcompl}. The grading is defined by the conformal degree.

Let ${\mathcal U}(W_\kappa({\mathfrak gl}_r))$ be the image of 
${\mathfrak U}(W_\kappa({\mathfrak gl}_r))$ in ${\rm End}(\pi_\beta)$ where $\beta$ is determined by the relations 
\be\label{eq:HweightsB} 
\langle b^{(a)}, \beta \rangle = -\kappa^{-1} \epsilon_a + (a-1)\kappa^{-1} \xi, \qquad 1\leq a \leq r. 
\ee
One of the main technical results in \cite{Cherednik_W}, Theorem 8.21, states that there is an embedding 
\[
\Theta^{(r)}: {\bf SH}^{(r)}_{K_{r+2}} \to {\mathcal U}(W_\kappa({\mathfrak gl}_r)) 
\]
of degreewise topological $K_{r+2}$-algebras with a degreewise dense image. 
Furthermore, Lemma 8.24 in loc. cit. proves that $\Theta^{(r)}$ extends to a surjective morphism of degreewise topological $K_{r+2}$-algebras 
\[
{\mathfrak U}({\bf SH}^{(r)}_{K_{r+2}}) \to {\mathcal U}(W_\kappa({\mathfrak gl}_r)). 
\]
Finally, Corollary 8.27 proves that the pull-back via $\Theta^{(r)}$ yields an equivalence from the category of admissible ${\mathcal U}(W_\kappa({\mathfrak gl}_r))$-modules to the category of admissible ${\bf SH}^{(r)}_{K_{{\bf T}_{r+2}}}$-modules. This equivalence intertwines between the free field realizations 
$\rho^{(1^r)}$ and $\pi_\beta$. 

The analogous statement for the ${\bf T}_0$-specialization is proven below.   First note that the Calabi-Yau specialization sets 
\[
\langle b^{(a)}, \beta \rangle =0, \qquad 1\leq a \leq r
\]
in equation \eqref{eq:HweightsB}. Therefore $\pi_\beta$ specializes to the Fock vacuum module $\pi_0$ of the Heisenberg algebra $\CH^{(r)}$. As above, let ${\mathcal U}_0(W_\kappa({\mathfrak gl}_r))$ denote the image of the $W$-algebra in ${\rm End}(\pi_0)$. 

\begin{lemm}\label{CYSHW} 
$(i)$ There is an embedding 
\be\label{eq:CYthetaA}
\Theta_0^{(r)}: {\bf SH}^{(r)}_{K_{0}} \to {\mathcal U}_0(W_\kappa({\mathfrak gl}_r)) 
\ee
of degreewise topological $K_{0}$-algebras with a degreewise dense image. 

$(ii)$ $\Theta_0^{(r)}$ extends to a surjective morphism of degreewise topological $K_{0}$-algebras 
\be\label{eq:CYthetaB}
{\mathfrak U}({\bf SH}^{(r)}_{K_{0}}) \to {\mathcal U}_0(W_\kappa({\mathfrak gl}_r)). 
\ee
\end{lemm} 

{\it Proof}. First recall the construction of the map 
 map $\Theta^{(r)}$, which is based on the free field realization 
\[
\rho^{(1^r)} : {\bf SH}^{(r)}_{K_{{r+2}}} \to {\rm End}({\bf L}^{(1)}_{K_{r+2}})^{\otimes r}. 
\]
Using the coproduct structure in Section \ref{coprodsect},
the morphism $\Theta^{(r)}$ is determined by $\Theta^{(1)}$ and $\Theta^{(2)}$. In order to summarize the explicit formulas for these maps, it will be helpful to recall that 
the algebra ${\bf SH}^{({\bf c})}$ is generated by the elements 
 ${\bf c}_l, D_{1,0}, D_{-1,0}, D_{0,2}$, as shown in 
\cite[Prop. 1.34]{Cherednik_W}, Moreover, it shown in Section 1.11 of loc. cit that the elements 
 \be\label{eq:HeisenbergA}
 b_l = (-x)^{-l}D_{-l,0}, \qquad b_{-l}=y^{-l}D_{l,0}, \qquad b_0= \kappa^{-1} E_1,
 \ee
 with $l\geq 1$, and ${\bf c}_0$ 
 define a Heisenberg subalgebra of ${\bf SH}^{\bf c}$. That is 
 \be\label{HeisenbergB}
 [b_l, b_k] = l \kappa^{-1} \delta_{l,k} {\bf c}_0.
 \ee
The analogous statements will hold for the ${\bf T}_{r+2}$ and ${\bf T}_0$ specializations. Then the construction of $\Theta^{(r)}$ proceeds as follows. 

For $r=1$ there is a single $W$-field 
\[ 
W_1=b(z) 
\]
and $\pi_\beta$ is a Verma module of the Heisenberg algebra $\CH^{(1)}$ with height weight 
\[
\beta = \kappa^{-1} \epsilon_1
\]
By Proposition 1.40 in \cite{Cherednik_W} there exists a unique isomorphism of $K_{{r+2}}$-vector spaces 
\[
{\bf L}_{K_{{r+2}}}^{(1)} \to \pi_\beta 
\]
mapping $[\alpha_\emptyset]$ to $|\beta\rangle$ which intertwines naturally between the action of Heisenberg subalgebra of 
${\bf SH}^{(1)}_{K_{{r+2}}}$ 
on ${\bf L}_K^{(1)}$
and 
the action of $\CH^{(1)}$ on $\pi_\beta$. 
As shown in Propositions 8.15 and  8.16 of loc. cit. this extends to an embedding $\Theta^{(1)}$ where 
\[ 
\Theta^{(1)}(b_l) = b_l,
\]
\[
\Theta^{(1)}(D_{0,1}) = \kappa\sum_{l\geq 1} b_{-l} b_l 
\]
\[
\Theta^{(1)}(D_{0,2}) = \kappa \Box - \epsilon_1 \Theta^{(1)}(D_{0,1}) 
\]
\[
\Box = \xi\sum_{l \geq 1} (l-1)b_{-l}b_l/2 + \kappa\sum_{l,k\geq 1} 
(b_{-l-k}b_lb_k+b_{-l}b_{-k}b_{l+k})/2. 
\]
The ${\bf T}_0$ specialization sets $\epsilon_1=0$, hence in this case $\pi_\beta = \pi_0$ is the Fock vacuum module of the Heisenberg algebra.  Clearly, the all the above formulas have well defined specialization. In particular, 
\[
\Theta^{(1)}_0(b_l) = b_l,
\]
\[
\Theta^{(1)}_0(D_{0,1})= \kappa\sum_{l\geq 1} b_{-l} b_l 
\]
and 
\[
\Theta^{(1)}_0(D_{0,2}) = \kappa\Box 
\]
define again an embedding of degreewise topological algebras. 

An analogous computation holds for $r=2$. In this case one identifies 
\[ 
\pi^{(1^2)} =\pi_\beta, \qquad \langle b^{(a)}, \beta \rangle = -\kappa^{-1} \epsilon_a+ (a-1)\kappa^{-1}\xi, \qquad 1\leq a\leq 2. 
\]
Note also that the $K_{{r+2}}$-vector space $\pi_\beta$ is canonically  isomorphic to 
$({\bf L}^{(1)}_{K_{{r+2}}})^{\otimes 2}$. Then $\Theta^{(2)}$ is determined in this case by the formulas 
\[ 
\Theta^{(2)}(b_l) = b_l, 
\]
\[
\bal
& \Theta^{(2)}(D_{0,2}) = \\
& {\kappa\over 2} 
\sum_{l\in \IZ} : W_{1,-l} W_{2,l}: + {\kappa^2\over 24} \sum_{k,l\in \IZ} 
: W_{1,-k-l}W_kW_l: + {\kappa\xi\over 4} \sum_{l\in \IZ} (|l|-1) : W_{1,-l}W_{1,l}: + \xi W_{0,2} + c.\\
\eal 
\]
where 
\[ 
c= p_3(\epsilon_1,\epsilon_2) + p_2(\epsilon_1,\epsilon_2) \xi/ 4\kappa - 
p_1(\epsilon_1, \epsilon_2) \xi^2/2\kappa + \xi^3 /12\kappa. 
\]
The ${\bf T}_0$ specialization sets
\[
\epsilon_1=0, \qquad \epsilon_2= \xi 
\]
and $\pi_\beta$ is again the Fock module $\pi_0$. Again, the above formulas well defined specializations. In fact the only $\epsilon$-dependence is through the constant $c$ which is polynomial in $\epsilon_1, \epsilon_2$. Therefore 
$\Theta^{(2)}_0$ is again well defined and yields an embedding of topological algebras. 

For $r\geq 2$, the map $\Theta^{(r)}$ is determined naturally by $\Theta^{(1)}$ and $\Theta^{(2)}$ using the coproduct \eqref{eq:coprodA} as shown in Theorem 8.21 of \cite{Cherednik_W}. Since the formulas \eqref{eq:coprodA} are independent of $\epsilon$, it follows that 
the ${\bf T}_0$-specialization of $\Theta^{(r)}$ is well defined and 
determines an embedding of topological algebras by analogy with 
Theorem 8.21 in \cite{Cherednik_W}.  
Note here that Theorem 8.23 used in loc. cit. also holds for $\pi^{\beta}=\pi_0$, as proven in Proposition 5.5 of \cite{W_coset}. 

The ${\bf T}_0$ specialization of Lemma 8.24 in \cite{Cherednik_W} also holds  since  the formulas proven in Claims 8.25 and 8.26 of loc. cit. are independent of $\epsilon$. 
Therefore the map $\Theta^{(r)}_0$ extends to a surjective morphism of 
degreewise topological algebras \eqref{eq:CYthetaB}.

\hfill $\Box$

Finally, note that the ${\bf T}_0$-specialization of Corollary 8.27 of \cite{Cherednik_W} also holds since the proof given in loc. cit. is completely independent of parameters. In conclusion one has: 

\begin{coro}\label{modequiv}
The  pull-back via $\Theta^{(r)}_0$ yields an equivalence from the category of admissible ${\mathcal U}_0(W_\kappa({\mathfrak gl}_r))$-modules to the category of admissible ${\bf SH}^{(r)}_{K_0}$-modules. This equivalence intertwines between the free field realizations 
$\rho_0^{(1^r)}$ and $\pi_0$. In particular, the vacuum vector $|0\rangle$ of $\pi_0$ is identified with the element $[\alpha_\emptyset]^{\otimes r}\in 
({\bf V}^{(1)}_{K_0})^{\otimes r}$. 
\end{coro} 

A second consequence of the proof of Lemma \ref{CYSHW} is recorded below.
\begin{coro}\label{Thetamap}
For any $d\geq 1$ let $W'_{d,0}$ in ${\mathfrak U}({\bf SH}^{(r)}_{K_0})$ be the element mapped by $\Theta_0$ to $W_{d,0}$ in ${\mathfrak U}(W_\kappa({\mathfrak gl}_r))$. 
Then $W'_{d,0}$ is a linear combination of monomials 
\[ 
D_{k_1,l_1} \cdots D_{k_r, l_r}, \qquad k_1+\cdots + k_r = 0
\]
where $D_{1,l}$, 
 $D_{0,l}$, $D_{-1,l}$, $l\geq 0$ appear exactly in this order from left to 
 right. 
 \end{coro} 
 
{\it Proof.} This is analogous to equation (8.123) in \cite{Cherednik_W}. 
It follows from the construction of the map $\Theta_0^{(r)}$ in Lemma \ref{CYSHW} using \cite[Proposition 8.3]{Cherednik_W} and Lemma \ref{SHgenlemm}.ii. 

\hfill $\Box$

\bibliography{Framed_Higgs_Ref}
 \bibliographystyle{halpha}

\bigskip
\end{document}